\documentclass{article}
\usepackage[utf8]{inputenc}
\usepackage{sty}
\usepackage{fullpage}
\usepackage{enumitem}
\usepackage{placeins}
\usepackage{subcaption}
\usepackage{graphicx}
\usepackage{caption}
\usepackage{color}
\usepackage{url}
\usepackage{relsize}
\usepackage{hyperref}

\title{Estimation under group actions: recovering orbits from invariants}

\usepackage{authblk}

\author[1]{Afonso S.\ Bandeira\thanks{email:~\texttt{bandeira@math.ethz.ch}. Most of this work was done while ASB was with the Department of Mathematics, Courant Institute of Mathematical Sciences and the Center for Data Science at NYU and supported partly by NSF grants DMS-1712730, DMS-1719545 and by a grant from the Sloan foundation. Part of this work was also done while this author was visiting the Simons Institute for the Theory of Computing; this visit was partially supported by the DIMACS/Simons Collaboration on Bridging Continuous and Discrete Optimization through NSF grant \#CCF-1740425.}}
\author[2]{Ben Blum-Smith\thanks{email:~\texttt{bblumsm1@jhu.edu}. Most of this work was done while BBS was affiliated with NYU.}}
\author[3]{Joe Kileel\thanks{email:~\texttt{jkileel@math.utexas.edu}.  Partially supported by the Simons Collaboration on Algorithms and Geometry.  Part of this work was done while with the Program in Applied and Computational Mathematics, Princeton University.}}
\author[4]{Jonathan Niles-Weed\thanks{email:~\texttt{jnw@cims.nyu.edu}. Parts of this work were completed while supported in part by NSF Graduate Research Fellowship DGE-1122374 and the Institute for Advanced Study.}}
\author[5]{\mbox{Amelia Perry}\thanks{Supported in part by NSF CAREER Award CCF-1453261 and a grant from the MIT NEC Corporation.}}
\author[6]{Alexander S.\ Wein\thanks{email:~\texttt{aswein@ucdavis.edu}. Part of this work was done while with the Courant Institute at NYU, partially supported by NSF grant DMS-1712730 and by the Simons Collaboration on Algorithms and Geometry. Part of this work was done while with the MIT Department of Mathematics. \\
\indent\textup{2020} \textit{Mathematics Subject Classification}. Primary 94A12, 94A16, 92C55, 62B10 Secondary 13A50, 13P25, 14Q20, 22E70}}
\affil[1]{Department of Mathematics, ETH Z\"urich}
\affil[2]{Department of Applied Mathematics and Statistics, Johns Hopkins University}
\affil[3]{Department of Mathematics and Oden Institute, University of Texas at Austin}
\affil[4]{Courant Institute of Mathematical Sciences and Center for Data Science, NYU}
\affil[5]{Department of Mathematics, Massachusetts Institute of Technology}
\affil[6]{Department of Mathematics, University of California, Davis}

\date{}

\begin{document}

\maketitle

\vspace{-5pt}

\begin{center}
In memory of Amelia Perry.
\end{center}

\vspace{15pt}

\begin{abstract}
We study a class of \emph{orbit recovery} problems in which we observe independent copies of an unknown element of $\RR^p$, each linearly acted upon by a random element of some group (such as $\ZZ/p$ or $\mathrm{SO}(3)$) and then corrupted by additive Gaussian noise. We prove matching upper and lower bounds on the number of samples required to approximately recover the group orbit of this unknown element with high probability. 
These bounds, based on quantitative techniques in  invariant theory, give a precise correspondence between the statistical difficulty of the estimation problem and algebraic properties of the group. Furthermore, we give computer-assisted procedures to certify these properties that are computationally efficient in many cases of interest.

The model is motivated by geometric problems in signal processing, computer vision, and structural biology, and applies to the reconstruction problem in cryo-electron microscopy (cryo-EM), a problem of significant practical interest. 
Our results allow us to verify (for a given problem size) that if cryo-EM images are corrupted by noise with variance $\sigma^2$, the number of images required to recover the molecule structure scales as $\sigma^6$.
We match this bound with a novel (albeit computationally expensive) algorithm for \textit{ab initio} reconstruction in cryo-EM, based on invariant features of degree at most 3. We further discuss how to recover multiple molecular structures from mixed (or heterogeneous) cryo-EM samples.
\end{abstract}

\section{Introduction}
Let $G$ be a compact Lie group acting continuously by orthogonal transformations on a real inner product space~$V$; for example, $G$ may be a finite group. 
Throughout, the space $V$ will be taken to be $\RR^p$, and by taking the quotient by the kernel of the action, we will regard the group $G$ as a subgroup of $\mathrm{O}(p)$, the orthogonal group.
Let $g$ be an $\mathrm{O}(p)$-valued random variable whose law is the normalized Haar measure on $G$, and let $\xi \sim \cN(0, \sigma^2 I)$ be a Gaussian vector independent of $g$.
Given a  vector $\theta \in V$, which we call the \emph{signal}, consider the \nolinebreak random \nolinebreak variable
\begin{equation}\label{eq:basic_sample}
    y \defeq g \cdot \theta + \xi.
\end{equation}
The law of $y$ depends only on the group orbit of $\theta$: if $\theta' \in V$ lies in the same $G$-orbit as $\theta$, then the right-invariance of the Haar measure implies that $y' \defeq g \cdot \theta' + \xi$ and $y$ have the same distribution.
However in all other cases, the laws of $y$ and $y'$ differ.
\begin{proposition} \label{prop:diff-dist}
If $\theta, \theta' \in V$ lie in different orbits of $G$, then $g \cdot \theta + \xi$ and $g \cdot \theta' + \xi$ have different distributions.
\end{proposition}
\begin{proof}
Suppose that $g \cdot \theta + \xi$ and $g \cdot \theta' + \xi$ have the same distribution.
Then their Fourier transforms
\begin{align*}
\EE \exp(i \langle (g \cdot \theta + \xi), t\rangle) & = \EE \exp(i \langle (g \cdot \theta), t\rangle) \exp( - \sigma^2 \|t\|^2/2) \\
\EE \exp(i \langle g \cdot \theta' + \xi), t\rangle) & = \EE \exp(i \langle (g \cdot \theta'), t \rangle) \exp( - \sigma^2 \|t\|^2/2)
\end{align*}
agree.
Since $\exp( - \sigma^2 \|t\|^2/2)$ is nonzero for all $t \in \RR^p$, this implies that $\EE \exp(i \langle (g \cdot \theta), t\rangle) = \EE \exp(i \langle (g \cdot \theta'), t\rangle)$ for all $t \in \RR^p$, and which implies that $g \cdot \theta$ and $g \cdot \theta'$ have the same law~\cite[Thm.~9.5.1]{Dud02}.
Since these measures are supported on the orbits of $\theta$ and $\theta'$, respectively, the orbits must agree.
\end{proof}

Thus given the distribution of $y$ in~\eqref{eq:basic_sample}, the orbit of $\theta$ is uniquely identified.
We call the problem of identifying the orbit of $\theta$ from $y$ the \emph{orbit recovery problem}.
This paper investigates a quantitative version of this problem: 
\begin{question}[informal]
Given $n$ independent and identically distributed (i.i.d.)\ realizations of~\eqref{eq:basic_sample}, how large must $n$ be so that this data approximately identifies the orbit of $\theta$ with high probability?
\end{question}
We will call the minimum possible value of $n$ the \emph{sample complexity}.
Our primary concern is to exhibit the optimal dependence of the sample complexity on the structure of the group $G$ and on the variance $\sigma^2$ of the Gaussian noise in the $\sigma \to \infty$ limit.

This question is motivated by many computational problems throughout the sciences, where measurements possess an algebraic structure.
An important example is cryo-electron microscopy (cryo-EM) \cite{AdrDubLep84,ss-cryo,nog-cryo}, an imaging technique in structural biology that was awarded the 2017 Nobel Prize in Chemistry. This technique seeks to estimate the structure of a large biological molecule, such as a protein, from many noisy 2-dimensional images of the molecule from random unknown directions in 3-dimensional space.
The experimenter therefore observes data corrupted by the action of $\mathrm{SO}(3)$, and estimating the 3-dimensional structure of the molecule amounts to approximately recovering the orbit of the molecule under this action.\footnote{As stated,~\eqref{eq:basic_sample} neglects two important practical aspects of cryo-EM: the fact that only 2-dimensional images of the rotated molecule are available, and the fact that chemically identical proteins in a single laboratory sample can appear in several different shapes, called ``conformations.'' We prove our results for a generalization of \eqref{eq:basic_sample}, introduced in Section~\ref{sec:gen-prob}, which includes both phenomena.}

The most extensively investigated prior example of the orbit recovery problem is what is known as \emph{multireference alignment}.
This is the setting where $G$ is the cyclic group $\mathbb{Z}/p$, acting on $V = \RR^p$ by cyclically permuting the coordinates.
Multireference alignment was proposed in the late 1990s~\cite{Sig98} as a simplified model of cryo-EM.
It was established~\cite{PerWeBan17} that, as $\sigma \to \infty$, the sample complexity is $n = \Theta(\sigma^6)$ for Lebesgue-almost every $\theta \in \mathbb R^p$.\footnote{We alert the reader to the fact that we will use $\mathrm{O}(p)$ to refer to the group of orthogonal matrices in dimension $p$ and will use $O(g(n))$ in the standard ``big-O'' sense: $f(n)=O(g(n))$ if there exists a constant $C$ such that $f(n)\leq C g(n)$ for all $n$ sufficiently large. In this case we also write $g(n) = \Omega(f(n))$. We use $f(n)=\Theta(g(n))$ to denote that both $f(n)=O(g(n))$ and $f(n)=\Omega(g(n))$. It will be clear from context which meaning of $O(\cdot)$ is meant.} This sample complexity rate is established by a two-part argument: \textbf{(a)}~showing that in the $\sigma\to\infty$ regime the sample complexity is tightly connected to the number of moments of the random variable $y$ needed to uniquely determine the orbit of $\theta$; and \textbf{(b)}~showing that for multi-reference alignment three moments are necessary and sufficient for Lebesgue-almost every $\theta \in \mathbb R^p$.

In this paper, we show that argument \textbf{(a)} is optimal for general orbit recovery problems, and show that this question is related to properties of the algebra of invariant polynomials (see Section~\ref{sec:stat}). More significantly, we leverage tools from algebraic geometry to create a streamlined process to answer question \textbf{(b)} in the general orbit recovery problem, replacing the ad hoc approach used for multi-reference alignment. Indeed, the problem of understanding when a subring of invariant polynomials resolves an orbit is a classical question in invariant theory (see for example~\cite{pop-vin, kac-notes, dolgachev, sturmfels, dk-book}). Leveraging these tools we relate the property of resolving orbits with classical algebraic notions and use this connection to create computational tests that compute the sample complexity of orbit recovery in general, these (efficient) computational tests are the main technical contribution of this paper. As we will describe below there are four natural notions of recovery (based on different combinations of choices for criteria \textbf{(i)} and \textbf{(ii)} in Section~\ref{sec:4recoveryconditions}). We catalogue in Table~\ref{table:ListJoe} numerous significant specific orbit recovery problems and their respective sample complexity obtained from our tools.

\subsection{Recovery criteria and main results}
\label{sec:4recoveryconditions}

It is natural to consider different success criteria for the orbit recovery problem. Consider the following two decisions:

\begin{enumerate}

\item[\textbf{(i)}] Do we study \emph{generic} signals $\theta$, or do we allow for a \emph{worst-case} signal? (Here \emph{generic} means that there can be a measure-zero set of disallowed signals.)

\item[\textbf{(ii)}] Do we ask to output a $\theta'$ such that $\theta'$ (approximately) lies in the orbit of $\theta$ (\emph{unique recovery}), or simply a finite list $\theta_1,\ldots,\theta_s$ of candidates such that one of them (approximately) lies in the orbit of $\theta$ (\emph{list recovery})?

\end{enumerate}

\noindent The terminology ``list recovery'' is borrowed from the idea of \emph{list decoding} in the theory of error-correcting codes~\cite{Eli57}. By taking all combinations of the two options above, there are four different recovery criteria. Strikingly, these different recovery criteria can be very different in terms of sample complexity, as the following examples show (see Section~\ref{sec:examples} for more details):

\begin{itemize}

\item \emph{Multi-reference alignment (MRA)}:
Recall the setup from the motivating example from the introduction.
It is known \cite{PerWeBan17} that if $\theta$ is generic then unique recovery is possible with $O(\sigma^6)$ samples. However, for a worst-case $\theta$, many more samples are required (even for list recovery); as shown in \cite{BanRigWee17}, there are some very particular infinite families of signals that cannot be distinguished without $\Omega(\sigma^{2p})$ samples, where $p$ is the length of the signal. This illustrates a large gap in difficulty between the generic and worst-case problems.

\item \emph{Learning a rigid body}: 
Let $G$ be the rotation group $\mathrm{SO}(p)$ acting on the matrix space $V = \RR^{p \times m}$ by left multiplication. (This problem, which we call \emph{learning a rigid body}, is further described in Section~\ref{sec:rigid-body}.)
With $O(\sigma^4)$ samples it is possible to recover the rigid body up to reflection, so that list recovery (with a list of size 2) is possible. However, unique recovery (even for a generic signal) requires drastically more samples: $\Omega(\sigma^{2p})$.

\end{itemize}

We address all four recovery criteria but our main focus is on the case of \emph{generic list recovery}, as it is algebraically the most tractable to analyze. For the following reasons we also argue that it is perhaps the most practically relevant case. First, it is reasonable to assume that real-world signals are generic, since almost all physical processes are subject to small random perturbations. Recovering a list of candidate signals is also a reasonable goal. In some practical applications, unique recovery is actually impossible; for instance, in cryo-EM it is impossible to determine the chirality of the molecule. (However, we can hope for unique recovery if we work over the group $\mathrm{O}(3)$ instead of $\mathrm{SO}(3)$.) Furthermore, one could hope to use application-specific clues to pick the true signal out from a finite list; for instance, in cryo-EM we might hope that the list contains only one, or very few, solutions in our finite list that look like ``reasonable'' molecules.

For the case of generic list recovery, we give a simple efficient algorithm for determining the degree of invariant polynomials required (which in turn dictates the optimal sample complexity) for any given orbit recovery problem. The main step of this algorithm is to compute the rank of a particular Jacobian matrix.

Furthermore, in certain cases, including a variant of cryo-EM without projections, we are able to extend our results to generic \emph{unique} recovery, where the output is a single solution instead of a finite list. We also show that if unique recovery is possible for some problem then unique recovery is also possible for the heterogeneous version of that problem, provided that a particular efficiently-testable criterion holds.

More precisely, for each of the four recovery conditions our main results are as follows.

\begin{itemize}
    \item \textbf{Generic list recovery}: we provide an exact algebraic characterization (in terms of transcendence degree) of the asymptotic order of the sample complexity, and an efficient algorithm to compute it (Theorems~\ref{thm:generic-list} and~\ref{thm:alg}). 
    \item \textbf{Generic unique recovery}: we provide a computationally inefficient algorithm, based on Gr\"obner bases, to compute an upper bound for the sample complexity based on a field generation criterion (Corollary~\ref{cor:generic-unique}).
    We give a method to determine the exact asymptotic order of the sample complexity, albeit at the cost of a cylindrical algebraic decomposition which is more involved (Theorem~\ref{thm:generic-unique}).
    \item \textbf{Worst-case unique recovery}: we provide a computationally inefficient algorithm, involving Gr\"obner bases and oracle access to symbolic semidefinite programming, to compute the exact order of the sample complexity (Proposition~\ref{prop:alg-separate}).
    \item \textbf{Worst-case list recovery}: we provide a computationally inefficient algorithm, based on Gr\"obner bases, to compute an upper bound for the sample complexity based on module-finiteness of a ring (Proposition~\ref{prop:alg-worst-list}).
\end{itemize}
For some specific orbit recovery problems, we improve upon the above generic methods and obtain computationally efficient tests to determine the exact order of the sample complexity (see Table~\ref{table:ListJoe}).   

We stress that all of the computationally efficient algorithms described above are methods \emph{to determine the sample complexity}, given as input the dimension of the vector space and a description of the group action. For most of the orbit recovery problems studied below, obtaining a computationally efficient algorithm to actually recover the signal from noisy samples is an open question. We indicate orbit recovery problems for which efficient recovery algorithms are known by check marks in the last column of Table~\ref{table:ListJoe}. These algorithms were known prior to our work, with the notable exception of unprojected cryo-EM, a.k.a. cryo-electron tomography.
For that orbit recovery problem, we give a new, efficient method to recover the signal itself from enough samples (Section~\ref{sec:cryoET}).

\subsection{Extensions to the orbit recovery problem}\label{sec:extensions}

We will also consider extensions to the basic orbit recovery problem defined above, motivated by cryo-EM:

\begin{enumerate}

\item {\bf Projection}: In cryo-EM, we do not observe a noisy 3-dimensional model of the rotated molecule; we only observe a 2-dimensional projection of it. We will model this projection by a linear map $\Pi: \RR^p \to \RR^q$ that maps a 3-dimensional model to its 2-dimensional projection (from a fixed viewing direction). The samples are then distributed as $y = \Pi(g \cdot \theta) + \xi$ where $\xi \sim \mathcal{N}(0,\sigma^2 I)$.

\item {\bf Heterogeneity}: In cryo-EM we observe images of many different copies of the same molecule, each rotated differently. However, if our sample is not pure, we may have a mixture of different molecules and want to recover the structure of all of them. We will model this by taking $K$ different unknown signals $\theta_1, \ldots, \theta_K$ along with positive mixing weights $w_1,\ldots,w_K$ which sum to 1. Each sample is distributed as $y = g \cdot \theta_k + \xi$ where $k$ is chosen at random according to the mixing weights.

\end{enumerate}

\noindent In Section~\ref{sec:gen-prob} we will formally define a generalization of the orbit recovery problem that allows for either (or both) of the above extensions. All of our methods can be adapted to handle this general case. 
For heterogenous problems, it is natural to ask whether the invariants of a mixture can be decomposed into invariants of the individual signals in a generically unique way, in which case the heterogeneous problem can be decoupled into multiple homogeneous problems. We call the property \emph{generic unique de-mixing} (Definition~\ref{def:demix}).

\begin{itemize}
    \item {\bf Generic unique de-mixing}: we provide an efficiently testable sufficient condition for generic unique de-mixing to hold. The algorithm involves computing the rank of a certain Hessian matrix (see Theorem~\ref{thm:hessian}).
\end{itemize}
We conjecture that this method proves that unique de-mixing holds for cryo-EM for regimes of interest, and we verified this for small problem sizes (see Section~\ref{subsec:het-cryo}).

\begin{table}
\begin{small}
\begin{tabular}{ |p{3.15cm}||p{3.7cm}|p{1.6cm}| p{2.3cm} | p{1.4cm} | p{1.5cm}|}
 \hline
 \multicolumn{6}{|c|}{ \textbf{\normalsize{Examples of Orbit Recovery Problems and (Conjectured) Sample Complexity  $\Theta\left(\sigma^{2d}\right)$}}} \\[1em]
 \hline
 \textbf{problem name} & \textbf{$G \curvearrowright V$}  & \textbf{nontrivial \mbox{projection}} & \textbf{\mbox{generic list} recovery} & \textbf{generic unique  \mbox{recovery}} & \textbf{efficient recovery known}\\[1em]
 \hline
 \mbox{bag of numbers} (Sec.~\ref{ex:bag})  & $S_p \curvearrowright \mathbb{R}^{p}$ (permutation)&   no & $d=p$ & $d=p$   & \hspace{0.5cm} \checkmark  \\[1em]
 \hline
 \mbox{rigid body} (Sec.~\ref{sec:rigid-body}) &  ${\mathrm{SO}(p) \curvearrowright \mathbb{R}^{p \times m}}$ \mbox{(left multiplication)}   & no & $d=2$  & $d=p$  & \hspace{0.5cm} \checkmark \\[1em]
 \hline
 \mbox{functions on a finite} group (Sec.~\ref{sec:reg_rep}) & ${G \curvearrowright \RR^G}$ \mbox{(regular representation)} &  no & $d = 3$ unless $G = (\ZZ / 2 \ZZ)^k$ & $d = 3$  & \hspace{0.5cm} \checkmark  \\[1em]
 \hline
 \mbox{MRA (Sec.~\ref{subsec:MRA})} \cite{PerWeBan17,BanRigWee17}  & $\mathbb{Z}/p \curvearrowright \mathbb{R}^p$ (cyclic shift)&  no & $d=3$ & $d=3$  & \hspace{0.5cm} \checkmark 
 \\[1em]
 \hline
 \mbox{projected MRA} (Sec.~\ref{subsubsec:MRA-proj}) &   $\mathbb{Z}/p \curvearrowright \mathbb{R}^p $ ($p$ odd)    & Eq.~\eqref{eq:mra-proj} & $d=3$ & & \\[1em]
 \hline
 \mbox{heterogenous MRA} (Sec.~\ref{sec:mra-het}) & $(\mathbb{Z}/p)^{K}  \wr  S_K  \curvearrowright (\mathbb{R}^p)^{\oplus K}$   & no& $d=3$ & ${d=3}$ \linebreak \mbox{($K \leq p/6$)} & \\[1em]
 \hline
  $S^2$-registration (Sec.~\ref{sec:ex-s2}) & ${\mathrm{SO}(3) \curvearrowright \oplus_{\ell=1}^F V_{\ell}}$  \mbox{(irreducible decomposition)} & no & ${d=3}$ ($F \geq 10$) & & \\[1em]
 \hline
 cryo-EM (Sec.~\ref{sec:ex-cryo}) & $\mathrm{SO}(3) \curvearrowright \left(\oplus_{\ell=1}^F V_{\ell}\right)^{\oplus S}$ & Eq.~\eqref{eq:cryo-projection} & ${d=3}$ ($S, F \geq 2$) & & \\[1em]
 \hline
 \mbox{unprojected cryo-EM} \linebreak (cryo-ET) (Sec.~\ref{sec:cryoET}) & $\mathrm{SO}(3) \curvearrowright \left(\oplus_{\ell=1}^F V_{\ell}\right)^{\oplus S}$& no & $d=3$ & ${d=3}$   \linebreak ($S \geq 3$) & \hspace{0.5cm} \checkmark   \\[1em]
 \hline
 \mbox{heterogeneous cryo-EM} (Sec.~\ref{subsec:het-cryo}) & ${\mathrm{SO}(3)^K \wr S_K\curvearrowright}$ $\left(\oplus_{\ell=1}^F V_{\ell}\right)^{\oplus S}\otimes \mathbb{R}^K$ & Eq.~\eqref{eq:cryo-projection} & ${d=3}$  \linebreak $(K \leq S^2/4)$ & & \\[1em]
 \hline
 \mbox{cryo-EM symmetric} molecules (Sec.~\ref{sec:sym-hom}) & ${\mathrm{SO}(3) \curvearrowright} \left(\left(\oplus_{\ell=1}^F V_{\ell}\right)^{\oplus S} \right)^{\mathbb{Z}/L}$ & Eq.~\eqref{eq:cryo-projection}  & ${d=2}$ ($F < L$) \linebreak ${d=3}$ ($F \geq L$)  & &\\[1em]
 \hline
 \mbox{heterogeneous cryo-EM} \mbox{symmetric molecules} (Sec.~\ref{sec:sym-het}) & ${\mathrm{SO}(3)^K \wr S_K \curvearrowright}$ ${\left(\left(\oplus_{\ell=1}^F V_{\ell}\right)^{\oplus S} \right)^{\mathbb{Z}/L}  \otimes \mathbb{R}^K}$ & Eq.~\eqref{eq:cryo-projection} & $d=3$  ($F > L$, up to critical $K$) & & \\[1em]
 \hline 
\end{tabular}
\end{small}
\caption{Examples of orbit recovery problems\label{table:ListJoe}. For most of these problems, the sample complexity $\Theta\left(\sigma^{2d}\right)$ described is conjectural and computer-verified for small cases using our methods. For the first four rows, the sample complexity results have been proved rigorously for all $p$. For some problems, the degree $d$ drops when $p \leq 2$. For some problems, it is also possible to efficiently solve the orbit recovery problem itself (i.e., recover the signal from noisy samples) as indicated in the last column.
}
\end{table}

\subsection{Relation to prior art in  Cryo-EM and other orbit recovery problems}

Several special cases of the orbit recovery problem have been studied for their theoretical and practical interest. Besides cryo-EM, another such problem is the first example problem mentioned above, \emph{multi-reference alignment} (MRA)~\cite{mra-sdp,AbbPerSin17,BanRigWee17,PerWeBan17}.
MRA is motivated by questions in signal processing~\cite{ZwaHeiGel03,PitZurAmo05} with further relevance to structural biology~\cite{Dia92,TheSte12}. 
Since the cyclic group $\ZZ/p$ is simpler than $\mathrm{SO}(3)$, understanding MRA has been seen as a useful stepping stone towards a full statistical analysis of cryo-EM.

Much of the work on orbit recovery problems has focused on a framework known as \emph{synchronization}.
This is an approach used in the analysis of cryo-EM images in practice~\cite{vainshtein-goncharov,vanheel,ss-cryo} based on comparing pairs of images and aligning (``synchronizing'') them based on shared features.
This alignment task is in general computationally infeasible; however, this approach has led to a theoretical analysis of methods to perform synchronization over various groups by means such as spectral methods \cite{singer-angular,ss-cryo}, semidefinite programming \cite{singer-angular,ss-cryo,mra-sdp,afonso-thesis,bandeira2020non}, and approximate message passing (AMP) \cite{pwbm-amp}.
A general Gaussian model for synchronization problems over any compact group is studied in \cite{pwbm-contig,pwbm-amp}. However, theoretical results on synchronization are largely confined to the setting where the variance $\sigma^2$ of the noise $\xi$ is sufficiently small.
Indeed, when $\sigma^2$ is large, reliably synchronizing the samples is impossible, even with unlimited computational power~\cite{AguDelBar16}. Moreover, it is the high-noise regime that is the practically relevant one for many applications, including cryo-EM, where the presence of large noise is a primary obstruction to current techniques~\cite{Sig16}. As a result, recent work has focused on analyzing~\eqref{eq:basic_sample} in the large-noise limit.

Compared to the low-noise regime, the situation in the high-noise is fundamentally different.
In fact, as we saw above, for multi-reference alignment when $\sigma \gg 1$, the sample complexity typically scales as $\sigma^6$ (first observed in the late 1990s~\cite{Sig98} and shown more recently~\cite{BanRigWee17,AbbPerSin17}) while, when $\sigma^2$ is smaller than some threshold, only $O(\sigma^2)$ samples are required.
From an algorithmic viewpoint, the connection with invariant polynomials suggests a class of recovery algorithms: use the data to estimate invariant polynomials and find the orbit that matches the values of those polynomials; since the variance of the noise of a polynomial of degree $d$ will be of order $\sigma^{2d}$, this approach succeeds at the optimal sample complexity scale. This is known as the \emph{invariant features} framework~\cite{AbbPerSin17,BanRigWee17,BenBouMa17,PerWeBan17,mra-het,cryo-2clean,spectral-bispectrum}. This approach has a long history in the signal processing literature~\cite{Kam80,Sad89,SadGia92} and is analogous to the well known ``method of moments'' in statistics~\cite{vdV98}. Note how the invariant features approach bypasses entirely the problem of estimating the group elements and focuses instead on estimating features of the signal which are preserved by the action of the group, and then recover the orbit from this information.

The application of invariant features to cryo-EM dates back to 1980 with the work of Kam \cite{Kam80}, who partially solved cryo-EM by means of degree-$2$ invariant features, reducing the unknown molecule structure to a collection of unknown orthogonal matrices. Kam proposed to use the degree-$3$ invariant features to determine the orthogonal matrices, in a heuristic procedure without algebraic uniqueness guarantees. 
Subsequent work has explored methods to estimate these orthogonal matrices \cite{cryo-orthogonal}, including recent work showing how two noiseless tomographic projections suffice to recover these orthogonal matrices \cite{cryo-2clean}. 
Our theory's contribution to cryo-EM can be viewed as rigorously solidifying Kam's idea to solve cryo-EM using the degree-$3$ invariants, while also circumventing the issue of missing orthogonal matrices, without requiring any noiseless projections. The approach is \emph{ab initio}, i.e., it does not require an initial guess of what the molecule looks like and thus cannot suffer from \emph{model bias}, which is a documented phenomenon~\cite{Cohen_Model_Bias} where the initial guess can have a significant effect on the result. \textit{Ab initio} estimates are particularly useful to serve as a model-free starting point for popular iterative refinement algorithms such as RELION~\cite{Scheres_Relion}.

We clarify that in this paper when we discuss computational procedures we are often describing a procedure to determine the optimal sample complexity (given a description of the group action as input), and not a recovery algorithm. Some of the procedures we develop of the former type are computationally efficient, but the corresponding procedure for actually solving an instance of the orbit recovery problem (recovering the orbit, e.g. a molecule, given samples as input) is not in general efficient because it involves solving a particular polynomial system whose size depends on the group action. There are fast non-convex heuristic methods to solve these systems in practice \cite{BenBouMa17,mra-het}, but we leave for future work the question of analyzing such methods rigorously and exploring whether they reach the information-theoretic limits determined in this paper. In some special cases, the polynomial system can be solved efficiently using tensor decomposition; this is true for MRA \cite{PerWeBan17} and certain other orbit recovery problems over finite groups (see Section~\ref{sec:jennrich}). Other efficient methods for MRA include frequency marching \cite{BenBouMa17}, which has also been applied to obtain a heuristic algorithm for cryo-EM~\cite{BarGre17}, and spectral decomposition of the bispectrum \cite{spectral-bispectrum}. Alternatively, homotopy continuation \cite{SW05} can be used to solve polynomial systems, where a proxy for  cost is the number of paths tracked during the calculation \cite{Lai16}.

\paragraph{Follow-up work.} Since the first version of this paper appeared, a number of works have made significant progress in both algorithmic and statistical aspects of the orbit recovery problem. 
These include works relating the sample complexity for the method of moments to that of maximum likelihood estimation by analyzing the non-convex likelihood landscape, using techniques from invariant theory introduced in the present 
paper \cite{fan2020likelihood,fan2021maximum,katsevich2020likelihood}.  In particular,  \cite{fan2021maximum} proved some algebraic conjectures we made based on results in small dimensions.
Another development has been the analysis of algorithms for signal recovery in unprojected cryo-EM, see \cite{LM-so3} which refines our analysis in Section~\ref{sec:cryoET}.
For a few other specific orbit recovery problems, new algorithms have been given for recovering the signal efficiently, for instance based on overcomplete tensor decomposition~\cite{alex-thesis} and tensor networks~\cite{tensor-net}.
There has also been work on high-dimensional multi-reference alignment, where the signal size grows with the noise variance \cite{doi:10.1137/20M1354994,DouFanZho22}. In contrast to our setting, the sample complexity exhibits a phase transition under this model.
Other work has imposed sparsity of the signal in cryo-EM, and shown that the sample complexity can improve  \cite{bendory2022sparsity}.

\subsection{Relation to invariant theory}

Theorems~\ref{thm:statistical-upper} and \ref{thm:statistical-lower} below show that the sample complexity of an orbit recovery task is controlled by the degrees of {\em polynomials invariant under the group action} that are needed to achieve a corresponding orbit-separation property. The study of polynomials invariant under a group action is the purview of invariant theory, an old and well-developed branch of algebra. 

The field's traditional goal is to describe explicitly the algebra of all invariant polynomial functions. Since the 19th century, culminating in the pioneering work of Hilbert \cite{hilbert1,hilbert2}, it has been known that the invariant algebra is finitely generated over the ground field in many cases of interest, so a finite description of the set of all invariant functions is possible. As algebraic geometry developed, this description task came to be seen as an algebraic approximation to the geometric task of constructing the action's orbit space (i.e., the space whose points parametrize the orbits of $G$ on $V$) as an algebraic variety. {\em Geometric Invariant Theory} \cite{GIT}, or {\em GIT}, an important research monograph originally published by David Mumford in 1965, worked out in detail to what extent and under what conditions the algebra of invariant polynomials accurately captures the geometry of the orbit space. Although we do not make use of the technical machinery of {\em GIT}, the present work is permeated by the underlying geometric point of view, which has become standard in the field. Below in Section~\ref{sec:it-ag-setup}, we explicate it in detail. Then, in the remainder of Section~\ref{sec:algebraic}, we use it to study the orbit recovery problem in general (with several proofs relegated to Section~\ref{sec:algebraic-proofs}), and in Section~\ref{sec:examples}, we work out many examples.

Most of the work in Section~\ref{sec:algebraic} is about providing algorithmic tests of the algebraic criteria we give for the four recovery conditions described in Section~\ref{sec:4recoveryconditions}, and discussing tools and considerations for applying them in practice to orbit recovery problems. (In the case of generic unique recovery, where the algebraic criterion is sufficient but not necessary, we also give an algorithmic test of the recovery condition itself; see Theorem~\ref{thm:generic-unique}.) The algebraic criteria themselves are all essentially direct corollaries of the standard geometric picture described in Section~\ref{sec:it-ag-setup}, though they vary in the amount of work involved: the characterization of generic list recovery (Theorem~\ref{thm:generic-list}) takes some care because we need to make sure that we have sufficient control over the $\RR$-points that the relevant standard geometric results over $\CC$ also work over $\RR$; on the other hand, the sufficient criteria for generic unique and worst-case list recovery (Corollary~\ref{cor:generic-unique} and Proposition~\ref{prop:worst-list}) are immediate inferences from the discussion in Section~\ref{sec:it-ag-setup}, and the characterization of worst-case unique recovery (Proposition~\ref{prop:worst-case-unique}) merely restates the definition of a separating set.

Although the primary aim of this paper is to understand the statistical properties of orbit recovery problems, and the application of invariant theory (and other algebraic techniques) is a means to that end, there are some results we hope will be of interest to invariant theorists and other algebraists. In particular:
\begin{itemize}
    \item Given an oracle to produce a strictly feasible point of a semidefinite program in exact arithmetic, we show that there is an algorithm to decide whether a given proposed set of invariants for the action of a compact Lie group is a separating set over $\RR$ (Proposition~\ref{prop:alg-separate}).
    \item We adapt previous algebro-geometric work on a test (based on computing the rank of a Hessian) to determine identifiability of the secant variety of a Veronese embedding of projective space, to a broader context. Specifically, we modify the test to apply to secant varieties of arbitrary varieties parametrized by maps that are not necessarily injective. We use this to relate the heterogeneous orbit recovery problem to the homogeneous problem (Theorem~\ref{thm:hessian}).
    \item We work out the degrees of polynomials needed to achieve various orbit recovery-related tasks for many examples of group actions (see Section~\ref{sec:examples}), revealing that there is typically a dramatic difference in the needed degree for generic vs.\ for worst-case tasks. Two cases of particular interest:
    \begin{itemize}
        \item For a large class of representations of $\mathrm{SO}(3)$, we show that, modulo a linear-algebra condition we have verified in the first fifteen cases, the polynomials of degree $\leq 3$ generically separate orbits (Theorem~\ref{thm:cryo-et}).
        \item For the regular representation of any finite abelian group, we give a Galois-theoretic proof that the field of invariant rational functions is generated by polynomials of degree $\leq 3$, regardless of the group (Theorem~\ref{thm:reg-rep-galois}).
    \end{itemize}
\end{itemize}

The results quoted in the last bullet connect this work with a long line of research in invariant theory, and also motivate some relatively unexplored questions. Toward the fundamental goal of describing invariant algebras, a much-studied problem has been to produce {\em a priori} bounds on the degrees of the generators. Such bounds are immediately relevant to the present setting, because if one needs to separate orbits using invariant polynomials, access to {\em all} the invariant polynomials is the best one can hope to do. Noether \cite{noether} gave a general bound for finite groups: over a field of characteristic zero, the invariant algebra is generated in degree at most the group order.\footnote{The restriction on the field characteristic has been lifted to the extent possible: the bound holds as long as the characteristic is prime to the group order \cite{fleischmann, fogarty}, and if this condition is not met, there is no uniform bound on the degrees of the invariants needed to generate the invariant algebra \cite{richman1996invariants}.} This bound has been lowered in important cases, e.g.,~\cite{gobel, sezer2002sharpening, cziszter-thesis, cziszter}. 

Still, for finite groups, the known bounds typically grow with the size of the group (for example, for cyclic groups, the Noether bound is sharp). For infinite groups, no bound depending only on $G$ is possible \cite{derksen2003global, bryant2005global}; known bounds grow with the degree of the representation \cite[Section~4.7]{dk-book}.\footnote{Again, better bounds exist in special cases, e.g., \cite{wehlau1993constructive}.} Thus, the control afforded by these bounds (in conjunction with Theorem~\ref{thm:statistical-upper}) 
over the sample complexity of orbit recovery problems is weak overall. The present application therefore motivates the search for degree bounds on sets of invariants satisfying weaker criteria than algebra generation, that still yield desirable separation properties. 

In particular, all of the algebraic criteria given below in Section~\ref{sec:algebraic} for the recovery conditions discussed in Section~\ref{sec:4recoveryconditions} are weaker than algebra generation. 
Thus the present work motivates degree bounds for these criteria.

For the worst-case recovery conditions, this is already an area of active research. In the original 2002 edition of \cite{dk-book}, Derksen and Kemper introduced the notion of a separating set---a subset of the invariant ring that separates all orbits of the group action which are separated by the full invariant algebra. A separating set is thus as effective as an algebra generating set at distinguishing orbits. In the present case of a compact group, the full invariant ring separates all orbits, so a separating set achieves worst-case unique recovery; see Section~\ref{sec:worst-unique} for more details. The question of degree bounds for separating sets has begun to be studied \cite{kohlskraft, kemper-sep, domokos} (and see also \cite{kadish, cahill2022group}, which undertake related inquiries for non-polynomial functions). In \cite{domokos} it is shown that, for finite abelian groups, separating sets can usually be found in lower degree than a full generating set. 

Meanwhile, for worst-case list recovery, the study of degree bounds for the corresponding algebraic condition on a set of invariants (namely, finiteness of the invariant algebra over the subalgebra it generates; see Section~\ref{sec:worst-unique} for details) reaches further back. It was originally studied because it is related to the problem of bounding the degrees of algebra generators \cite[Section~4.7]{dk-book}, and results in this direction go back to Hilbert \cite{hilbert2}. General degree bounds were given by Popov \cite{popov1981constructive, popov1982}, Hiss \cite{hiss}, and Derksen \cite{derksen2001polynomial}; see also \cite[Propositions~4.7.12 and 4.7.16]{dk-book}. This question continues to be a topic of contemporary research \cite{cziszter2013generalized, elmer2014zero, elmer2016zero}.\footnote{In these latter papers, the question is studied in an equivalent (see \cite[Lemma~2.5.5]{dk-book}) form: what degree invariants are needed to separate zero from other orbits (to the extent possible)?}

On the other hand, degree bounds for the algebraic criteria associated with the generic recovery conditions appear to have received comparatively little attention. Generic unique recovery is achieved by a set of invariants that generates the {\em field} of rational invariants (see Section~\ref{sec:generic-unique} for details). Tight degree bounds for generators of invariant fields do not appear to have been systematically studied in invariant theory. It has been an active research direction to develop algorithms for finding field generators, based on Gr\"{o}bner bases \cite{MQB, hubert-kogan, kemper-field}. However, these works do not attend to the degrees of the generators they produce (which in any case are not necessarily polynomials). Hubert and her collaborators have given alternative, Gr\"{o}bner-basis-free methods in several important special cases \cite{hubert-labahn1, hubert-labahn2, gor-hub-pap}, and Fleischmann, Kemper, and Woodcock have given a general explicit construction for finite groups \cite{FKW}. While \cite{hubert-labahn2} and \cite{FKW} do extract degree bounds from their constructions---both variants on the Noether bound---it is not a goal of any of these works to minimize the degrees of the generators. The results quoted in the last bullet above (and proven in Section~\ref{sec:examples}) seem to indicate that there is much more to learn.

Meanwhile, generic list recovery is achieved by a transcendence basis (see Section~\ref{sec:generic-list} for details). The work mentioned above on the algebraic condition associated with worst-case list recovery provides {\em a priori} upper bounds for transcendence bases as well, since generic list recovery is a relaxation of worst-case list recovery---but the point would be to improve on these bounds. We are not aware of work in invariant theory that pursues this goal, although we do not claim an exhaustive review. Potentially relevant results are found in work on constructive approaches to identifying when the invariant field of a finite group is purely transcendental \cite{kemper1996constructive}.

We hope that the present work stimulates further investigation of both of these questions.

\subsubsection*{Organization of the paper}\label{sec:org}

In Section~\ref{sec:prob-statement}, we define the orbit recovery problem and specify the basic algebraic objects that will be central to our results. In Section~\ref{sec:gen-prob} we define a generalization of the problem that captures projection and heterogeneity. In Section~\ref{sec:stat}, we formally state statistical upper and lower bounds which connect the optimal sample complexity to certain algebraic properties of the group action. In Section~\ref{sec:algebraic}, we establish our basic algebraic results and specify the algebraic properties that correspond to the different recovery criteria discussed earlier. We show how to efficiently test some of these criteria, namely the criterion for generic list recovery and the criterion for generic unique de-mixing. In Section~\ref{sec:examples}, we apply our results to several examples of the orbit recovery problem, giving concrete results for various problems including MRA, cryo-EM and various extensions. This includes proofs of generic unique recovery in a few special cases of the orbit recovery problem: the regular representation of a finite group (Section~\ref{sec:reg_rep}), and unprojected cryo-EM (Section~\ref{sec:cryoET}). We conclude in Section~\ref{sec:open-questions} with questions for future work.

Sections~\ref{sec:proofs-stat} and~\ref{sec:algebraic-proofs} contain proofs of results from preceding sections. Appendix~\ref{app:algprimer} contains a primer on algebra and invariant theory. Appendix~\ref{app:so3} contains an account of the invariant theory of $\mathrm{SO}(3)$.

\section{Orbit recovery and invariant features}

\subsection{Basic problem statement}\label{sec:prob-statement}

For ease of notation, we begin with the problem statement of the basic orbit recovery problem without the projection and heterogeneity extensions. These extensions are captured by the general problem statement that we give later in Section~\ref{sec:gen-prob}.

Recall the following setup from the introduction. Throughout, we consider a compact Lie group $G$ acting continuously and orthogonally on a finite-dimensional real vector space $V = \RR^p$. In other words, $G$ acts on $V$ via a linear representation $\rho: G \to \mathrm{O}(V)$, and the action map $(g,v)\mapsto \rho(g)(v)$ is continuous; equivalently, $\rho$ itself is continuous. Here $\mathrm{O}(V)$ denotes the space of real orthogonal $p \times p$ matrices. In what follows, $G$ only enters through its action on $V$. Thus, we may take it to be a subgroup of $\mathrm{O}(V)$.

Let $\Haar(G)$ denote the normalized Haar measure on $G$, the unique invariant Radon probability measure on $G$.\footnote{See \cite[Section~1.4.2]{sepanski} for a systematic development of the Haar measure on a compact Lie group.} We define the \emph{orbit recovery} problem as follows.

\begin{problem}[orbit recovery]
\label{prob:orbit}
Let $V = \RR^p$ and let $\theta \in V$ be the unknown signal. 
Let $G$ be a compact group that acts linearly, continuously, and orthogonally on $V$.
For $i \in [n] = \{1,2,\ldots,n\}$ we observe
$$y_i = g_i \cdot \theta + \xi_i$$
where $g_i \sim \Haar(G)$ and $\xi_i \sim \cN(0,\sigma^2 I_{p \times p})$, all independently. The goal is to estimate the orbit of $\theta$.
\end{problem}

Our primary goal is to study the sample complexity of the problem: how must the number of samples $n$ scale with the noise level $\sigma$ (as $\sigma \to \infty$ with $G$ and $V$ fixed) in order for orbit recovery to be statistically possible?

In practical applications, $\sigma$ is often known in advance and, when it is not, it can generally be estimated accurately on the basis of the samples. We therefore assume throughout that $\sigma$ is known and do not pursue the question of its estimation in this work.

The objective of Problem~\ref{prob:orbit} is to estimate $\theta$ given the samples $\{y_i\}$. Note, however, that we can only hope to recover $\theta$ up to action by $G$; thus we aim to (approximately) recover the \emph{orbit} of $\theta$, defined as follows.

\begin{definition}
For $\theta_1, \theta_2 \in V$, define an equivalence relation $\stackrel{G}{\sim}$ by writing $\theta_1 \stackrel{G}{\sim} \theta_2$ if there exists $g \in G$ such that $g \cdot \theta_1 = \theta_2$. The \emph{orbit} of $\theta$ (under the action of $G$) is the equivalence class of $\theta$ under $\stackrel{G}{\sim}$, i.e.\ the set $\{g \cdot \theta \;:\; g \in G\}$. Denote by $V/G$ the set of orbits of $V$, that is, the equivalence classes of $V$ modulo the relation $\stackrel{G}{\sim}$.
\end{definition}

The following definition captures the notion of \emph{approximately} recovering the orbit of $\theta$.

\begin{definition}
For $\theta_1, \theta_2 \in V$, let
$$d_G(\theta_1,\theta_2) = \min_{g \in G} \|\theta_1 - g \cdot \theta_2\|_2.$$
This pseudometric induces a metric on the quotient space $V/G$ in the obvious way, so we can write  $d_G(\mathfrak{o}_1,\mathfrak{o}_2)$ for $\mathfrak{o}_1,\mathfrak{o}_2 \in V/G$.
By slight abuse of notation, we write $d_G(\theta_1, \mathfrak{o}_2)$ for $d_G(\mathfrak{o}_1, \mathfrak{o}_2)$, where $\mathfrak{o}_1$ is the orbit of $\theta_1$.
\end{definition}

Our techniques rely on estimation of the following moments:

\begin{definition}[moment tensor]\label{def:mom-ten}
For an integer $d \ge 1$, the \emph{order-$d$ moment tensor} is $T_d(\theta) \defeq \EE_g[(g \cdot \theta)^{\otimes d}]$ where $g \sim \Haar(G)$.\footnote{This definition is specific to the present setting, in which $V$ is a real vector space on which $G$ acts orthogonally; the order-$d$ moment tensor for a compact group acting unitarily on a complex vector space is $\EE_g[(g\cdot\theta)^{\otimes d-1}\otimes \overline{g\cdot \theta}]$.}
\end{definition}

\noindent We can estimate $T_d(\theta)$ from the samples by computing $\frac{1}{n} \sum_{i=1}^n y_i^{\otimes d}$ plus a correction term to cancel bias from the noise terms (see Section~\ref{sec:proofs-stat} for details). The moments $T_d(\theta)$ are related to polynomials that are invariant under the group action, which brings us to the fundamental object in invariant theory:

\begin{definition}[invariant ring] \label{def:inv-ring}
Let ${\bf x} = (x_1,\ldots,x_p)$ be a set of coordinate functions on $V = \RR^p$, i.e., a basis for the dual $V^*$. (If we fix a basis for $V$, we can think of ${\bf x}$ as indeterminate variables corresponding to the entries of $\theta \in V$.) Then $\Rx \defeq \RR[x_1,\ldots,x_p]$ is the ring of polynomial functions $V \to \RR$.\footnote{Here and throughout, we use the standard convention that the notation $\RR[\cdot]$ indicates a ring of polynomial functions with coefficients in $\RR$, while $\RR(\cdot)$ indicates a field of rational functions with coefficients in $\RR$; similarly for $\CC[\cdot]$ and $\CC(\cdot)$, although complex coefficients will be used only very sparingly.} We have an action of $G$ on $\Rx$ given by $(g \cdot f)(\cdot) = f(g^{-1}(\cdot))$.
The \emph{invariant ring} $\RxG \subseteq \Rx$ is the ring consisting of polynomials $f$ that satisfy $g \cdot f = f$ for all $g \in G$. An element of the invariant ring is called an \emph{invariant polynomial} (or simply an \emph{invariant}).
Invariant polynomials can be equivalently characterized as polynomials of the form $\EE_g[g \cdot f]$ where $f \in \Rx$ is any polynomial and $g \sim \Haar(G)$.
\end{definition}

\noindent The two objects above are equivalent in the following sense. The moment tensor $T_d(\theta)$ contains the same information as the set of evaluations $f(\theta)$ for all $f \in \RxG$ that are homogeneous of degree $d$. In particular, for any such polynomial $f$, $f(\theta)$ is a linear combination of the entries of $T_d(\theta)$, and every entry of $T_d(\theta)$ is such a polynomial.

The main observation motivating our work is the following connection between the sample complexity of orbit recovery and the algebraic structure of the invariant polynomials.
Concretely, we show that the answer to the statistical question is completely characterized by the answer to the following algebraic question:
given a specified set of invariant polynomials, when do the values of these polynomials, evaluated on $\theta$, determine the orbit of $\theta$ (in the appropriate sense)?
The following informal version of our main statistical results makes this connection precise.

\begin{theorem-nolabel}[informal version of Theorems \ref{thm:statistical-upper} and \ref{thm:statistical-lower}]\label{thm:recovery-upper}
Fix $\theta \in V$ and $\varepsilon > 0$. Let $d^*$ be the smallest $d$ for which the moments $T_1(\theta), \cdots, T_d(\theta)$ uniquely determine the orbit of $\theta$. Then $\Theta(\sigma^{2d^*})$ samples are necessary and sufficient in order to produce an estimator $\hat \theta$ such that $d_G(\theta,\hat\theta) \le \eps$ with high probability.
\end{theorem-nolabel}

\noindent Here, the asymptotic notation $\Theta(\cdot)$ pertains to the high-noise limit $\sigma \to \infty$ with all other parameters held fixed. In other words, $\Theta(\cdot)$ suppresses constants depending on $G$ (and its action), $V$, $\theta$, and $\varepsilon$.

The analogous result holds for list recovery (see Section~\ref{sec:stat}): if the moments determine a finite number of possibilities for the orbit of $\theta$ then we can output a finite list of estimators, one of which is close to the orbit of $\theta$.

Thus, we have reduced our problem to the algebraic question of determining how many moments are necessary to determine the orbit of $\theta$ (either uniquely or in the sense of list recovery). In Section~\ref{sec:algebraic} we will use tools from invariant theory and algebraic geometry in order to address these questions.

\subsection{General problem statement}
\label{sec:gen-prob}

Our results will apply not only to the basic orbit recovery problem (Problem~\ref{prob:orbit}) but to a generalization (Problem~\ref{prob:gen-orbit} below) that captures the projection and heterogeneity extensions discussed in the introduction. We first define mixing weights for heterogeneous problems.

\begin{definition}[mixing weights]\label{def:mixing-weights}
Let $w = (w_1,\ldots,w_K) \in \Delta_K \defeq \{(z_1,\ldots,z_K) \,:\, z_k \ge 0 \;\;\forall k,\, \sum_{k=1}^K z_k = 1\}$. Let $k \overset{w}{\sim} [K]$ indicate that $k$ is sampled from $[K] = \{1,\ldots,K\}$ such that $k = \ell$ with probability $w_\ell$. We will sometimes instead parametrize the mixing weights by $\bar w_k = w_k - 1/K$ so that $\bar w$ lies in the vector space $\bar \Delta \defeq \{(z_1,\ldots,z_K) \;:\; \sum_{k=1}^K z_k = 0\}$.
\end{definition}

\noindent In a heterogeneous problem with $K$ different signals, we can only hope to recover the signals up to permutation. To formalize this, our compound signal will lie in a larger vector space $V$ and we will seek to recover its orbit under a larger group $G$.

\begin{definition}[setup for heterogeneity] \label{def:het-action}
Let $\tilde G$ be a compact Lie group acting linearly, continuously, and orthogonally on $\tilde V = \RR^p$.
Let $V = \tilde V^{\oplus K} \oplus \bar\Delta_K$, so that $\theta \in V$ encodes $K$ different signals along with mixing weights: $\theta = (\theta_1,\ldots,\theta_K,\bar w)$. We let an element $(g_1,\ldots,g_K,\pi)$ of the Cartesian product set $\tilde G^K \times S_K$ act on $V$ as follows: first, each $g_k$ acts on the corresponding $\theta_k$, and then $\pi$ permutes the $\theta_k$ and the coordinates of $\bar w$. Note that this action is linear and orthogonal (where $\bar \Delta$ uses the usual inner product inherited from $\RR^K$). There is a natural group structure $G$ on the set $\tilde G^K \times S_K$ such that the action just described is actually a group action by $G$: the semidirect product $G = \tilde G^K \rtimes_\varphi S_K$, where $\varphi$ denotes the action of $S_K$ on $\tilde G^K$ by permutations of the factors. This is also called the {\em wreath product} of $\tilde G$ by $S_K$ and written $\tilde G\wr S_K$. Like $G$ it is a compact Lie group (the topology is the product topology on $\tilde G^K \times S_K$, with $S_K$ having the discrete topology), and the action described above is continuous.
\end{definition}

\noindent By taking $K=1$ we recover the basic setup (without heterogeneity) as a special case. We are now ready to give the general problem statement.

\begin{problem}[generalized orbit recovery]
\label{prob:gen-orbit}
Let $\tilde V = \RR^p$ and $W = \RR^q$.
Let $\tilde G$ be a compact Lie group acting linearly, continuously, and orthogonally on $\tilde V$.
Let $\Pi: \tilde V \to W$ be a linear map. Let $\theta = (\theta_1,\ldots,\theta_K,\bar w) \in V \defeq \tilde V^{\oplus K} \oplus \bar\Delta_K$ be an unknown collection of $K$ signals with mixing weights $w \in \Delta_K$. For $i \in [n] = \{1,2,\ldots,n\}$ we observe
$$y_i = \Pi(g_i \cdot \theta_{k_i}) + \xi_i$$
where $g_i \sim \Haar(\tilde G)$, $k_i \overset{w}{\sim} [K]$, $\xi_i \sim \cN(0,\sigma^2 I_{q \times q})$, all independently. The goal is to estimate the orbit of $\theta$ under $G \defeq \tilde G^K \rtimes S_K$.
\end{problem}

\noindent Note that as in the basic setup, our goal here is to recover the orbit of a vector $\theta$ under the action of some compact group. In this case we do not necessarily have an analogue to Proposition~\ref{prop:diff-dist}: depending on the nature of the projection $\Pi$, the law of $y_i$ may or may not contain enough information to determine the orbit of $\theta$.

We now provide analogues of Definitions~\ref{def:mom-ten} and~\ref{def:inv-ring}.
The moments are defined as follows.

\begin{definition}[moment tensor] \label{def:mom-ten-gen}
For the generalized orbit recovery problem (Problem~\ref{prob:gen-orbit}), the \emph{order-$d$ moment tensor} is $T_d(\theta) \defeq \EE_{g,k}[(\Pi(g \cdot \theta_k))^{\otimes d}]$ where $g \sim \Haar(\tilde G)$ and $k \overset{w}{\sim} [K]$. Equivalently, $T_d(\theta) = \sum_{k=1}^K w_k \,\EE_g[(\Pi(g \cdot \theta_k))^{\otimes d}]$.
\end{definition}

\noindent The invariant ring is defined as in Definition~\ref{def:inv-ring} but now for the larger group $G$ acting on the larger $V$:

\begin{definition}[invariant ring]
Note that $\dim(V) = Kp + K-1$ and let ${\bf x} = (x_1,\ldots,x_{\dim(V)})$ be a basis for $V^*$; here the last $K-1$ variables correspond to $\bar\Delta$, e.g., they can correspond to $\bar w_1,\ldots,\bar w_{K-1}$. We then let $\RxG \subseteq \Rx$ be the polynomials in ${\bf x}$ that are invariant under the action of $G$ (as in Definition~\ref{def:inv-ring}).
\end{definition}

Recall that in the basic orbit recovery problem, $T_d(\theta)$ corresponds precisely to the homogeneous invariant polynomials of degree $d$; now $T_d(\theta)$ corresponds to a subspace of the homogeneous invariant polynomials of degree $d$. Specifically, the method of moments gives us access to the following polynomials (evaluated at $\theta$):

\begin{definition}\label{def:U}
Let $U^T_d$ be the subspace (over $\RR$) of the invariant ring $\RxG$ consisting of all $\RR$-linear combinations of entries of $T_d({\bf x})$. Let $U^T_{\le d} = U^T_1 \oplus \cdots \oplus U^T_d \subseteq \RxG$. Here we write $T_d({\bf x})$ for the collection of polynomials (one for each entry of $T_d(\theta)$) that map $\theta$ to $T_d(\theta)$.
\end{definition}

We will be interested in whether the subspace $U^T_{\le d}$ contains enough information to uniquely determine the orbit of $\theta$ (or determine a finite list of possible orbits) in the following sense.

\begin{definition}
A subspace $U \subseteq \RxG$ \emph{resolves} $\theta \in V$ if there exists a unique $\mathfrak{o} \in V/G$ such that $f(\theta) = f(\mathfrak{o})$ for all $f \in U$.
Similarly, $U$ \emph{list-resolves} $\theta$ if there are only finitely many orbits $\mathfrak{o}_1,\ldots,\mathfrak{o}_s$ such that $f(\theta) = f(\mathfrak{o}_i)$ for all $f \in U$.
\end{definition}

\noindent Here we have abused notation by writing $f(\mathfrak{o})$ to denote the (constant) value that $f$ takes on every $\theta \in \mathfrak{o}$. The following question is of central importance.

\begin{question}\label{que:alg}
Fix $\theta \in V$. How large must $d$ be in order for $U^T_{\le d}$ to uniquely resolve $\theta$? How large must $d$ be in order for $U^T_{\le d}$ to list-resolve $\theta$?
\end{question}

The answer depends on $G$ and $V$ but also on whether $\theta$ is a generic or worst-case signal, and whether we ask for unique recovery or list recovery. (Note that in the heterogeneous situation, $\theta$ includes the mixing weights as well as the signals; thus a generic $\theta$ means a generic $K$-tuple of signals with generic mixing weights.) Our statistical results in Section~\ref{sec:stat} will show that the sample complexity of the generalized orbit recovery problem is $\Theta(\sigma^{2d})$ where $d$ is the minimal $d$ from Question~\ref{que:alg}. More specifically, the recovery procedure that obtains this bound is based on estimating the moments $T_1(\theta), \ldots, T_d(\theta)$ and solving a system of polynomial equations to (approximately) recover $\theta$. Our algebraic results in Section~\ref{sec:algebraic} will give general methods to answer Question~\ref{que:alg} for any $G$ and $V$.

\subsection{Statistical results}
\label{sec:stat}

In this section, we state upper and lower bounds on the performance of optimal estimators for the orbit recovery problem. Proofs are deferred to Section~\ref{sec:proofs-stat}. Our approach will be based on the classical \emph{method of moments}, which is the basis for the link between the orbit recovery problem and invariant theory explored in the remainder of the paper.

We assume for normalization purposes that there exists a constant $c \geq 1$ such that $c^{-1} \leq \|\theta\| \leq c$, so that $\sigma$ captures entirely the signal-to-noise ratio of the problem. We denote by $\Theta$ the subset of $V$ consisting of vectors satisfying this requirement.

Denote by $\mathrm{P}_{\theta}$ the distribution of a sample arising from the generalized orbit recovery problem (Problem~\ref{prob:gen-orbit}) with parameter $\theta$.

\begin{definition}\label{def:order-d-set}
Given $\theta \in \Theta$, the \emph{order-$d$ matching set for $\theta$}, $\Mthetad$, is the set consisting of all $\tau \in V$ such that $f(\tau) = f(\theta)$ for all $f \in U^T_{\leq d}$.
\end{definition}

We note that $U^T_{\leq d}$ resolves $\theta$ exactly when $\Mthetad$ contains a single orbit, and $U^T_{\leq d}$ list-resolves $\theta$ when $\Mthetad$ is the union of a finite number of orbits.

We are now ready to state formal theorems justifying the informal statement in Section~\ref{sec:prob-statement}, above.
The following theorem establishes that we can approximately learn the order-$d$ matching set for $\theta$ with probability at least $1-\delta$ on the basis of $O(\sigma^{2d}\log(1/\delta))$ samples.

\begin{theorem}\label{thm:statistical-upper}
Suppose that $\Mthetad$ is the union of $M$ orbits $\mathfrak{o}_1, \dots, \mathfrak{o}_M$, where $M$ is finite.
There exists an $\eps_\theta$ such that, for $\eps < \eps_\theta$, if $n \geq c_{\theta,\eps,d} \log(1/\delta) \sigma^{2d}$, then on the basis of $n$ i.i.d.\ samples from $\mathrm{P}_\theta$ we can produce $M$ estimators $\hat \theta_1$, \dots $\hat \theta_M$ such that, with probability at least $1-\delta$, there exists a permutation $\pi: [M] \to [M]$ satisfying
\begin{equation*}
d_G(\hat \theta_i, \mathfrak{o}_{\pi(i)}) \leq \eps
\end{equation*}
for all $i \in [N]$.
\end{theorem}

Theorem~\ref{thm:statistical-upper} essentially follows from the observation that, since the variance of $y_i$ is $O(\sigma^2)$, a degree-$d$ polynomial in the entries of $y_i$ has variance $O(\sigma^{2d})$.
This implies that for any $f \in U^T_{\leq d}$, the evaluation $f(\theta)$ can be accurately estimated on the basis of $O(\sigma^{2d})$ samples.
By inverting a suitable polynomial system, we can thereby identify $\Mthetad$, at least approximately.
A full proof of Theorem~\ref{thm:statistical-upper} appears in Section~\ref{sec:proofs-stat}.

We now prove a lower bound showing that the dependence on $\sigma$ in Theorem~\ref{thm:statistical-upper} is tight.
We show that if $U_{\leq d-1}^T$ fails to resolve (or list-resolve) $\theta$, then $\Omega(\sigma^{2d})$ samples are necessary to recover (or list-recover) the orbit of $\theta$. Together with Theorem~\ref{thm:statistical-upper}, this lower bound implies that if $d^*$ is the smallest positive integer for which $U_{d^*}^T$ resolves (or list-resolve) $\theta$, then $\Theta(\sigma^{2d^*})$ samples are required to recover (or list-recover) the orbit of $\theta$.
We make this lower bound precise in Theorem~\ref{thm:statistical-lower}.

\begin{theorem}\label{thm:statistical-lower}
For any positive integer $d$, there exists a constant $c_d$ such that if $\tau_1$ and $\tau_2$ are elements in $\mathcal{M}_{\theta, d-1}$ lying in different orbits, then no procedure can distinguish between $\mathrm{P}^n_{\tau_1}$ and $\mathrm{P}^n_{\tau_2}$ with probability greater than $2/3$ if $n \leq c_d \sigma^{2d}$.
\end{theorem}
Via a standard argument known as Le Cam's method~\cite{LeC73}, Theorem~\ref{thm:statistical-lower} translates into lower bounds on the number of samples necessary to recover~$\theta$ in the $\sigma \to \infty$ limit.
This lower bound is of a slightly different form from the upper bound: Theorem~\ref{thm:statistical-upper} establishes that there exists a sequence of estimators such that, for a fixed $\theta \in V$, if $U_{\leq d}^T$ resolves $\theta$ and $n = \omega(\sigma^{2d})$ as $\sigma \to \infty$, then these estimators converge to $\theta$.
Theorem~\ref{thm:statistical-lower} implies that, if $U_{\leq d-1}^T$ fails to resolve $\theta$ and $n = o(\sigma^{2d})$ as $\sigma \to \infty$, then no sequence of estimators can succeed uniformly on all signals in the order-$(d-1)$ matching set of $\theta$.

The proof of Theorem~\ref{thm:statistical-lower} relies on a tight bound for the Kullbeck-Leibler divergence between the distributions $\mathrm{P}_{\mathfrak{o}_1}$ and $\mathrm{P}_{\mathfrak{o}_2}$ established in~\cite{BanRigWee17}. More details appear in Section~\ref{sec:proofs-stat}.

\section{Algebraic and algorithmic criteria}
\label{sec:algebraic}

In this section, we consider the four recovery criteria defined in Section~\ref{sec:4recoveryconditions}, and give sufficient algebraic conditions for their satisfaction in each case. The results of Section~\ref{sec:stat} imply that it suffices to focus our attention on deciding when a subspace $U$ of invariant polynomials resolves (or list-resolves) a parameter $\theta$. We show below how to answer this question by purely algebraic means. Moreover, for generic list recovery, we show how this question can be answered algorithmically in polynomial time. For the other recovery criteria, we also give algorithms to decide the corresponding algebraic condition; however, these algorithms are not efficient.

Throughout, we assume the setup defined in Section~\ref{sec:gen-prob} for the generalized orbit recovery problem. In particular, $G$ is a compact Lie group acting continuously and orthogonally on a finite-dimensional real vector space $V$. (It costs no generality to also assume the action is faithful, because the group only enters the picture through its action on $V$.)
We have the invariant ring $\RxG$ corresponding to the action of $G$ on $V$, and a subspace $U \subseteq \RxG$ of invariants that we have access to (e.g.\ $U^T_{\le d}$, the invariants of degree at most $d$). In this section we require $U$ to be finite-dimensional over $\RR$ and to be {\em graded}, i.e., to have a basis consisting of homogeneous elements. We are interested in whether the values $f(\theta)$ for $f \in U$ determine the orbit of $\theta \in V$ under $G$. The specific structure of $G$ and of $U^T_{\le d}$ (as defined in Section~\ref{sec:gen-prob}) will largely be unimportant in this section.  The exception is for heterogeneous problems with $G = \tilde{G} \wr S_{K}$, where it pays off to relate these problems to their homogeneous (single signal) versions.

An important intended audience of this paper is researchers interested in the orbit recovery problems treated here, who may not have a background in algebraic geometry. While this section uses some algebraic geometry and invariant theory, we have made an effort to write the statements and arguments in such a way that they are readable without this background. But we have also tried to avoid belaboring information well-known to more algebraically-oriented readers. With this in mind, we sometimes write sentences that begin, ``The reader without a background in algebraic geometry should picture...," or similar, and also include references to Appendix~\ref{app:algprimer}, which serves as a primer on algebra and invariant theory. The reader trained in algebraic geometry can ignore these. Meanwhile, sentences that begin, ``Formally, ..." are intended to pin down precise definitions, but can safely be ignored by readers without an algebraic geometry background. Hopefully we have struck a balance that is reasonable to all intended audiences.

\subsection{Invariant theory and algebraic geometry setup}\label{sec:it-ag-setup}

This section develops the algebraic-geometric setup for the proofs of our algebraic results. Except for Definitions~\ref{def:reynolds} and \ref{def:hilbertseries}, it can be skipped by readers interested in the statements but not the proofs. (For the results on generic list recovery, Theorems~\ref{thm:generic-list} and \ref{thm:bound-list}, alternative proofs that avoid algebraic geometry are given in Sections~\ref{sec:generic-proof} and \ref{sec:generic-converse}.)

Here is a succinct if imprecise summary of the main idea of the proofs in this section. The point of view here is standard in invariant theory; see for example \cite{GIT, pop-vin}. Recall that $V$ is the vector space of signals, and that $U$ is a finite-dimensional graded vector space of invariant polynomials to which we have access. We construct a space $X$ that parametrizes the orbits, along with a map
\[
\pi:V\rightarrow X
\]
that sends each orbit in $V$ to the corresponding point of $X$. Then we construct another space $Y$ parametrizing the subsets of $V$ that the functions in $U$ are able to distinguish, along with a map 
\[
\varphi:X\rightarrow Y
\]
sending two points of $X$ (corresponding to two orbits in $V$) to different points in $Y$ if and only if some function in $U$ evaluates differently on them. Then the ability of $U$ to resolve orbits is determined by properties of the map $\varphi$. If it is injective, $U$ resolves orbits uniquely. If it has finite fibers, $U$ list-resolves orbits. On the other hand, if $\varphi$ has infinite fibers, then $U$ cannot list-resolve orbits. Etc.

This summary is qualified with the word ``imprecise" because it elides significant details. The spaces $X$ and $Y$ will be algebraic varieties defined over $\RR$ (also known as $\RR$-varieties). These have $\CC$-points and $\RR$-points. Readers without a background in algebraic geometry should picture them as subsets of $\CC^N$ defined by the vanishing of sets of polynomials whose coefficients lie in $\RR$; then the $\CC$-points (written $X(\CC)$ for a variety $X$) are all the points,\footnote{The $\RR$-variety $X$ itself is distinguished from its set of $\CC$-points $X(\CC)$. The latter is a bare topological space, while the former also contains the data of which functions on this space are restricted from polynomials on $\CC^N$, and among these, which are real polynomials.} while the $\RR$-points (written $X(\RR)$) are the ones that lie within $\RR^N\subset \CC^N$. On the one hand, if the picture described above were concerned with all the $\mathbb{C}$-points, then information about the size of the fibers of $\varphi$ would readily be accessed using standard theorems of algebraic geometry. For example if the dimension of $Y$ is lower than that of $X$, then the fibers of the map $\varphi$ have positive dimension as varieties, and positive-dimensional varieties have infinitely many $\CC$-points. On the other hand, we are not actually concerned with all the $\CC$-points. It will only be the subset of $X$ consisting of the image of the $\RR$-points of $V$ under $\pi$ that actually parametrizes the orbits of $G$ on $V$. Likewise, only the subset of $Y$ consisting of the image of the $\RR$-points of $V$ under $\varphi\circ\pi$ actually parametrizes the classes of orbits distinguishable by elements of $U$. The main work in most of the proofs, and particularly the proof of Theorem~\ref{thm:generic-list}, is to make sure there is enough control over the $\RR$-points that general theorems about the varieties can be translated into desired conclusions about these subsets.

We preface the discussions of the individual recovery criteria with some general geometric and algebraic considerations on the problem of distinguishing orbits with invariants, which we articulate as a series of discussions labeled ``Setup" for later cross-referencing. 

\begin{setup}\label{rmk:AGsetup} Because $G$ is a compact subgroup of $GL(V)$, 
it is the set of real points of a linear algebraic group defined over $\RR$ \cite[Chapter~3, \S~4, Theorem~5]{on-vin}. 
View $V$ as an affine $\RR$-variety,\footnote{Formally, when we write ``affine $\RR$-variety", we mean a reduced affine scheme of finite type over $\RR$, so that we have access to the functor of points.} on which $G$ acts algebraically by $\RR$-automorphisms. (We abuse notation by using the same symbol $V$ for both the original $\RR$-vector space, and the $\RR$-variety [which also has $\CC$-points]. The points of the original vector space $V$ can be recovered as the $\RR$-points of the variety $V$.) Because $G$ is compact, a Reynolds operator $\mathcal{R}:\Rx\rightarrow\RxG$ (see Definition~\ref{def:reynolds} below) exists, and it follows by a standard argument (see for example \cite[Theorem~2.2.10]{dk-book}, or \ref{thm:finite-gen} in the Appendix) that the ring $\RxG$ is finitely generated. Thus it defines an affine $\RR$-variety, the {\em categorical quotient} of $V$ by $G$, which we will call $X$. (For the reader without a background in algebraic geometry, a summary of the construction of $X$ from $\RxG$ is given in Appendix~\ref{app:someAG}.) The subspace $U\subset \RxG$ being finite-dimensional, the subring $\RR[U]$ of $\RxG$ that it generates over $\RR$ is also finitely generated, so it too defines an affine $\RR$-variety $Y$, and the ring inclusions 
\[
\RR[U]\subset\RxG\subset\Rx
\]
define dominant morphisms of $\RR$-varieties
\[
\begin{CD}
V @>\pi>> X @>\varphi>> Y.
\end{CD}
\]
The varieties $V,X,Y$ are irreducible. (The construction of these maps from the ring inclusions $\RR[U]\subset\RxG\subset\Rx$, and the notions of dominant morphism and irreducible variety, are also discussed in Appendix~\ref{app:someAG}.)

The $G$-action on $\Rx$ is trivial when restricted to $\RxG$ (by definition of the latter); thus for any $g\in G$, we have an equality of morphisms $\pi = \pi\circ g$. The morphisms $\pi,\varphi$ induce maps of the $\RR$-points, which we denote by $\pi_\RR, \varphi_\RR$:
\[
\begin{CD}
V(\RR) @>\pi_\RR>> X(\RR) @>\varphi_\RR>> Y(\RR).
\end{CD}
\]
The equation $\pi = \pi\circ g$ implies that the fibers of the map $\pi_\RR:V(\RR)\rightarrow X(\RR)$ are unions of $G$-orbits. 
\end{setup}

\begin{setup}\label{rmk:separation-in-terms-of-AG}
With this setup, the statement that $U$ distinguishes between two orbits of $G$ is precisely the statement that these orbits lie in distinct fibers of the composed map $\varphi_\RR\circ\pi_\RR:V(\RR)\rightarrow Y(\RR)$. Because $G$ is compact, the full invariant ring $\RxG$ resolves every orbit of $G$ on the original vector space $V$, i.e., the map $\pi_\RR$ distinguishes between all pairs of $G$-orbits on $V(\RR)$ (see, e.g., \cite[Chapter~3, \S~4, Theorem~3]{on-vin}, or \ref{thm:RxG-separates} in the Appendix).\footnote{This statement depends both on the compactness of $G$ and the fact that it is acting on a real vector space. It does not hold for noncompact groups, even reductive ones. Nor does it hold when the action is by unitary transformations on a $\CC$-vector space, even when the group is compact (unless one enriches the ring of polynomial functions by including the antiholomorphic polynomials, essentially reverting to the real case, before taking invariants). The famous example of a reductive (but noncompact) group whose invariants fail to separate orbits, is the action of $\RR^\times=GL_1(\RR)$ by scalars on $\RR^2$. The only invariant polynomials under this action are the constants, so no two orbits are distinguished by them. Then $X(\RR)$ consists of a single point, although there are infinitely many orbits. A slight variation on this example yields a compact group acting on a complex vector space such that the polynomial invariants fail to separate orbits, namely $SU(1)$ acting on $\CC^2$ by scalars. Again, there are no non-constant invariants.}\label{note:noncompact-do-not-separate} Put another way, any orbit of $G$ on $V(\RR)$ corresponds uniquely to an $\RR$-point of $X$ lying in the image of $\pi_\RR$, with the correspondence given by $\pi_\RR$. Therefore, the ability of a subset $U\subset\RxG$ to distinguish a pair of orbits comes down to the ability of the map $\varphi_\RR:X(\RR)\rightarrow Y(\RR)$ to distinguish the corresponding pair of elements of $X(\RR)$. So we will be interested, in the below, in controlling the size of the fibers of $\varphi:X(\RR)\rightarrow Y(\RR)$. For example, if $\varphi_\RR$ is injective then $U$ must resolve every orbit of $G$ in $V$ uniquely, and if $\varphi_\RR$ has finite fibers, then $U$ must list-resolve every orbit of $G$.
\end{setup}

\begin{setup}\label{rmk:pathologies}
By Setup \ref{rmk:separation-in-terms-of-AG}, the image of the map $\pi_\RR$ is set-theoretically the quotient of $V(\RR)$ by $G$. However, $\pi_\RR$ need not be surjective onto $X(\RR)$, and in general $\varphi_\RR$ is neither injective nor surjective, even in the case that $U$ distinguishes all orbits. Here is an illustrative example: Let $G=\ZZ/2$, and let $V=\RR^1$, with $G$ acting by the sign representation (i.e., the nontrivial group element reverses the sign of $v\in V$). Then $\RxG=\RR[x^2]$. Let $U$ be the subspace of $\RxG$ spanned by $x^4$. Then $\Rx$, $\RxG$, and $\RR[U]$ are all univariate polynomial rings, so $V$, $X$, and $Y$ are all isomorphic to the affine line over $\RR$, and $\pi$ and $\varphi$ are both the squaring map, which is neither injective nor surjective when restricted to the real points. Nonetheless, if $x$ is real, knowing $x^4$ specifies $x$ up to sign, so $U$ uniquely resolves every orbit of $G$ in this case.
\end{setup}

\begin{setup}\label{rmk:X(C)etc}
We will use general theorems from algebraic geometry to first get information about the $\CC$-points of $X$ and $Y$, on the way to the information we need about the $\RR$-points.  The sets $X(\CC)$ and $Y(\CC)$ carry Zariski topologies. (Formally, $X(\CC)$ and $Y(\CC)$ are in bijection with the closed points of the schemes $X\times_\RR\CC$ and $Y\times_\RR\CC$; the Zariski topologies on the latter then induce topologies on the former via this bijection. For the reader without a background in algebraic geometry, the Zariski topology is described in Definition~\ref{def:zariski} in Appendix~\ref{app:someAG}.) As with $\RR$, we refer to the maps induced by $\pi$ and $\varphi$ on $V(\CC)$ and $X(\CC)$ as $\pi_\CC$ and $\varphi_\CC$. 

An example of how working over $\CC$ simplifies the geometry is that $\pi_\CC$, unlike $\pi_\RR$, is always surjective. The complexification $G_\CC$ of $G$ is reductive \cite[Chapter~5, \S~2, Theorem~12]{on-vin}. 
The invariant ring $\Cx^{G_\mathbb{C}}$ is nothing but $\CxG = \CC\otimes_\RR \RxG$.\footnote{The fact that $G$ and $G_\CC$ have the same invariant ring in $\Cx$ is not completely obvious because $G_\CC$ is a larger Lie group than $G$, so it might seem its invariant ring should be smaller. However, the image of $G$ under the natural embedding of Lie groups $G\hookrightarrow G_\CC$ is Zariski-dense because $\RR$ is a {\em large field} (see Lemma~\ref{lem:largefield}). If $g\in G_\CC$, then for fixed $f\in \Cx$, the equation $gf=f$ is a closed condition on $g$. Therefore, if $f$ is fixed under the action of $G$, it is fixed under the action of all elements of $G_\CC$.}\label{note:same-invariant-ring} So surjectivity of $\pi_\CC$ follows from \cite[Theorem~4.6]{pop-vin}.

There are natural injections $i_X:X(\RR)\rightarrow X(\CC)$ and $i_Y:Y(\RR)\rightarrow Y(\CC)$, and they make the following diagram commute:
\[
\begin{CD}
X(\RR) @>\varphi_\RR>> Y(\RR)\\
@Vi_XVV @VVi_YV\\
X(\CC) @>>\varphi_\CC> Y(\CC)
\end{CD}
\]
(Formally, $i_X$ is composition of an $\RR$-point $p:\operatorname{Spec}\RR\rightarrow X$ with the canonical morphism $\operatorname{Spec}\CC\rightarrow \operatorname{Spec}\RR$, and $i_Y$ is defined similarly. The reader can picture the obvious inclusion of the set of points with only real coordinates, $X(\RR)$, into the set of all the points $X(\CC)$.) Because $i_X,i_Y$ are injections, they induce subspace topologies on $X(\RR)$ and $Y(\RR)$. For the same reason, $\varphi_\RR$ can be viewed as the restriction of $\varphi_\CC$ to $X(\RR)$; thus an upper bound on the size of the fibers of $\varphi_\CC$ immediately yields the same upper bound on the size of the fibers of $\varphi_\RR$.

(In a similar way, $V(\CC)$ and $V(\RR)$ can be viewed as topological spaces. The topological space $V(\RR)$ is exactly the underlying set of the original $\RR$-vector space $V$, equipped with the Zariski topology.)
\end{setup}

\begin{setup}
The following three lemmas are the basic tools we use to gain adequate control of the $\RR$-points of $V$, $X$, and $Y$.

\begin{lemma}\label{lem:largefield}
If $X$ is an irreducible variety defined over $\RR$, and $X(\RR)$ contains a smooth point of $X$, then $X(\RR)$ is Zariski-dense in $X$.
\end{lemma}

\begin{proof}
The field $\RR$ is a {\em large field} by \cite[Part~I, Section~1, A2 and A3]{pop}. For any large field $k$, if $X$ is an irreducible $k$-variety with a smooth $k$-point, then the $k$-points are Zariski-dense \cite[Proposition~2.6]{pop}.
\end{proof}

To apply Lemma~\ref{lem:largefield}, we will need to be able to guarantee the existence of smooth real points in fibers of maps. This will be possible because after removing a Zariski-closed set from the domain, every map of algebraic varieties has fibers that are completely smooth.

\begin{lemma}[generic smoothness]\label{lem:generic-smoothness-of-fibers}
If $f: X \rightarrow Y$ is a dominant morphism of irreducible varieties defined over $\RR$, then there is an open subvariety $U$ of $X$, defined over $\RR$, such that the restriction of $f$ to $U$ has smooth fibers. 
\end{lemma}

\begin{proof}
First, let $V$ be the regular (nonsingular) locus of $Y$; it is a nonempty open subset. Then $f^{-1}(V)$ is open in $X$, and because $f$ is dominant, it is nonempty. Replace $Y$ with $V$ and $X$ with $f^{-1}(V)$, so $Y$ is now a regular variety. The extension $\RR(X)/\RR(Y)$ of function fields is separable because $\RR(Y)$ is of characteristic zero. Thus the hypotheses of \cite[Corollary~5.4.3]{mumford-oda} are met, so there is a point $x\in X$ at which $f$ is smooth, whereupon it follows that there exists an open subset $U\subset X$ containing $x$ (and defined over $\RR$) such that the morphism $f|_U:U\rightarrow Y$ is smooth \cite[Proposition~5.3.2]{mumford-oda}. Since smooth morphisms have smooth fibers \cite[Proposition~5.3.3(i)]{mumford-oda}, this is the desired $U$.
\end{proof}

\begin{lemma}\label{lem:denseimage}
If $f:X\rightarrow Y$ is a dominant morphism of varieties, and $S\subset X$ is a Zariski-dense set of points of $X$, then $f(S)$ is Zariski-dense in $Y$.
\end{lemma}

\begin{proof}
Let $U$ be any nonempty open subset of $Y$. Since $f$ is dominant, $U$ meets $f(X)$, so $f^{-1}(U)$ is a nonempty open subset of $X$. Because $S$ is dense in $X$, $f^{-1}(U)$ then meets $S$. Therefore $U$ meets $f(S)$. We conclude $f(S)\subset Y$ is dense.
\end{proof}

\begin{remark}
There is nothing specifically algebro-geometric about Lemma~\ref{lem:denseimage}. The same argument shows in general that if $f:X\rightarrow Y$ is a continuous map of topological spaces with dense image, then $f$ maps dense subsets to dense subsets.
\end{remark}
\end{setup}

\begin{setup}
We define two fundamental concepts from invariant theory that are used in the sequel: the Reynolds operator and the Hilbert series.
\begin{definition}[Reynolds operator] \label{def:reynolds} For a compact Lie group $G$ acting gradedly on $\Rx$, the \emph{Reynolds operator} $\mathcal{R}: \Rx \to \RxG$ is defined by $$\mathcal{R}(f) = \Ex_{g \sim \Haar(G)} [g \cdot f].$$
\end{definition}

\noindent Note that the Reynolds operator is a linear projection from $\Rx$ to $\RxG$ that preserves the degree of homogeneous polynomials (i.e., a homogeneous polynomial of degree $d$ gets mapped either to a homogeneous polynomial of degree $d$ or to zero).

\begin{observation} \label{obs:basis}
Let $\RxG_d$ denote the vector space consisting of homogeneous invariants of degree $d$. We can obtain a basis for $\RxG_d$ by applying $\mathcal{R}$ to each monomial in $\RR[{\bf x}]$ of degree $d$. (This yields a spanning set which can be pruned to a basis if desired.)
\end{observation}

\begin{definition}[Hilbert series]\label{def:hilbertseries}
Let $\Rx_d^G$ be the subspace (over $\RR$) of $\RxG$ consisting of homogeneous invariants of degree $d$. The \emph{Hilbert series} of $\RxG$ is the formal power series
$$H(t) \defeq \sum_{d=0}^\infty \dim(\Rx_d^G)\, t^d.$$
\end{definition}
\end{setup}

\begin{remark}
In what follows we will be discussing algorithms that take the problem setup as input (including $\tilde G$ and its action on $\tilde V$, along with $\Pi,K$) and decide whether or not $U^T_{\le d}$ (for some given $d$) is capable of a particular recovery task (e.g.\ list recovery of a generic $\theta \in V$). We will always assume that these algorithms have a procedure to compute a basis for $U^T_d$ (for any $d$) in exact symbolic arithmetic. It is always possible in principle to do this if the group and the representation are specified in exact symbolic arithmetic, since this implies they are defined over the algebraic closure $\overline\QQ$ of $\QQ$, in which case, in view of the fact that $G$ is reductive (Setup~\ref{rmk:X(C)etc}) and therefore linearly reductive (as $\overline\QQ$ has characteristic zero), the Reynolds operator is defined over $\overline\QQ$ as well \cite[Theorem~2.2.5]{dk-book}. For the important case of $\mathrm{SO}(3)$, see Appendix~\ref{app:so3} for an explicit basis for the invariants.
\end{remark}

\subsection{Generic list recovery}\label{sec:generic-list}

\noindent We now present our algebraic characterization of the generic list recovery problem. For a subset $B\subset \Rx$, we denote by $\trdeg(B)$ the transcendence degree of $B$, i.e., the maximum cardinality of subsets of $B$ that are algebraically independent over $\RR$. Recall that a set has \emph{full measure} if its complement has measure zero. 

\begin{theorem}[generic list recovery]\label{thm:generic-list}
Let $U \subseteq \RxG$ be a finite-dimensional subspace. If $\trdeg(U) = \trdeg(\RxG)$ then there exists a set $S \subseteq V$ of full measure such that if $\theta \in S$ then $U$ list-resolves $\theta$. Conversely, if $\trdeg(U) < \trdeg(\RxG)$ then there exists a set $S \subseteq V$ of full measure such that if $\theta \in S$ then $U$ does not list-resolve $\theta$.
\end{theorem}

\begin{remark}
A {\em minimal} subset $B$ of $\RxG$ such that $\trdeg(B)=\trdeg(\RxG)$ is called a {\em transcendence basis}, thus Theorem~\ref{thm:generic-list} could have alternatively been stated that $U$ achieves generic list recovery if and only if it contains a transcendence basis. The intuition behind Theorem~\ref{thm:generic-list} is that $\trdeg(\RxG)$ is the number of degrees of freedom that need to be pinned down in order to learn the orbit of $\theta$, and so we need this many algebraically independent constraints (invariant polynomials). If the orbits were parametrized by the $\CC$-points of $X$, this would essentially already be a proof: by standard results in algebraic geometry, the dimensions of $X$ and $Y$ are the transcendence degrees of $\RxG$ and $\RR[U]$, and the dimension of the generic fiber of $\varphi:X\rightarrow Y$ is the difference between them. Thus if $U$ has full transcendence degree, the generic fiber is an affine variety of dimension zero, which is finite; while if $U$ has deficient transcendence degree, the generic fiber has positive dimension so is infinite. The reason the proof below requires more care is that the orbits are actually parametrized not by $X(\CC)$ but only by the subset of $X(\RR)$ lying in the image of $\pi_\RR$.
\end{remark}

We give here a proof based on algebraic geometry. An alternative proof, also encompassing Theorem~\ref{thm:bound-list} below (which serves as a lemma for Theorem~\ref{thm:generic-list}), can be found in Section~\ref{sec:algebraic-proofs}. The alternative proof of Theorem~\ref{thm:bound-list} is via elementary field theory and is found in Section~\ref{sec:generic-proof}. The alternative proof of the rest of Theorem~\ref{thm:generic-list} is via differential topology, in Section~\ref{sec:generic-converse}.

\begin{proof}[Proof of Theorem~\ref{thm:generic-list}]
We use the notations of Setup~\ref{rmk:AGsetup}--\ref{rmk:X(C)etc}. It suffices to demonstrate the existence of nonempty Zariski-open subsets $S$ of $V(\RR)$ with the stated properties, as such subsets are of full measure.

We have $\trdeg_\RR(\RxG) = \dim X$ and $\trdeg_\RR(U)=\dim Y$. We write $\RR(X),\RR(Y)$ for the fields of rational functions on $X,Y$ with real coefficients. (These are just the fraction fields of $\RxG,\RR[U]$). Let 
\[
r \defeq\trdeg_{\RR(Y)} \RR(X) = \dim X - \dim Y.
\]
With this notation, the statement is that $U$ generically list-resolves $\theta\in V$ if $r=0$, and generically fails to list-resolve $\theta$ if $r>0$.

First suppose $r>0$. By generic smoothness (Lemma \ref{lem:generic-smoothness-of-fibers}) applied to the morphism $\varphi\circ\pi:V\rightarrow Y$, there is a nonempty open subvariety $S'\subset V$ defined over $\RR$, such that the restricted morphism $\varphi\circ \pi|_{S'}:S'\rightarrow Y$ has smooth fibers. Let $S\defeq S'(\RR)\subset V(\RR)$ be the set of real points of $S'$. This is nonempty because $V(\RR)$ is dense in $V$ (by Lemma~\ref{lem:largefield}). We will show that for $\theta\in S$, $U$ fails to list-resolve $\theta$, so that $S$ is the claimed nonempty open subset of $V(\RR)$.

Suppose $\theta\in S$. Let $y=\varphi_\RR\circ\pi_\RR(\theta)\in Y(\RR)$. Consider the fiber $V_y\defeq(\varphi\circ\pi)^{-1}(y)$, which is a closed subvariety of $V$ defined over $\RR$. Because $\pi:V\rightarrow X$ is surjective (see Setup~\ref{rmk:X(C)etc}), the restriction of $\pi$ to $V_y$ is surjective and therefore dominant onto the fiber $X_y\defeq\varphi^{-1}(y)\subset X$.

The fiber $S'_y=(\varphi\circ\pi|_{S'})^{-1}(y)$ of the restriction of $\varphi\circ\pi$ to $S'$ is an open subvariety of $V_y$ defined over $\RR$. By construction it is smooth and contains the real point $\theta$. Therefore, by Lemma~\ref{lem:largefield}, the set of real points $S'_y(\RR)$ forms a dense subset of $S'_y$. Since $S'_y$ is Zariski-open in $V_y$, $S'_y(\RR)$ is also dense in $V_y$. It follows that the image $\varphi_\RR(S'_y(\RR))$ is Zariski-dense in $X_y$ (per Lemma~\ref{lem:denseimage}). Because $X_y$ has dimension at least $r$ (e.g., \cite[Corollary~4.5.5]{mumford-oda}), which is positive by assumption, it follows that $\varphi_\RR(S'_y(\RR))$ is infinite. By the considerations in Setup~\ref{rmk:separation-in-terms-of-AG}, each of the infinitely many points of this set comes from a distinct $G$-orbit in $V(\RR)$, but they all have the same image $y$ in $Y$, i.e., the invariants in $U$ cannot distinguish them from each other. Thus $\theta$ is not list-resolved by $U$. Since $\theta\in S$ was arbitrary, $S$ is the claimed nonempty open subset of $V(\RR)$ on which $U$ fails to list-resolve.

Now consider the case that $r=0$. In this case, the field extension $\RR(X)/\RR(Y)$ is algebraic. Since, in addition, $\RR(X)$ is finitely generated as a field over $\RR$, it follows that $\RR(X)$ has finite degree as an extension of $\RR(Y)$. Therefore, the desired statement is immediate from the following sharper result, which gives an upper bound on the list size.
\end{proof}

\begin{theorem}\label{thm:bound-list}
Let $U$ be a finite-dimensional subspace of the invariant ring $\RxG$. Let $\RR(X)$ be the field of fractions of $\RxG$ and let $\RR(Y)$ be the field of fractions of $\RR[U]$. If the extension degree $[\RR(X) : \RR(Y)] = D < \infty$, then there exists a set $S \subseteq V$ of full measure such that for any $\theta \in S$, $U$ list-resolves $\theta$ with a list of size $\le D$.
\end{theorem}

\noindent Here $[\RR(X) : \RR(Y)]$ denotes the degree of a field extension. Appendix~\ref{app:fields} collects all the definitions and facts from field theory used in this paper.

\begin{proof}
The field extension $\RR(X)/\RR(Y)$ is separable since $\RR(Y)$ has characteristic zero. It follows (e.g., \cite[Theorem~5.1.6(iii)]{springer}) that there is a nonempty open subvariety $S'\subset X$ such that for every $x\in S'(\CC)$, the cardinality of the fiber $\varphi_\CC^{-1}(\varphi_\CC(x))$ is $D$, the extension degree. By the diagram in Setup~\ref{rmk:X(C)etc}, $\varphi_\RR$ can be viewed as a restriction of $\varphi_\CC$; thus if $x\in S'(\RR)$, the cardinality of $\varphi_\RR^{-1}(\varphi_\RR(x))$ is less than or equal to $D$. 

The preimage $\pi^{-1}(S')$ is a nonempty open subvariety of $V$, and it has real points because $V(\RR)$ is dense in $V$. The set $S$ of these real points is the desired nonempty open subset of $V(\RR)$. Indeed, for any $\theta\in S$, let $y=\varphi_\RR\circ\pi_\RR\in Y(\RR)$. Then the orbits that $U$ cannot distinguish from $\theta$'s biject with the points in the intersection of $\pi_\RR(V(\RR))$ with $\varphi_\RR^{-1}(y)$ (see Setup~\ref{rmk:separation-in-terms-of-AG}), but we have just seen that the latter set has cardinality $\le D$, so the intersection does too. Thus $\theta$ is list-resolved by $U$ with a list size of $\le D$.
\end{proof}

In order for Theorem~\ref{thm:generic-list} to be useful, we need a way to compute the transcendence degree of both $\RxG$ and $U$. In what follows, we discuss methods for both of these: in Section~\ref{sec:trdeg-R} we collect together ways of computing $\trdeg(\RxG)$, and in Section~\ref{sec:trdeg-U} we give an efficient algorithm to compute $\trdeg(U)$ for a subspace $U$. By taking $U = U^T_{\le d}$ this yields an efficient algorithm to determine the smallest degree $d$ at which $U^T_{\le d}$ list-resolves a generic $\theta$ (thereby answering Question~\ref{que:alg} for the case of generic list recovery).

\subsubsection{Computing the transcendence degree of $\RxG$} \label{sec:trdeg-R}

This subsection collects together tools that can be used in practice to find the transcendence degree of $\RxG$ in the context of a given orbit recovery problem.

Intuitively, the transcendence degree of $\RxG$ is the number of parameters required to describe an orbit of $G$; in other words, the dimension of the orbit space $V/G$. For finite groups, this is simply the dimension of $V$:

\begin{proposition}[\cite{sturmfels} Proposition~2.1.1]
\label{prop:trdeg-finite}
If $G$ is a finite group, $\trdeg(\RxG) = \dim(V)$.
\end{proposition}

For infinite groups, the situation may be slightly different; the transcendence degree will not  be determined solely by the dimension of $V$ or even the dimensions of $V$ and $G$. For instance, if $\mathrm{SO}(3)$ acts on $V = \RR^3$ in the canonical way (rotations in 3 dimensions), then a generic orbit is a sphere, with dimension two. This means there is only one parameter to learn, namely the 2-norm, and $\RxG$ has transcendence degree $\dim V-2=1$ accordingly. On the other hand, if $\mathrm{SO}(3)$ acts on a rich class of functions $S^2 \to \RR$ (as in the $S^2$ registration problem; see Section~\ref{sec:ex-s2}) then the generic orbit resembles a copy of $\mathrm{SO}(3)$, which has dimension 3, so the orbit space will have dimension $\dim V - 3$. The difference in the size of a generic orbit between the two actions corresponds with a difference in the size of a generic point stabilizer. This is formalized in the following.

\begin{proposition}
\label{prop:trdeg-general}
If $G$ is an algebraic group, then
$$\trdeg(\RxG) = \dim(V) - \dim(G) + \min_{v \in V} \dim(G_v),$$
where $G_v$ is the stabilizer of $v$ for the action of $G$ (that is, the subgroup of all $g \in G$ fixing $v$).
\end{proposition}

\begin{proof}
    The statement would be precisely \cite[Corollary~6.2]{dolgachev}, if the field of invariant rational functions $\RR(\mathbf{x})^G$ had appeared in the place of the ring $\RxG$ of invariant polynomials; thus the proposition follows if we show that $\trdeg(\RxG)=\trdeg(\RR(\mathbf{x})^G$. For real actions of compact groups, this holds, because in fact $\RR(\mathbf{x})^G$ is the field of fractions of $\RxG$. This is seen as follows:

    If $G^0$ is the connected component of the identity in $G$, then any continuous homomorphism $\chi:G^0\rightarrow \RR^\times$ is trivial, because the only compact, connected subgroup of the multiplicative group of reals is the trivial subgroup. By a slight variant of the argument in \cite[Theorem~3.3(b)]{pop-vin}, it then follows that $\RR(\mathbf{x})^G = \operatorname{Frac} \RxG$.
\end{proof}

An alternate approach to the transcendence degree of $\RxG$ uses the Hilbert series $H(t)$ of $\RxG$ (recall Definition~\ref{def:hilbertseries}). One can extract $\trdeg(\RxG)$ from $H(t)$ as follows.

\begin{proposition}\label{prop:pole-order}
The order of the pole at $t=1$ of $H(t)$ is equal to $\trdeg(\RxG)$.
\end{proposition}
\noindent This is a standard fact about finitely generated graded algebras over a field. One can find a proof, e.g., in \cite{dk-book}. For a given $G$ acting on $V$, there is an explicit formula (\emph{Molien's formula}) for the Hilbert series:

\begin{proposition}[\cite{kac-notes} Remark~3-1.8] \label{prop:molien}
Let $\rho: G \to \mathrm{GL}(V)$ be the representation by which $G$ acts on $V$. Then for $|t| < 1$, $H(t)$ converges and we have
\[
H(t) = \Ex_{g \sim \Haar(G)} \det(I-t \,\rho(g))^{-1}.
\]
\end{proposition}

\noindent This formula is tractable to compute, even for complicated groups; see Section~\ref{sec:ex-s2} for details in the case of $\mathrm{SO}(3)$.

For heterogeneous problems $(K > 1)$, the transcendence degree can be computed easily from the transcendence degree of the corresponding homogeneous $(K=1)$ problem.

\begin{proposition} \label{prop:trdeg-het}
Let $\tilde G$ be a compact group acting linearly and continuously on $\tilde V$, and let $G = \tilde G^K \rtimes S_K$ act on $V = \tilde V^{\oplus K} \oplus \bar\Delta_K$ as in Definition~\ref{def:het-action}. Let $\RxG$ be the invariant ring corresponding to the action of $G$ on $V$, and let $\RR[{\bf \tilde x}]^{\tilde G}$ be the invariant ring corresponding to the action of $\tilde G$ on $\tilde V$ (i.e.\ the $K=1$ problem). Then $\trdeg(\RxG) = K \cdot \trdeg(\RR[{\bf \tilde x}]^{\tilde G}) + K-1.$
\end{proposition}

The argument uses some commutative algebra. For readers without this background, the principles used are drawn out in Appendix~\ref{app:comm-alg}.

\begin{proof}
For a finitely generated algebra over a field, such as $\RxG$ or $\RR[{\bf \tilde x}]^{\tilde G}$, the transcendence degree is equal to the Krull dimension (this follows from Noether normalization). So it suffices to prove the statement for Krull dimension.

Now $\tilde G^K$ is a normal subgroup of $G$. Thus the invariant ring $\Rx^{\tilde G^K}$ is stable setwise under the action of $G$, the restriction of the $G$-action to this subring factors through $G/\tilde G^K \cong S_K$, and we have
\[
\RxG = \left(\Rx^{\tilde G^K}\right)^{S_K}.
\]
Because $S_K$ is a finite group, this equation implies $\Rx^{\tilde G^K}$ is integral over $\RxG$, and therefore they have the same Krull dimension. So we have to prove that $\dim_{\mathrm{Krull}} \Rx^{\tilde G^K} = K\cdot \dim_{\mathrm{Krull}} \RR[{\bf \tilde x}]^{\tilde G} + K-1$. But, indeed, as each of the $K$ factors of $\tilde G^K$ acts separately on a corresponding copy of $\tilde V$ in $V=\tilde V^{\oplus K}\oplus \overline\Delta_K$, the ring $\Rx^{\tilde G^K}$ is precisely the tensor product (over $\RR$) of $K$ copies of $\RR[{\bf \tilde x}]^{\tilde G}$ and a polynomial ring in $K-1$ ($=\dim_\RR \overline\Delta_K$) variables. In general, the Krull dimension of a tensor product of finitely generated algebras over a field is the sum of their individual dimensions. (Again, this follows from Noether normalization.) This completes the proof.
\end{proof}

\subsubsection{Algorithm for transcendence basis of $U$} \label{sec:trdeg-U}

We turn to the practical problem of finding the transcendence degree of $\RR[U]$. In this section we prove the following.
\begin{theorem}\label{thm:alg}
There is an efficient algorithm to perform the following task. Given a basis $\{u_1,\ldots,u_s\}$ for a finite-dimensional subspace $U \subseteq \Rx$, output a transcendence basis for $U$.
\end{theorem}

The proof exhibits the algorithm explicitly. It is given at the end of the section. First we recall the following simple classical test for algebraic independence (see, e.g.,~\cite{ER,blackbox} for a proof), and some of its consequences.

\begin{definition}[Jacobian]
Given polynomials $f_1,\ldots,f_m \in \RR[{\bf x}] = \RR[x_1,\ldots,x_p]$, the \emph{Jacobian matrix} $J_{\bf x}(f_1,\ldots,f_m) \in (\Rx)^{m \times p}$ is defined by $(J_{\bf x}(f_1,\ldots,f_m))_{ij} = \partial_{x_j} f_i$ where $\partial_{x_j}$ denotes formal partial derivative with respect to $x_j$.
\end{definition}

\begin{proposition}[Jacobian criterion for algebraic independence]
\label{prop:jac}
Polynomials ${\bf f} = (f_1,\ldots,f_m)$ are algebraically independent if and only if the Jacobian matrix $J_{\bf x}({\bf f})$ has full row rank (over the field $\RR({\bf x})$).
\end{proposition}

\noindent It suffices to test the rank of the Jacobian at a generic point ${\bf x}$.

\begin{corollary}\label{cor:subs}
Fix ${\bf f} = (f_1,\ldots,f_m)$. Let $z \sim \cN(0,I_{p \times p})$. If ${\bf f}$ is algebraically dependent then $J_{\bf x}({\bf f})|_{{\bf x} = z}$ does not have full row rank. If ${\bf f}$ is algebraically independent, then $J_{\bf x}({\bf f})|_{{\bf x} = z}$ has full row rank with probability 1.
\end{corollary}
\begin{proof}
An $m \times p$ matrix has deficient row rank if and only every $m\times m$ submatrix has determinant zero; thus the set $T$ of $z$-values for which $J_{\bf x}({\bf f})|_{{\bf x} = z}$ has deficient row rank is a Zariski-closed subset $Z$ of $\RR^p$. If the generic Jacobian $J_{\bf x}({\bf f})$ already has less than full row rank, then $Z=\RR^p$, so $z\in Z$ with certainty. If the generic Jacobian has full row rank, $Z$ is a proper closed subset, so of Lebesgue measure zero. Since $\cN(0,I_{p \times p})$ is absolutely continuous with respect to Lebesgue measure, $z$ avoids $Z$ with probability 1.
\end{proof}

\begin{remark}
In practice we may choose to plug in random \emph{rational} values for ${\bf x}$ so that the rank computation can be done in exact symbolic arithmetic.
If we manage to find \emph{any} value of ${\bf x}$ for which the Jacobian has full row rank, this constitutes a proof of algebraic independence.
On the other hand, if the Jacobian repeatedly fails to have full row rank for many random rational values of ${\bf x}$, we can conclude with high confidence that an algebraic dependence exists.
\end{remark}

\begin{remark}\label{rem:jac-num}
In some cases (e.g.,\ if the polynomials involve irrational values) it may be slow to compute the Jacobian rank in exact symbolic arithmetic. We can alternatively compute the singular values numerically and count how many are reasonably far from zero. This method works reliably in practice (i.e., it is extremely clear how to separate the zero and nonzero singular values) but does not constitute a rigorous proof of algebraic independence.
\end{remark}

Although the Jacobian criterion gives an efficient test for algebraic dependence, it is much harder ($\#P$-hard) to actually find the algebraic dependence (i.e., the polynomial relation) when one exists \cite{complexity-ann}.

The Jacobian criterion implies the well-known fact that the collection of algebraically independent subsets of $\Rx$ form a \emph{matroid}; this is called an \emph{algebraic matroid} (see, e.g.,\ \cite{schrijver}). In particular, we have the following exchange property:

\begin{proposition}\label{prop:exchange}
Let $I,J$ be finite subsets of $\Rx$, each algebraically independent. If $|I| < |J|$ then there exists $f \in J \setminus I$ such that $I \cup \{f\}$ is algebraically independent.
\end{proposition}

We next note that in the task from Theorem~\ref{thm:alg}, a transcendence basis can always be taken from the basis $\{u_1,\ldots,u_s\}$ itself.

\begin{lemma}\label{lem:basis}
Let $U$ be a finite-dimensional subspace of $\Rx$ with basis $B = \{u_1,\ldots,u_s\}$. If $U$ contains $r$ algebraically independent elements, then so does $B$.
\end{lemma}
\begin{proof}
Let $B' \subseteq B$ be a maximal set of algebraically independent elements of $B$. For a contradiction, suppose that $|B'| < r$. Then by the exchange property (Proposition~\ref{prop:exchange}) there exists $v \in U \setminus B'$ such that $B' \cup \{v\}$ is algebraically independent. Write $v = \sum_{i=1}^s \alpha_i u_i$. Because $B'$ is maximal, for every $u_i$ not already in $B'$, the set $B'\cup \{u_i\}$ is algebraically dependent. It follows from the Jacobian criterion (Proposition~\ref{prop:jac}) that the row vector $J_{\bf x}(u_i)$ lies in the $\RR({\bf x)}$-span of $\mathcal{B} \defeq \{J_{\bf x}(b)\}_{b \in B'}$. This is also trivially true for $J_{\bf x}(u_i)$ if $u_i$ does lie in $B'$; thus it holds for all $u_i$ with $1\leq i\leq s$. But this means that $J_{\bf x}(v) = \sum_{i=1}^s \alpha_i J_{\bf x}(u_i)$ lies in the $\RR({\bf x})$-span of $\mathcal{B}$. By the Jacobian criterion this contradicts the algebraic independence of $B' \cup \{v\}$.
\end{proof}

\begin{proof}[Proof of Theorem~\ref{thm:alg}]
\leavevmode 
Let $\{u_1,\ldots,u_s\}$ be a basis (or spanning set) for $U$. From above we have that the transcendence degree of $U$ is the row rank of the Jacobian $J_{\bf x}(u_1,\ldots,u_s)$ evaluated at a generic point ${\bf x}$. A transcendence basis for $U$ is the set of $u_i$ corresponding to a maximal linearly independent set of rows.

We can use the following simple greedy algorithm to construct a transcendence basis. As input, receive a list of polynomials $\{u_1,\ldots,u_s\}$. Initialize $I = \emptyset$. For $i = 1, \dots, s$, add $\{u_i\}$ to $I$ if $I \cup \{u_i\}$ is algebraically independent, and do nothing otherwise. (Note that this condition can be efficiently tested by Corollary~\ref{cor:subs}.) Output the resulting set $I$.

We now show correctness. Let $I_i$ be the set after item $u_i$ has been considered (and possibly added), and set $I_0 = \emptyset$. It suffices to show that for each $i \in \{0, \dots, s\}$, $I_i$ is a maximal independent subset of $\{u_1, \dots, u_i\}$. We proceed by induction. The claim is vacuously true when $i=0$. Assume it holds for $i-1$. If $I_i$ is not a maximal independent subset of $\{u_1, \dots, u_i\}$, then there exists an independent set $J \subseteq \{u_1, \dots, u_i\}$ with $|J| > |I|$, so by the exchange property (Proposition~\ref{prop:exchange}) there exists a $u_j$ with $j \leq i$ such that $u_j  \notin I_i$ and $I_i \cup \{u_j\}$ is independent. In particular, the subset $I_{j-1} \cup \{u_j\}$ of $I_i \cup \{u_j\}$ is independent. But the fact that $u_j$ was not added at the $(j-1)$th step implies that $I_{j-1} \cup \{u_j\}$ is not independent, a contradiction. So $I_i$ is indeed maximal.

We obtain that $I = I_s$ is a maximal independent subset of $\{u_1, \dots, u_s\}$, and hence by Lemma~\ref{lem:basis} a transcendence basis of $U$.
\end{proof}

\subsection{Generic unique recovery}
\label{sec:generic-unique}

If the subspace $U$ generates $\RR(X)$, the field of fractions of $\RxG$, as a field then $\RR(Y)=\RR(X)$, i.e. $[\RR(X):\RR(Y)]=1$. Thus the following is immediate from Theorem~\ref{thm:bound-list}.

\begin{corollary}[generic unique recovery]\label{cor:generic-unique}
If $U$ generates the field of fractions of $\RxG$ (as a field), then there exists a set $S \subseteq V$ of full measure such that if $\theta \in S$, then $U$ resolves $\theta$ uniquely.\qed
\end{corollary}

\noindent The intuition here is that we want to be able to learn every invariant polynomial by adding, multiplying, and dividing polynomials from $U$ (and scalars from $\RR$). We need $\theta$ to be generic so that we never divide by zero in the process. 

\begin{remark}
Corollary~\ref{cor:generic-unique} is in essence a special case of a celebrated theorem of Rosenlicht. Recall that $G$, because it is compact and acts continuously and faithfully on a finite-dimensional real vector space, is necessarily a linear algebraic group defined over $\RR$. {\em Rosenlicht's theorem} \cite[Theorem~2]{rosenlicht} asserts that, for any algebraic group acting on any irreducible variety, a finite generating set for the field of invariant rational functions (exists and) separates the orbits on a nonempty $G$-stable Zariski-open subset of the variety. A priori, applying this in the present context would seem to require a stronger hypothesis than that of Corollary~\ref{cor:generic-unique}, because an invariant rational function need not, at the generality of algebraic groups, be a quotient of {\em invariant} polynomials; i.e., $\RR(X)=\operatorname{Frac}\RxG$ might be properly contained in the field $\RR(\mathbf x)^G$ of invariant rational functions. Thus field generators for $\RR(X)$ may not be enough to generate $\RR(\mathbf x)^G$. However, because $G$ is compact, we do have $\RR(\mathbf{x})^G = \operatorname{Frac} \RxG = \RR(X)$ in the present situation; this is established above in the proof of Proposition~\ref{prop:trdeg-general}.  Thus if $U$ generates $\RR(X)=\operatorname{Frac}\RxG$ as a field, Rosenlicht's theorem promises a nonempty Zariski-open subset $S\subset V$ on which any two distinct orbits are distinguished by $U$.\footnote{This is {\em almost} Corollary~\ref{cor:generic-unique}; to extract the full strength of the statement, one would need to do some bookkeeping to guarantee that $S$ can be chosen so that it itself is distinguished from its complement by $U$. This is likely possible but we have not pursued it.}
\end{remark}

\begin{remark}\label{rmk:field-gen-not-nec}
The hypothesis of Corollary~\ref{cor:generic-unique} is sufficient but not necessary for generic unique recovery, as can be seen from the example in Setup~\ref{rmk:pathologies}. An even simpler example is the case that $G$ is the trivial group and $V=\RR^1$ (so $\RxG=\Rx=\RR[x]$), and $U$ is spanned by $x^3$. The function $x^3$ distinguishes every real number even though it does not generate the fraction field $\RR(x)$ over $\RR$. This non-necessity of field generation for generic unique recovery is an artifact of working over a non-algebraically closed field; over $\CC$, generic unique recovery would imply field generation \cite[Lemma~2.1]{pop-vin}.
\end{remark}

\begin{theorem} \label{thm:alg-ext-deg}
For a finite-dimensional subspace $U \subseteq \RxG$, there is an algorithm to compute the degree of the field extension from Theorem~\ref{thm:bound-list}. As input, the algorithm requires a basis for $U$ and the ability to compute the Reynolds operator (Definition~\ref{def:reynolds}).
\end{theorem}

\noindent We give the algorithm and proof in Section~\ref{sec:compute-degree}. The algorithm uses Gr\"obner bases and is often inefficient to run in practice.

Per Remark~\ref{rmk:field-gen-not-nec}, the criterion of Corollary~\ref{cor:generic-unique} is sufficient but not necessary for generic unique recovery. At the cost of an even more inefficient algorithm based on cylindrical algebraic decompositions, we can test generic unique recovery in an if-and-only-if way:

\begin{theorem}\label{thm:generic-unique}
For a finite-dimensional subspace $U \subseteq \mathbb{R}[\textup{\textbf{x}}]^G$, there is an algorithm to decide whether or not $U$ satisfies generic unique recovery. 
\end{theorem}
\begin{proof}
Fix a basis 
$\{u_1, \ldots, u_s\}$ for $U$.
Let $F : V \rightarrow \mathbb{R}^m$ be the polynomial map evaluating $\theta \in V$ on the basis, i.e. $F(\theta) = (u_1(\theta), \ldots, u_s(\theta))$.
The set of signals which are uniquely recovered by $U$ may be written \nolinebreak as:
\begin{equation} \label{eq:def-W}
    W = \left\{ \theta \in V : \forall \theta' \in V~\text{s.t.}~ F(\theta') = F(\theta)~\exists g \in G~\text{s.t.}~\theta' = g \cdot \theta \right\}.
\end{equation}
We seek an algorithm to decide whether or not $W$ has full Lebesgue-measure in $V$.

Note that $W$ is a semi-algebraic subset of $V$ \cite[Definition~2.1.3]{bochnak2013real}. 
Here semi-algebraicity is due to \cite[Proposition~2.2.4]{bochnak2013real} and the fact that $g \in G$ is characterized by polynomial constraints on entries of $g$.
However in general, a semi-algebraic set has full Lebesgue-measure if and only if it is dense with respect to the Euclidean topology, since any semi-algebraic set is a union of finitely many smooth submanifolds \cite[Proposition~9.1.8]{bochnak2013real}. %The cited work says they are unions of finitely many Nash manifolds; Nash manifolds are smooth. 
Thus it suffices to decide whether $V = \overline{W}$ holds, where overline is Euclidean closure.

By definition, 
    $\overline{W} = \left\{ \theta \in V : \forall \varepsilon > 0~\exists \theta' \in W~\text{with}~\|\theta' - \theta\| \leq \epsilon \right\}.$
Inserting \eqref{eq:def-W}, this reads
\begin{equation*} \label{eq:generic-unique-real}
    \overline{W} = \left\{ \theta \in V : \forall \varepsilon > 0~\exists \theta'\in V~\text{with}~\|\theta' - \theta\| \leq \varepsilon~\text{s.t.}~\forall \theta'' \in V~\text{s.t.}~F(\theta'')=F(\theta')~\exists g \in G~\text{s.t.}~\theta''=g \cdot \theta'\right\},
\end{equation*}
whence it is clear ``$V = \overline{W}$" may be expressed as a sentence over the field $\mathbb{R}$ \cite[Section~2.3]{basu-book}.
Consequently, we can decide whether or not ``$V = \overline{W}$" holds by running \cite[Algorithm~11.14]{basu-book}, where the input is  the prenex normal form of said sentence.  
This proves Theorem~\ref{thm:generic-unique}.
\end{proof}

\subsection{Generic unique de-mixing}

For list recovery in heterogeneous problems, there is a shortcut for bounding the size of the list.  Provided an efficiently verifiable condition holds, the list size for a heterogenous problem relates simply to the list size for the corresponding homogeneous ($K=1$) problem.

In what follows, we use $F_G$ to refer to the fraction field of $\RxG$ in order to make the dependence on $G$ explicit. (The dependence on $V$ remains implicit.)

\begin{theorem} \label{thm:hessian}
Let $\tilde{G}$ be a compact group acting linearly and continuously on $\tilde{V}$, and let $G = \tilde{G}^{K} \rtimes S_{K}$ act on $V = \tilde{V}^{\oplus K} \oplus \overline{\Delta}_{K}$ as in Definition~\ref{def:het-action}.  Suppose $U \subseteq \RxG$ is a finite-dimensional subspace of invariant polynomials, corresponding to a subspace $\tilde{U} \subseteq \RR[\tilde{{\bf x}}]^{\tilde{G}}$ in the following sense: ${\bf x} = ({\bf x}^{(1)}, \ldots, {\bf x}^{(K)}, w_{1}, \ldots, w_{K})$ and $U = \{ \sum_{i=1}^{K} w_{i} f({\bf x}^{(i)}) \,\, | \,\, f(\tilde{\bf x}) \in \tilde{U} \}$.  Assume that  $\textup{trdeg}(U) = \textup{trdeg}(\RR[{\bf x}]^{G})$, so that generic list recovery holds in the heterogeneous problem.  Furthermore, assume that $\textup{trdeg}(U) < \dim(\tilde U)$.  Then, there exists an efficient algorithm taking as input the integer $K$ together with a basis for $\tilde{U}$ and returning as output ``pass" or ``fail" such that 
when the test is passed, we have $[F_G : \RR(U)] = [F_{\tilde{G}} : \RR(\tilde{U})]^{K}$.  In particular, when the Hessian test is passed and $[F_{\tilde{G}} : \RR(\tilde{U})]=1$ (so that we have generic unique recovery for the homogeneous problem), then we have generic unique recovery for the heterogeneous problem.
\end{theorem}

\noindent We call the algorithm in Theorem~\ref{thm:hessian} the \emph{Hessian test}, because it involves calculating the rank of a certain Hessian matrix.
It is based on algebraic geometric methods developed for studying symmetric tensor decompositions \cite{COV2}.
Geometrically, the Hessian test bounds the dimension of the contact locus of an associated secant variety. 
See Section~\ref{sec:hessian} for more discussion and the proof of Theorem~\ref{thm:hessian}.

\begin{definition}\label{def:demix}
Passing of the Hessian test implies the following property, which we call \textit{generic unique de-mixing}.  For simplicity, assume $U = U^{T}_{\leq d}$.  Given
$(\Theta, \overline{w}) \in \tilde{V}^{\oplus K} \oplus \overline{\Delta}_{K}$ with $\Theta = (\theta_{1}, \ldots, \theta_{K})$ and $\overline{w} = (\overline{w}_{1}, \ldots, \overline{w}_K) = (w_1 - 1/K,\ldots,w_K - 1/K)$, one may form the concatenation of moment tensors,  $\big{(}T_{1}(\Theta, w) \oplus \ldots \oplus T_{d}(\Theta, w)\big{)}$. This list of moments for the heterogeneous problem equals a weighted average of $K$ lists of moments for the homogeneous problem:
$$w_{1}\big{(}T_{1}(\theta_{1}) \oplus \ldots \oplus T_{d}(\theta_{1})\big{)} + \ldots + w_{K}\big{(}T_{1}(\theta_{K}) \oplus \ldots \oplus T_{d}(\theta_{K})\big{)}.$$
Generic unique de-mixing entails that this splitting of heterogeneous moments into homogeneous moments is unique (up to permutation of summands), for $(\Theta, \overline{w})$ in a nonempty Zariski-open subset of $\tilde{V}^{K} \times \overline{\Delta}_{K-1}$.
\end{definition}

\begin{remark}
Our proof shows whenever the Hessian test is passed, generic unique de-mixing holds.  The proof leverages the observation that generic unique de-mixing is akin to uniqueness for tensor decomposition (where concatenated moments for the single signal problem play the role of rank 1 tensors). In terms of orbit recovery, bounding list sizes for the homogeneous problem then remains.
\end{remark}

\subsection{Worst-case unique and list recovery}\label{sec:worst-unique}

We conclude Section~\ref{sec:algebraic} by considering the worst-case list and unique recovery problems. For readers without an algebraic background, the main ideas from commutative algebra used in this section are collected in Appendix~\ref{app:comm-alg}.

By Setup~\ref{rmk:separation-in-terms-of-AG}, a sufficient criterion for worst-case unique recovery is that $U$ generates $\RxG$ as an $\RR$-algebra, for in this case, the map $\varphi_\RR$ is the identity (and therefore injective). This criterion is computationally accessible (although not efficiently, as it relies on Gr\"{o}bner bases): compute a finite generating set for $\RxG$ using one of the known algorithms (for example see Appendix~\ref{app:generator-alg}); then test whether each of these generators is in the algebra generated by $U$ (a standard Gr\"{o}bner basis test; see Appendix~\ref{app:grobner}).

The algorithm given in Appendix~\ref{app:generator-alg} for computing generators for $\RxG$ is originally due to Derksen \cite{derksen1999}. It, and the algorithms discussed below, require a presentation of $G$ as an algebraic variety, and a matrix of polynomial functions on $G$ describing its action on $V$. An example illustrating these data is given in Appendix~\ref{app:generator-alg}.

Generation of $\RxG$ by $U$ is not necessary for worst-case unique recovery, however. It is possible for a subspace $U$ of $\RxG$ to uniquely resolve every orbit without generating $\RxG$. A subset of $\RxG$ that can distinguish any pair of orbits distinguished by the entire ring $\RxG$ is called a \emph{separating set} for $\RxG$, and the subalgebra it generates a \emph{separating algebra}---see \cite[Section~2.4]{dk-book} and \cite{kemper-sep}. In the language of Setup~\ref{rmk:AGsetup}, $U$ is a separating set if $\varphi_\RR:X(\RR)\rightarrow Y(\RR)$ is an injection when restricted to the image of $\pi_\RR$.\footnote{There is a minor ambiguity in the definition of a separating set coming from the fact that most of the work on separating sets has presumed that the ground field is algebraically closed, so that it is not unreasonable to read ``$U$ is a separating set" to mean that $\varphi_\CC$ (rather than $\varphi_\RR$) is injective on the image of $\varphi_\CC$ (rather than $\pi_\RR$). Our usage follows Dufresne \cite{dufresne}, who distinguishes {\em geometric separating sets}, i.e., sets which remain separating after base change to the algebraic closure, from sets that are merely separating. In the case of trivial action on $\RR^1$, $x^3$ is separating, in this language, but not geometrically separating.} Since in the present case of compact $G$ the map $\varphi:V(\RR)\rightarrow X(\RR)$ uniquely resolves every orbit, injectivity of $\varphi:X(\RR)\rightarrow Y(\RR)$ on the image of $\pi_\RR:V(\RR)\rightarrow X(\RR)$ is sufficient and necessary for $U$ to uniquely resolve all orbits of $G$ on $V$, i.e., $U$ resolves every orbit if and only if it is a separating set. The notions of generating and separating sets do not always coincide, as illustrated by Setup~\ref{rmk:pathologies} and Remark~\ref{rmk:field-gen-not-nec} above, where $x^4$ and $x^3$ respectively are separating but not generating. Example~2.4.2 in \cite{dk-book} shows that they need not coincide even when working over an algebraically closed ground field.

We pull out the take-home message of the previous paragraph. Essentially it restates the definition of a separating set:

\begin{proposition}[worst-case unique recovery]\label{prop:worst-case-unique}
Let $U \subseteq \RxG$ be a finite-dimensional subspace. Then $U$ resolves every $\theta \in V$ if and only if it is a separating set for $\RxG$.\qed
\end{proposition}

\noindent Using ideas from computational invariant theory and computational real algebraic geometry, this criterion is computationally accessible given an oracle to find a strictly feasible point (i.e., a relative interior point of the feasible region) for an SDP in exact arithmetic.

\begin{proposition} \label{prop:alg-separate}
Given an oracle to produce, in exact arithmetic, a strictly feasible point of a semidefinite program specified in exact arithmetic, there is an algorithm to test whether or not $U$ is a separating set for $\RxG$. As input, the algorithm requires a basis for $U$, a presentation of $G$ as an algebraic variety, a polynomial description of $G$'s action on $V$, and the ability to compute the Reynolds operator (Definition~\ref{def:reynolds}).
\end{proposition}

\noindent The proof is given in Section~\ref{sec:separating-alg}. The algorithm is not efficient since it relies on Gr\"{o}bner bases. It is based on the main idea of an algorithm of Kemper for computing separating invariants over an algebraically closed field \cite{kemper2003} (and see the beautiful exposition in \cite{dk-book}, where it appears as Algorithm~4.9.4). It is simpler than Kemper's algorithm: because $G$ is {\em linearly} reductive (i.e., it has a Reynolds operator---see Definition~\ref{def:reynolds}), we can take advantage of Derksen's algorithm \cite{derksen1999} for generators of $\RxG$ (see Appendix~\ref{app:generator-alg}) as a subroutine. Also, because we are working over $\RR$ rather than the algebraically closed field $\CC$, we use the real radical membership test of Reid, Wang, and Wolkowicz  \cite{reid-wang-wolkowicz} in place of an ordinary radical membership test. This is the step that requires the SDP strictly-feasible point oracle. It actually may be within the reach of current computational methods, see Remark~\ref{rmk:sdp-numerics}, but this is beyond our scope.

\begin{remark}\label{rmk:apriori-bounds-sep}
In an orbit recovery problem without a projection $\Pi$, the subspace $U_{\leq d}^T$ is nothing but $\RxG_{\leq d}$, the space of invariants of degree less than or equal to $d$. In this situation, some a priori information is known about when $\RxG_{\leq d}$ is either a generating or a separating set for $\RxG$.

If $G$ is a finite group, bounding the degree $d$ in which $\RxG$ is generated is an old and well-studied problem. It is known that $\RxG$ has a generating set for which all elements have degree at most $|G|$ (this is \emph{Noether's degree bound}, originally proven in \cite{noether}; see Theorem~2.1.4 in \cite{sturmfels} for a modern exposition). It follows from Setup~\ref{rmk:separation-in-terms-of-AG} that $\RxG_{\le |G|}$ resolves every $\theta \in V$. Recall from Section~\ref{sec:prob-statement} that this is tight for MRA: degree $|G|$ is necessary for worst-case signals. However, this is the only scenario in which it is tight: if $G$ is a finite group that is not cyclic, $\RxG$ is already generated below degree $|G|$ \cite{schmid}; in fact, it is generated by elements of degree at most $\frac{3}{4}|G|$ \cite{domokos-hegedus}. If $G$ does not have an index $\leq 2$ cyclic subgroup and is not on a short list of small exceptional cases, then $\RxG$ is even generated by elements of degree at most $\frac{1}{2}|G|$ \cite{cziszter}.

If $G$ is not finite, some a priori bounds are still known---see \cite[Section~4.7]{dk-book}.

The problem of determining $d$ such that $\RxG_{\leq d}$ is a separating set has more recently begun to attract research attention; see \cite{kemper-sep, kohlskraft, domokos}.\footnote{A related inquiry is undertaken by Kadish \cite{kadish}, where the ring of polynomial functions is augmented by ``quasi-inverses". This destroys the notion of the degree of a function, but Kadish replaces it with the {\em length of a program to determine the function}, and he finds bounds on the lengths required for functions to be separating. Another related inquiry is found in \cite{cahill2022group}, which separates orbits using certain functions which are non-polynomial but have other good properties (convex, Lipschitz), and gives bounds on the complexity of computing them in certain special cases (see Section~3.2).} If $G$ is finite abelian and non-cyclic, then \cite[Corollary~3.11]{domokos} shows that, with a small number of exceptions, separation can be achieved in lower degree than that required for generation of $\RxG$.\footnote{The result in \cite{domokos} is stated over an algebraically closed field rather than $\mathbb{R}$, but if the action is defined over $\mathbb{R}$ then separation over $\mathbb{C}$ implies separation over $\mathbb{R}$, as can be seen from the diagram in Setup~\ref{rmk:X(C)etc} above (and see also \cite[Theorem~2.2]{dufresne}).}
\end{remark}

We also give a sufficient algebraic condition for worst-case list recovery:

\begin{proposition}[worst-case list recovery] \label{prop:worst-list}
Let $U \subseteq \RxG$ be a subspace with finite homogeneous basis $\{f_1,\ldots,f_m\}$. If $\RxG$ is finitely generated as a $\RR[f_1,\ldots,f_m]$-module, then $U$ list-resolves every $\theta \in V$.
\end{proposition}

\begin{proof}
The hypothesis is that $\varphi:X\rightarrow Y$ is a finite map. Thus $\varphi_\CC:X(\CC)\rightarrow Y(\CC)$ has finite fibers (e.g., \cite[Section~5.3]{shaf}), and it follows from the diagram in Setup~\ref{rmk:X(C)etc} that $\varphi_\RR:X(\RR)\rightarrow Y(\RR)$ does too. We conclude by the discussion in Setup~\ref{rmk:separation-in-terms-of-AG}.
\end{proof}

The criterion of Proposition~\ref{prop:worst-list} is computationally accessible:

\begin{proposition}\label{prop:alg-worst-list}
There is an algorithm to test if $\RxG$ is finitely generated as a module over $\RR[U]$. As input, the algorithm requires a basis for $U$, a presentation of $G$ as an algebraic variety, a polynomial description of $G$'s action on $V$, and the ability to compute the Reynolds operator (Definition~\ref{def:reynolds}).
\end{proposition}

\noindent The algorithm is similar to that of Proposition~\ref{prop:alg-separate}, although it is somewhat simpler, so we include it here. (The relationship is discussed in Section~\ref{sec:separating-alg}.) Like that algorithm, it is not efficient as it involves Gr\"{o}bner bases.

\begin{proof}
Because $U$ has a homogeneous basis, we can apply the following fact: $\RxG$ is module-finite over $\RR[U]$ if and only if the radical of the ideal generated by $U$ in $\Rx$ contains all homogeneous invariants of positive degree. (This follows from, e.g., \cite[Lemma~2.5.5 and Remark~4.7.2]{dk-book} in view of Hilbert's Nullstellensatz.) This condition can be tested by computing a homogeneous algebra generating set for $\RxG$ (for example using Derksen's method; see Appendix~\ref{app:generator-alg}), and then testing membership of each of these generators in the radical of the ideal generated by a basis for $U$. Testing membership in a radical is a standard Gr\"{o}bner basis procedure \cite[Chapter~4, Section~2, Proposition~8]{iva}; there is a slightly simplified version when everything in sight is homogeneous \cite[Lemma~4.9.3]{dk-book}.
\end{proof}

\noindent The radical membership tests consist of checking whether $1$ is in a certain ideal explicitly constructed from the input data, using the ideal membership test described in Appendix~\ref{app:grobner}.

\begin{remark}
The condition in Proposition~\ref{prop:worst-list} is sufficient but not quite necessary. Worst-case list recovery holds when the restriction of $\varphi_\RR$ to the image of $\pi_\RR$ has finite fibers. A priori, one can imagine several ways this might happen without $\varphi$ being finite: perhaps (i) $\varphi$ is only quasifinite (in other words, $\varphi_\CC$ has finite fibers, but not because $\RxG$ is module-finite over $\RR[U]$); or perhaps (ii) $\varphi_\CC$ has one or more infinite fibers, but they all lie outside of the image of $\pi_\RR$; or lastly, perhaps (iii) $\varphi_\CC$ has some infinite fibers that intersect the image of $\pi_\RR$, but they do so in only finitely many points. Because we are assuming $U$ is spanned by homogeneous invariants, the first two of these hypotheticals are actually impossible; but the third does happen. 

For example, take $V=\RR^2$ with $G$ the trivial group. Then $\Rx=\RxG=\RR[x,y]$, $V=X$ is a plane, and $\pi$ is the identity map. Let $U$ be the subspace of $\RR[x,y]$ generated by $x^2+y^2$ and $(x^2+y^2)x$. These are algebraically independent, so that $\RR[U]$ is isomorphic to a polynomial algebra $\RR[s,t]$ with $s=x^2+y^2$ and $t=(x^2+y^2)x$, and $Y$ is also a plane. The map $\varphi_\CC:X(\CC)\rightarrow Y(\CC)$ is given on points of $\CC^2$ by
\[
(x,y) \mapsto (x^2+y^2, (x^2+y^2)x).
\]
Over any point $(\alpha,\beta)\in \CC^2$ with $\alpha \neq 0$, the fiber of $\varphi_\CC$ is the solution set of the system  $x^2+y^2=\alpha$, $x(x^2+y^2)=\beta$, which is $x=\beta/\alpha$, $y=\pm \sqrt{\alpha - (\beta/\alpha)^2}$; this has cardinality 1 or 2. If $\alpha=0$ but $\beta\neq 0$, the fiber is empty. But if $\alpha=\beta=0$, the fiber is the solution set of $x^2+y^2=0$. This is infinite, as it is a pair of complex lines $\{x=\pm iy\}$. However, its only real point is $x=y=0$. Thus $\varphi_\RR$ has finite fibers in this situation, although $\varphi_\CC$ does not.

To see that hypotheticals (i) and (ii) are impossible on the assumption that $U$ is graded, we use the fact that, on this hypothesis, module-finiteness of $\RxG$ over $\RR[U]$ is equivalent to the condition that the ideal $\langle U \rangle_{\RxG}$ generated by $U$ in $\RxG$ has height equal to the Krull dimension of $\RxG$.\footnote{This is true because in one direction, if $\operatorname{ht} \langle U\rangle_{\RxG} = \dim_{\mathrm{Krull}} \RxG$, then the quotient ring $\RxG/\langle U\rangle_{\RxG}$ is a graded $\RR$-algebra of Krull dimension zero, and therefore finite-dimensional as a graded $\RR$-vector space. One can take a graded $\RR$-basis for $\RxG/\langle U\rangle_{\RxG}$ (which is therefore finite) and find homogeneous lifts of each member in $\RxG$; these will generate $\RxG$ as an $\RR[U]$-module by the graded Nakayama lemma. In the other direction, suppose $\RxG$ is module-finite over $\RR[U]$. Any prime $\mathfrak{p}$ of $\RxG$ that is minimal over $\langle U\rangle_{\RxG}$ pulls back to the ideal $\langle U\rangle_{\RR[U]}$ generated by $U$ in $\RR[U]$, which is evidently maximal as the quotient by it is $\RR$; because $\RxG$ is module-finite over $\RR[U]$, $\mathfrak{p}$ is therefore itself maximal. As $\RxG$ is an integral domain finitely generated over a field, all its maximal ideals have same height, and this is the Krull dimension of $\RxG$.} Because the latter is an integral domain which is a finitely generated algebra over a field, its Krull dimension is also the height of all of its maximal primes. In particular, if $\RxG$ is not module-finite over $\RR[U]$, then the quotient ring $\RxG/\langle U\rangle$ has positive Krull dimension. But this quotient ring is isomorphically the tensor product
\[
\RxG \otimes_{\RR[U]}\RR[U]/\langle U\rangle_{\RR[U]},
\]
which is precisely the coordinate ring of the fiber of $\varphi$ over the point of $Y$ corresponding to $0$ in $V$. As $0\in V(\RR)$, this point always lies in the image of $\pi_\RR$. Therefore, if $\varphi$ fails to be a finite map, then $\varphi_\CC$ has a positive-dimensional (and therefore infinite) fiber that contains the image of $0\in V(\RR)$ and thus intersects the image of $\varphi_\RR$. This rules out (i) and (ii).
\end{remark}

\begin{remark}
As in Remark~\ref{rmk:apriori-bounds-sep}, when working in the setup without a projection $\Pi$, there exist a priori bounds on the $d$ necessary so that $U=U_{\leq d}^T = \RxG_{\leq d}$ satisfies the criterion of Proposition~\ref{prop:worst-list}. These bounds require as input information about $G$ and its action on $V$.

Here is one such bound. The orbits of $G$ on the vector space $V$ are algebraic varieties defined over $\RR$ \cite[Chapter~3, \S~4, Theorem~4]{on-vin}. Given an orbit $O$, it has a {\em degree}: this can be characterized as the cardinality of the intersection of its set of $\CC$-points, $O(\CC)$, with a suitably generic complex affine subspace of $V(\CC)=\CC^p$ that is of dimension complementary to $O$. (For example, the nontrivial orbits for the canonical action of $G = \mathrm{O}(3)$ on $V=\RR^3$ are the spheres $x^2+y^2+z^2 = R$. These are all degree $2$ because a generic complex line in $\CC^3$ intersects $\{x^2+y^2+z^2 = R\}\subset\CC^3$ in $2$ points for any $R>0$.) Then if $d$ is the maximum degree of an orbit of $G$ on $V$, then $\RxG$ is module-finite over $\RR[U_{\leq d}^T]$, the subalgebra generated by invariants of degree $d$ \cite[Proposition~4.7.12]{dk-book}. 

This and other bounds are discussed in Section~4.7.1 of \cite{dk-book}, which is based on work of Popov \cite{popov1982} and Hiss \cite{hiss}.
\end{remark}

\section{Examples}
\label{sec:examples}

In this section, we work out the main motivating examples, determining the degree at which recovery is possible by the methods of Section~\ref{sec:algebraic}. The examples involve representations of the cyclic and special orthogonal groups.  Mostly, we focus on generic list recovery since our algorithms for certifying the other recovery criteria are unfortunately too slow even for small examples.  We obtain several recovery theorems for orbit recovery problems such as MRA and cryo-EM within finite ranges of parameters where we have verified the Jacobian criterion using a computer, and beyond these parameter ranges we state conjectural patterns.

The following themes emerge in the examples studied in this section. First, we see that many problems are possible at degree $3$, which is promising from a practical standpoint as the problems require a (hopefully) manageable number of samples in practice. Second, we do not encounter any unexpected algebraic dependencies among invariants, and so we are able to show that heuristic parameter-counting arguments are correct. In particular, we find that if there are enough linearly independent invariants, there are also enough algebraically independent invariants.

\subsection{Learning a bag of numbers} \label{ex:bag}

Let $G$ be the symmetric group $S_p$ acting on $V = \RR^p$ by permutation matrices. The invariant ring consists of the symmetric polynomials, which are generated by the elementary symmetric polynomials $e_1,\ldots,e_p$ where $e_i$ has degree $i$. Worst-case unique recovery is possible at degree $p$ since $\RxG_{\le p}$ generates the full invariant ring. Furthermore, degree $p$ is actually required, even for generic list recovery. This is because any invariant of degree $\le p-1$ can be expressed as a polynomial in $e_1,\ldots,e_{p-1}$ and thus $\trdeg(\RxG_{\le p-1}) = p-1$. So this problem has a steep sample complexity of order $\sigma^{2p}$. 

An efficient signal recovery algorithm at degree $d=p$ is obtained as follows.
Let $T$ be a formal variable and consider the univariate polynomial $f = T^p - e_1 T^{p-1} + e_2 T^{p-2} - e_3 T^{p-3} + \ldots + (-1)^{p} e_p$.  The multiset of its roots is the bag of numbers.  
Hence we may compute the signal by computing all roots of $f$, e.g.\ by calculating the eigenvalues of the associated companion matrix \cite{sturmfels2002solving}.

\subsection{Learning a rigid body}\label{sec:rigid-body}

Let $G$ be the rotation group $\mathrm{SO}(p)$ acting on the matrix space $\RR^{p \times m}$ by left multiplication. We imagine the columns of our matrix as vertices defining a rigid body; thus we observe random rotations of this rigid body (with vertices labeled) plus noise. Let $U \in \RR^{p \times m}$ be such a matrix signal. With $O(\sigma^4)$ samples, we can estimate the degree-2 Gram matrix $U^\top U$; taking a Cholesky factorization, we recover $U$ up to left action by an element of the larger group $\mathrm{O}(p)$. Thus we recover the rigid body up to a reflection ambiguity, demonstrating list recovery at degree 2 (with a list of size 2). See~\cite{procrustes} for a more detailed analysis of this algorithm. Surprisingly, assuming $m \geq p$, we do not uniquely resolve a generic signal until degree $p$, where with $O(\sigma^{2p})$ samples we can estimate a $p \times p$ minor of $U$, which is a degree-$p$ invariant that changes sign under reflection.

The impossibility of unique recovery until degree $p$ is a consequence of the ``first fundamental theorem'' for the special orthogonal group $\mathrm{SO}(p)$, which asserts that the invariant ring is generated by the entries of the Gram matrix $U^\top U$ together with the $p \times p$ minors of $U$ (see for instance \cite{kac-notes}); thus the invariants of degree $3,\ldots,p-1$ carry no information beyond the information carried by the degree-2 invariants.

\subsection{Regular representation of a finite group}\label{sec:reg_rep}

Let $G$ be a finite group. Let $V$ be the regular representation of $G$ over $\RR$, i.e.,\ $V = \mathbb{R}^{|G|}$ with basis indexed by the elements of $G$, and the action of $G$ permutes the basis elements $v_g$ according to group multiplication: $h \cdot v_g = v_{hg}$. In other words, the signal can be thought of as a real-valued function on $G$. Note that for the cyclic group $G = \mathbb{Z}/p$, this is precisely the MRA problem. It is known \cite{complete-bispectrum} that for this setup, the \emph{triple correlation} (a collection of degree-$3$ invariants) is sufficient to resolve a generic signal, and thus generic unique recovery is possible at degree $3$. We give two alternative proofs of this fact, as we believe the proof techniques are instructive, and also each proof yields some new information. In Section~\ref{sec:galois}, we restrict to the special case where $G$ is abelian and use Galois theory to show that the invariant field is generated at degree 3. In Section~\ref{sec:jennrich} we use Jennrich's algorithm for tensor decomposition to show that generic unique recovery is not only possible, but achievable with an efficient algorithm.

\subsubsection{Field generation via Galois theory}
\label{sec:galois}

\begin{theorem}\label{thm:reg-rep-galois}
Let $G$ be a finite abelian group and let $V$ be the regular representation of $G$ over $\RR$. Then $\RR(U^T_{\le 3})$ is equal to the field of fractions of $\RxG$.
\end{theorem}

\noindent Recall that the equality of the two fields above implies generic unique recovery at degree 3 by Corollary~\ref{cor:generic-unique}. The proof uses Galois theory and Pontryagin duality; statements of and references for the needed theorems are in Appendices~\ref{app:fields} and \ref{app:pontryagin}.

\begin{proof}
First we argue it is sufficient to show $\CC(U^T_{\le 3}) = \CC(\RxG)$.  Suppose this holds, and $f$ is a nonzero element of $\RR(\RxG)$. Let $f = P/Q$ with $P,Q\in \RxG$. Since also $f\in\CC(\RxG) = \CC(U^T_{\le 3})$, we also have $f=L/M$ where $L,M$ are complex polynomial expressions in the elements of the real vector space $U^T_{\le 3}$. We have $L = L_{\text{re}}+iL_{\text{im}}$, and similarly $M = M_{\text{re}} + iM_{\text{im}}$, with $L_{\text{re}},L_{\text{im}},M_{\text{re}}, M_{\text{im}}\in \RR[U^T_{\le 3}]$. Without loss of generality, $M_{\text{re}}\neq 0$ (if not, multiply numerator and denominator of $L/M$ by $i$). Then
\[
P(M_{\text{re}}+i M_{\text{im}}) = Q(L_{\text{re}}+iL_{\text{im}})
\]
in the field $\CC(x)$. We find that $PM_{\text{re}} = QL_{\text{re}}$ by equating real and imaginary parts. Thus $f=P/Q = L_{\text{re}}/M_{\text{re}}$ expresses $f$ as an element of $\RR(U^T_{\le 3})$. So $\RR(\RxG)\subset\RR(U^T_{\le 3})$, and the opposite inclusion is already clear. 

Thus our goal is to show that $\CC(U_{\le 3}^T)$ exhausts $\CC(\RxG)$. We view all the action as taking place inside the field $\CC(\mathbf{x})$, of which these are subfields and $\Rx, \RxG$ are subrings. The action of $G$ on $\Rx$ extends in a natural way to an action on $\CC(\mathbf{x})$. Let $p = |G| = \dim(V)$.

Because the regular representation is self-dual, the induced action of $G$ on the degree-1 component $\langle \mathbf{x}\rangle_\CC$ of $\Cx\subset \CC(\mathbf{x})$ is isomorphic to the regular representation. Because $G$ is abelian and $\CC$ is algebraically closed, it follows that $\langle \mathbf{x}\rangle_\CC$ splits into the direct sum of $p$ one-dimensional subspaces, each an eigenspace for a unique one of $G$'s $p$ distinct characters. Let $\hat G$ be the character group of $G$, and for $\chi\in \hat G$, let $y_\chi\in \langle \mathbf{x}\rangle_\CC$ be an eigenvector of $\chi$, i.e., an element such that $g \cdot y_\chi = \chi(g) y_\chi$ for all $g \in G$. We have
\[
\CC(\mathbf{x}) = \CC(\{y_\chi\}_{\chi\in\hat G})
\]
because the $y_\chi$'s form a basis for $\langle\mathbf{x}\rangle_\CC$. For brevity, we abbreviate $\{y_\chi\}_{\chi\in \hat G}$ as $\{y_\chi\}$; the set braces will implicitly index over $\hat G$ in the below.

Observe that for any $\chi \in \hat G$, 
\[
y_\chi y_{\chi^{-1}}
\]
is a degree-2 invariant, and for any $\chi_1,\chi_2 \in \hat G$, 
\[
y_{\chi_1}y_{\chi_2}y_{(\chi_1\chi_2)^{-1}}
\]
is a degree-3 invariant. For each $\chi \in \hat G$ we have
\[
y_\chi^p \prod_{\chi' \in \hat G} y_{\chi'} y_{\chi'^{-1}} =\prod_{\chi' \in \hat G} y_{\chi}y_{\chi'}y_{(\chi\chi')^{-1}},
\]
since multiplication by $\chi$ permutes the elements of $\hat G$. It follows that $y_\chi^p$ lies in the field generated by the invariants of degree $\le 3$. Therefore, $\CC(U^T_{\leq 3})$ contains the field $\CC(\{y_\chi^p\})$.

The field extension $\CC(\{y_\chi\})/\CC(\{y_\chi^p\})$ is Galois. An automorphism $\phi$ of this extension multiplies each $y_\chi$ by a $p$th root of unity, and any choice of a $p$th root of unity for each $y_\chi$ yields an automorphism. Thus, the Galois group is $\Gamma =\prod_{\chi\in \hat G} \mu_p^\times$, where $\mu_p^\times$ denotes the group of $p$th roots of unity in $\CC$. Given an element $\phi\in \Gamma$, denote by $\phi_\chi$ the $p$th root of unity such that $\phi\cdot y_\chi = \phi_\chi y_\chi$.

By the fundamental theorem of Galois theory, intermediate extensions of the above Galois extension are in bijection with subgroups of $\Gamma$. Because $G$ consists of automorphisms of $\CC(\mathbf{x})=\CC(\{y_\chi\})$ under which each $y_\chi^p$ is invariant, $G$ is a subgroup of $\Gamma$. 
The field $\CC(\RxG)$ is its corresponding fixed field under the Galois correspondence: we see this because first, it is evidently equal to $\CC(\CxG)$, and second, since $G$ is finite we have $\CC(\CxG) = \CC(\mathbf{x})^G$ (for example by \cite[Lemma~3.2]{pop-vin}). The field $\CC(U^T_{\le 3})$ is an intermediate extension contained in $\CC(\RxG)$, so its corresponding subgroup, call it $G_1$, contains $G$. We will show it equals $G$, thus proving $\CC(U^T_{\le 3}) = \CC(\RxG)$ as desired. 

We need to show that there are no elements of $G_1$, i.e., elements of $\Gamma$ that fix $U^T_{\le 3}$, beyond those that lie in $G$. Let $\phi \in G_1$. Then $\phi$ fixes each degree-2 invariant $y_\chi y_{\chi^{-1}}$, so we must have 
\[
\phi_{\chi^{-1}} = \phi_\chi^{-1}
\]
for all $\chi$. Because $\phi$ also fixes each degree-3 invariant $y_{\chi_1}y_{\chi_2}y_{(\chi_1\chi_2)^{-1}}$, we must have 
\[
\phi_{\chi_1} \phi_{\chi_2} = \phi_{(\chi_1 \chi_2)^{-1}}^{-1} = \phi_{\chi_1 \chi_2},
\]
for all $\chi_1,\chi_2$, where the second equality uses what we learned from the degree-2 invariants.
In particular, the map
\begin{align*}
\hat G &\to \CC^\times\\ 
\chi &\mapsto \phi_\chi
\end{align*}
is a group homomorphism. Thus, any $\phi\in G_1$ can be identified with a character of $\hat G$. By Pontryagin duality, the group of these is isomorphic to $G$. In particular there are only $p$ such elements $\phi$, so $G$ exhausts $G_1$. This completes the proof.
\end{proof}

\begin{remark}\label{rmk:nonmodular-field-gen}
In the above statement, one can replace $\RR$ with any field of characteristic prime to the group order, and in the proof, one can replace $\CC$ with any field extension of it containing all the $|G|$th roots of unity. (The assumption about the characteristic is so that these roots of unity are distinct.) Essentially the same proof goes through. The argument in the first paragraph with real and imaginary parts in essence only uses the fact that $\CC$ is free as a module over $\RR$ (and therefore $\Cx$ is a free module over $\Rx$); a similar situation holds for any field extension.
\end{remark}

\begin{remark}
The argument here can be seen as an algebraization of the ``frequency marching'' method used in \cite{BenBouMa17} to algorithmically achieve generic unique recovery at degree 3 for MRA. Rewriting $\langle\mathbf{x}\rangle_\CC$ in terms of the $y_\chi$ is essentially passing to the discrete Fourier transform; finding that $y_\chi^p\in \CC(U_{\le 3}^T)$ is analogous to determining the magnitudes of the Fourier coefficients; and the argument that the automorphisms of $\CC(\mathbf{x})$ fixing $U_{\le 3}^T$ biject with characters of $\hat G$ corresponds with the ``frequency-marching" equations that relate the phases of the Fourier coefficients to each other using the bispectrum.

As mentioned above, generic unique recovery at degree 3 is actually true in the more general setting where $V$ is the regular representation over $\RR$ of any finite group; this follows from work of Kakarala \cite{complete-bispectrum},\footnote{Kakarala's result is actually a still more general statement about any compact group, but the regular representation is only finite-dimensional, and so only falls within the scope of the present work, when $G$ is finite.} and see below for a different proof. Kakarala proves generic unique recovery for the regular representation over $\CC$ using functions of degree three that are composed both of polynomials and antiholomorphic polynomials; they become true polynomials only when restricted to real points. Per Remark~\ref{rmk:field-gen-not-nec}, generic unique recovery over $\RR$ does not establish field generation, so a priori, Kakarala's result does not show that $U^T_{\le 3}$ generates $\RR(\RxG)=\RR(X)$ for the regular representation. Thus Theorem~\ref{thm:reg-rep-galois} seems to be new.

However, Kakarala's argument is based on Tannaka-Krein duality, which generalizes Pontryagin duality to nonabelian groups, and uses it in a manner analogous to our use of Pontryagin duality here. This suggests it may be possible to integrate the proof of Theorem~\ref{thm:reg-rep-galois} with Kakarala's method in order to obtain a field generation result that generalizes Theorem~\ref{thm:reg-rep-galois} to all finite groups. It is not trivial because there is no analogy in Kakarala's argument to the step of showing that each $y_\chi^p\in \CC(U_{\le 3}^T)$, which we needed here to establish that $\CC(U_{\le 3}^T)$ is an intermediate field of a Galois field extension with known Galois group. This would be an interesting avenue for further work.
\end{remark}

\subsubsection{Unique recovery via tensor decomposition}
\label{sec:jennrich}

The ideas in this section are inspired by prior work on MRA \cite{PerWeBan17}. In the \emph{symmetric third-order tensor decomposition problem}, we are given a tensor of the form
$$T = \sum_{i=1}^r a_i^{\otimes 3}$$
and the goal is to recover the vectors $a_1,\ldots,a_r \in \mathbb{R}^p$ (up to permutation). It is a classical result that if $\{a_i\}$ are linearly independent then there exists a unique decomposition (i.e.\ $a_1,\ldots,a_r$ are uniquely determined from $T$, up to permutation) and furthermore, it can be found in polynomial time using Jennrich's algorithm (see Chapter 3 of \cite{moitra-notes}).

For an orbit recovery problem over a finite group (without heterogeneity or projection), the third-order moment tensor is
$$T_3(\theta) = \EE_g[(g \cdot \theta)^{\otimes 3}] = \frac{1}{|G|}\sum_{i=1}^r (g_i \cdot \theta)^{\otimes 3}$$
where $|G| = r$ and $G = \{g_1,\ldots,g_r\}$. Thus, if the vectors $(g_i \cdot \theta)$ are linearly independent, tensor decomposition gives unique recovery of (the orbit of) $\theta$. When $V$ is the regular representation, we can show that this is the case:

\begin{theorem}\label{thm:generic-linear-indep}
Let $G$ be a finite group and let $V$ be the regular representation of $G$ over $\RR$. Then for a generic signal $\theta$, $(g \cdot \theta)_{g \in G}$ are linearly independent. Thus, generic unique recovery is possible at degree 3.
\end{theorem}
\begin{proof}
Let $\tilde \theta \in V$ be arbitrary and let $\theta$ be a perturbed signal obtained from $\tilde \theta$ by adding a small random value $\eta$ to the first coordinate, say $\eta$ is uniform on $[-\epsilon,\epsilon]$ for some $\epsilon > 0$. For each $i = 1,\ldots,p$, let $g_i$ be the unique group element that permutes the coordinates of the regular representation $V$ in such a way that index 1 maps to index $i$. Let $M$ be the matrix with columns $g_1 \cdot \theta, \ldots, g_p \cdot \theta$ and let $\tilde M$ have columns $g_1 \cdot \tilde \theta, \ldots, g_p \cdot \tilde \theta$. We have $M = \tilde M + \eta I$. Our goal is to show that $M$ is full rank. If this were not the case then there would exist a vector $v$ for which $M v = 0$, i.e.\ $\tilde M v = -\eta v$ and so $-\eta$ is an eigenvalue of $\tilde M$. However, since $\tilde M$ has a finite number of eigenvalues, this occurs with probability zero over the choice of $\eta$.

It follows from the above that for any $\tilde \theta$, the neighborhood $\tilde \theta + [-\epsilon,\epsilon]^p$ contains only a measure zero set of ``bad'' signals $\theta$ for which $(g \cdot \theta)_{g \in G}$ are not linearly independent. Thus the set of all ``bad'' $\theta \in V$ is measure zero.
\end{proof}

The same statement follows immediately for any representation $V$ of $G$ that contains the regular representation as a subrepresentation $V'$. Indeed, in this situation, Maschke's theorem implies there is a $G$-equivariant projection $V\rightarrow V'$, whereupon a generic vector $\theta\in V$ projects to a generic vector $\theta'$ in $V'$, and linear independence of $(g\cdot\theta')_{g\in G}$ (which holds by Theorem~\ref{thm:generic-linear-indep}) then implies linear independence of $(g\cdot\theta)_{g\in G}$.

Note, however, that if $V$ does not contain the regular representation (for example if $|G| > \dim(V)$), then it is impossible for $(g \cdot \theta)_{g \in G}$ to be linearly independent; indeed, if they are linearly independent then they span a copy of the regular representation. (This observation of Matthew Satriano was relayed to us by an anonymous referee.)

In this situation, it may be possible to show generic unique recovery using results on uniqueness of \emph{overcomplete} tensor decomposition (e.g.\ \cite{kruskal-tensor,strassen-tensor,kruskal-concise, COV1}), but these would not yield an efficient algorithm for performing the decomposition.

\begin{remark}
In this section we have shown that given exact access to the third-order moment tensor, Jennrich's algorithm can be used to efficiently recover the signal. As in \cite{PerWeBan17}, a robust analysis of Jennrich's algorithm (based on \cite{gvx}) shows that Jennrich's algorithm recovers a vector close to the truth when given a noisy estimate of the moments.
\end{remark}

\subsubsection{List recovery is impossible at degree 2}

To complement the previous two subsections, we show that list recovery is impossible at degree below three for the regular representation of a finite group $G$, except in the special case that $G\cong (\ZZ/2)^\ell$ for some $\ell \in \mathbb{N}$. Thus the results of these past two subsections are sharp for most groups, with respect not only to generic unique recovery but even to generic list recovery.\footnote{This is not a special property of the regular representation; generic list recovery is usually impossible in degree below 3. For finite abelian groups, this follows from \cite[Proposition~3.5]{blum2023degree} (in view of Theorem~\ref{thm:generic-list}): generic list recovery is impossible in degree below 3 unless the characters occurring in the representation are all involutions in the character group, from which it follows (assuming a faithful action, as we have) that $G\cong (\ZZ/2)^\ell$ for some $\ell \in \mathbb{N}$. Beyond $(\ZZ/2)^\ell$, notable exceptions to the general principle are $\mathrm{O}(p)$ and $\mathrm{SO}(p)$ themselves, acting in their canonical representations on a vector or a tuple of vectors; in all these cases, worst-case list recovery (and for $\mathrm{O}(p)$ even worst-case unique recovery) is achieved in degree 2. An anonymous referee suggested that one may work out a general statement for finite-dimensional representations of compact groups by analyzing the {\em group of ambiguities} introduced in \cite{bendory2022sample}, with appropriate adjustments.}

Let $V$ be the regular representation of $G$ over $\mathbb{R}$. With respect to the defining basis of $V$ consisting of the elements of $G$, an element $g\in G$ acts as a permutation matrix $\rho(g)$ consisting of a product of $|G|/\operatorname{ord}g$ disjoint cycles of length $\operatorname{ord} g$, where $\operatorname{ord} g$ is the order of $g\in G$. It follows that 
\[
\det(I-t \,\rho(g)) = (1-t^{\operatorname{ord} g})^{|G|/\operatorname{ord}g}.
\]
Thus the power series expansion of $\det(I-t \,\rho(g))^{-1}$ contains a $t^2$ term if and only if $\operatorname{ord} g = 1$ or $2$. In particular, if we let $\tau$ represent the number of involutions in $G$, then applying Molien's formula (Proposition~\ref{prop:molien}) yields that the $\mathbb{R}$-dimension of the space of degree-2 invariants in $\RxG$ is given by the coefficient of $t^2$ in 
\[
\frac{1}{|G|}\left(\frac{1}{(1-t)^{|G|}} + \frac{\tau}{(1-t^2)^{|G|/2}}\right),
\]
as the elements of order $>2$ do not contribute. Extracting this coefficient, we obtain the formula
\begin{align*}
\dim_\RR U^T_2 &= \frac{1}{|G|}\left(\binom{|G|+1}{2} + \tau \binom{|G|/2}{1}\right) \\ 
&=\frac{|G| + \tau + 1}{2}.
\end{align*}
Now if $\tau < |G|-1$, this yields $\dim_\RR U_2^T < |G| = \dim V$. In this situation, there are fewer than $|G|$ linearly independent elements in $U_2^T$. Furthermore there is only one linearly independent element in $U_1^T$ (this holds for every transitive permutation action), and its square is in $U_2^T$. Because algebraic independence is impossible without linear independence, we can conclude that there are fewer than $|G|=\dim V = \operatorname{tr.deg.}\RxG$ algebraically independent elements in $U^T_1$ and $U^T_2$ combined. It follows from Theorem~\ref{thm:generic-list} that list recovery is impossible at degree 2, as claimed.

On the other hand, if $\tau = |G|-1$ (note this is the largest possible value), then every nonidentity element in $G$ has order 2 (by definition of $\tau$). This implies $G$ is abelian, whereupon the classification theorem for finite abelian groups shows it is isomorphic to $(\ZZ/2)^\ell$. This completes the argument.

In the exceptional case that $G\cong (\ZZ/2)^\ell$, the regular representation achieves worst-case (not just generic) list recovery in degree 2. This can be seen by working in the basis of eigenvectors $y_\chi \in \Rx$ for the action, as in the proof of Theorem~\ref{thm:reg-rep-galois}. Because $G$, and therefore the character group $\widehat G$, has exponent $2$, $y_\chi^2$ is an invariant for all $\chi\in\widehat G$. (Incidentally, for the same reason, it is unnecessary to base-change to $\CC$ in order to work with the $y_\chi$'s, as $\RR$ already contains the square roots of unity.) In other words, the invariant ring $\RxG$ contains the subring $\RR[\{y_\chi^2\}]$, which is evidently generated by invariants of degree 2. But the entire ambient ring $\Rx = \RR[\{y_\chi\}]$ is integral over the subring $\RR[\{y_\chi^2\}]$, so the invariant ring $\RxG$ is as well. Because it is finitely generated, it is thus module-finite over this subring, so worst-case list recovery is achieved at degree 2 by Proposition~\ref{prop:worst-list}.

\subsection{Multi-reference alignment (MRA)} \label{subsec:MRA}

Recall that this is the case of $G = \mathbb{Z}/p$ acting on $V = \RR^p$ by cyclic shifts. It is already known that for the basic MRA problem without projection or heterogeneity, generic unique recovery is possible at degree $3$ for any $p$~\cite{BanRigWee17,PerWeBan17}. Indeed, as this action is exactly the regular representation of $\ZZ/p$, this follows from~\cite{complete-bispectrum} and from the work in the previous section.

We can also verify that for MRA with $p \ge 3$, generic list recovery is impossible at degree $2$. This follows from Theorem~\ref{thm:generic-list} because $\trdeg(\RxG) = p$ (since $G$ is finite) but the number of algebraically independent invariants of degree $\le 2$ is $\le \lfloor p/2 \rfloor + 1$ (by the formula in the previous section), which is $<p$ for $p\ge 3$. In fact, it is exactly $\lfloor p/2\rfloor + 1$.

\begin{lemma}
In MRA, the number of algebraically independent invariants of degree $\le 2$ is exactly $\lfloor p/2\rfloor + 1$.
\end{lemma}

\begin{proof}
A basis for the invariants of degree $\le 2$ is $\{\mathcal{R}(x_1),\mathcal{R}(x_1^2),\mathcal{R}(x_1 x_2),\mathcal{R}(x_1 x_3),\ldots,\mathcal{R}(x_1 x_s)\}$ with $s = \lfloor p/2 \rfloor + 1$. Here $\mathcal{R}$ denotes the Reynolds operator, which averages over cyclic shifts of the variables. For instance, $\mathcal{R}(x_1 x_2) = \frac{1}{p}(x_1 x_2 + x_2 x_3 + x_3 x_4 + \cdots + x_p x_1)$. Note that the basis above has size $\lfloor p/2 \rfloor + 2$, but there is an algebraic dependence within it because $\mathcal{R}(x_1)^2$ can be written in terms of the other basis elements. 

Meanwhile, the remaining basis elements are algebraically independent. Fix a monomial order (see Definition~\ref{def:monomial-order}) in which $x_1\ge \dots \ge x_p$. With respect to this order, the leading monomials of 
\[
\mathcal{R}(x_1),\mathcal{R}(x_1x_2),\dots,\mathcal{R}(x_1x_s)
\]
are 
\[
x_1,x_1x_2,\dots,x_1x_s.
\]
These leading monomials are evidently algebraically independent because each contains a new indeterminate. And a nontrivial algebraic relation between the above basis elements would imply a nontrivial algebraic relation between their leading monomials. So the transcendence degree of the subfield generated by invariants of degrees 1 and 2 is exactly $\lfloor p/2\rfloor + 1$.

An alternative argument is to note that there is no harm in base-changing to $\CC$, since algebraic dependence over $\RR$ is equivalent to algebraic dependence over $\CC$, whereupon, in the notation of the proof of Theorem~\ref{thm:reg-rep-galois}, a basis for the degree 1 and 2 invariants is given by $y_\mathbf{1}, \{y_\chi y_{\chi^{-1}}\}$, where $\mathbf{1}$ is the trivial character of $G\cong \ZZ/p$, and $(\chi,\chi^{-1})$ ranges over pairs of inverse characters. Again, there are $\lfloor p/2\rfloor + 2$ elements in this basis, but again, $y_\mathbf{1}^2$ is algebraically dependent on the others (in this case, only $y_\mathbf{1}$ is needed for the algebraic dependence). The remaining invariants in the basis are algebraically independent because they do not involve any of the same indeterminates.
\end{proof}

Generic list recovery is possible at degree $1$ for $p=1$ and at degree $2$ for $p = 2$. (This is true even for worst-case unique recovery; recall from Section~\ref{sec:worst-unique} that degree $|G|$ is always sufficient for this.)

We now move on to variants of the MRA problem.

\subsubsection{MRA with projection} \label{subsubsec:MRA-proj}

We now consider MRA with a projection step. We imagine that the coordinates of the signal are arranged in a circle so that $G$ acts by rotating the signal around the circle. We then observe a projection of the circle onto a line so that each observation is the sum of the two entries lying ``above'' it on the circle. This is intended to resemble the tomographic projection in cryo-EM. We formally define the setup as follows.

\begin{problem}[MRA with projection]
Let $p \ge 3$ be odd. Let $V = \RR^p$ and $G = \ZZ/p$ acting on $V$ by cyclic shifts. Let $q = (p-1)/2$ and $W = \RR^q$. Let $\Pi: V \to W$ be defined by
\begin{equation} \label{eq:mra-proj}
\Pi(v_1,\ldots,v_p) = (v_1 + v_p, v_2 + v_{p-1}, \ldots, v_{(p-1)/2} + v_{(p+3)/2}).
\end{equation}
We call the associated generalized orbit recovery problem (Problem~\ref{prob:gen-orbit}) \emph{MRA with projection}. (We consider the homogeneous case $K = 1$.)
\end{problem}

Note that since $p$ is odd, there is one entry $v_{(p+1)/2}$ which is discarded by $\Pi$. The reason we consider the odd-$p$ case rather than the seemingly more elegant even-$p$ case is because generic list recovery is actually impossible in the even-$p$ case. This is because the signals $\theta$ and $\theta + (c,-c,c,-c,\ldots)$ cannot be distinguished from the samples, even if there is no noise.

Restricting now to odd $p$, note that we cannot hope for generic unique recovery because it is impossible to tell whether the signal is wrapped clockwise or counterclockwise around the circle. In other words, reversing the signal via $(\theta_1,\ldots,\theta_p) \mapsto (\theta_p,\ldots,\theta_1)$ does not change the distribution of samples. We can still hope for generic list recovery, hopefully with a list of size exactly 2. This degeneracy is analogous to the chirality issue in cryo-EM: it is impossible to determine the chirality of the molecule (i.e.\ if the molecule is reflected about some 2-dimensional plane through the origin, this does not change the distribution of samples).

It appears that, as in the basic MRA problem, generic list recovery is possible at degree 3. We proved this for $p$ up to 21 by checking the Jacobian criterion (see Section~\ref{sec:generic-list}) on a computer, and we conjecture that this trend continues.

\begin{conjecture}
For MRA with projection, for any odd $p \ge 3$, generic list recovery is possible at degree 3.
\end{conjecture}

\noindent Note that generic list recovery is impossible at degree 2 because the addition of the projection step to basic MRA can only make it harder for $U^T_{\le d}$ to list-resolve $\theta$.

\subsubsection{Heterogeneous MRA}\label{sec:mra-het}

We now consider heterogeneous MRA, i.e.\ the generalized orbit recovery problem (Problem~\ref{prob:gen-orbit}) with $\tilde G = \ZZ/p$ acting on $\tilde V = \RR^p$ by cyclic shifts, $K \ge 2$ heterogeneous components, and no projection (i.e., $\Pi$ is the identity).

We will see that generic \textit{unique} recovery is possible at degree 3 provided that $p$ is large enough compared to $K$. First note that the number of degrees of freedom to be recovered is $\trdeg(\RxG) = Kp + K-1$ (see Propositions~\ref{prop:trdeg-finite} and \ref{prop:trdeg-het}). Let us now count the number of distinct entries of $T_d({\bf x})$ for $d \le 3$. Note that $T_d({\bf x})$ is symmetric (under permutations of indices) but we also have additional symmetries given by cyclic shifts, e.g.\ $(T_3({\bf x}))_{i,j,k} = (T_3({\bf x}))_{i+c,j+c,k+c}$ where $c$ is an integer and the sums $i+c,j+c,k+c$ are computed modulo $p$. One can compute that $T_1({\bf x})$ has 1 distinct entry, $T_2({\bf x})$ has $\lfloor p/2 \rfloor + 1$ distinct entries, and $T_3({\bf x})$ has $p + \lceil (p-1)(p-2)/6 \rceil$ distinct entries. The total number of distinct entries is
$$\mathcal{U} \defeq p + 2 + \lfloor p/2 \rfloor + \lceil (p-1)(p-2)/6 \rceil.$$

\noindent By Theorem~\ref{thm:generic-list}, \textit{list recovery} is impossible when $\mathcal{U} < Kp + K-1$. By testing the Jacobian criterion, we observe that the converse also appears to hold. Moreover, using the Hessian test (Theorem~\ref{thm:hessian}), we found generic unique de-mixing (defined in Definition~\ref{def:demix}) holds when $\mathcal{U} > Kp + K-1$.  Since in homogeneous MRA there is field generation, i.e.\ $F_{G} =  \RR(U^{T}_{\leq d})$ (Section~\ref{sec:galois}), we deduce generic unique recovery for heterogeneous MRA (for the range of parameters tested). By testing the Jacobian criterion and Hessian test in exact arithmetic on a computer, we have rigorously verified the following conjecture up to $K = 15$ and up to the corresponding critical $p$ value.

\begin{conjecture}
For heterogeneous ($K \ge 2$) MRA, generic \textup{unique} recovery is possible at degree 3 if $\mathcal{U} > Kp + K - 1$. Generic \textup{list} recovery is possible at degree 3 precisely if $\mathcal{U} \ge Kp + K -1$.  The latter condition on $\mathcal{U}$ can be stated more explicitly as follows:
\begin{itemize}
    \item $K = 2$ requires $p \ge 1$.
    \item $K = 3$ requires $p \ge 12$.
    \item $K = 4$ requires $p \ge 18$.
    \item Each $K \ge 5$ requires $p \ge 6K-5$.
\end{itemize}
\end{conjecture}
We expect that generic unique recovery is also possible in the case of equality ($\mathcal{U} = Kp + K -1$), but the Hessian test does not apply to this case.

Recent work \cite{mra-het} also studies the heterogeneous MRA problem. Similarly to the present work, they apply the method of moments and solve a polynomial system of equations in order to recover the signals. To solve the system they use an efficient heuristic method that has no provable guarantees but appears to work well in practice. Their experiments suggest that if the heterogeneous signals have i.i.d.\ Gaussian entries and uniform mixing weights, this method succeeds only when (roughly) $K \le \sqrt{p}$ instead of the condition (roughly) $K \le p/6$ that we see above (and that \cite{mra-het} also identified based on parameter-counting). Exploring this discrepancy is an interesting direction for future work.

One question of particular interest is whether this example evinces a statistical-computational gap, whereby all polynomial-time methods fail to succeed once $K$ exceeds $\sqrt{p}$. Some evidence for why we might expect this is the following analogy to tensor decomposition. Recent work on tensor decomposition \cite{tensor-decomp-sos} gives a polynomial-time algorithm to decompose a third order tensor of the form $\sum_{i=1}^N a_i^{\otimes 3} + E$ where $a_i \in \RR^p$ are i.i.d.\ from the unit sphere and $E$ is small noise, provided $N \le p^{1.5}$ (up to factors of $\log p$), and this is conjectured to be the optimal threshold for efficient decomposition. Similarly to Section~\ref{sec:jennrich}, the heterogeneous MRA problem can be cast as such a tensor decomposition problem with $N = Kp$ components $a_i$; the components are the $p$ cyclic shifts of each of the $K$ signals. Although these $a_i$ are not independent, we expect that if the signals are random then the $a_i$ are ``random enough'' for the same tensor decomposition result to hold, which exactly yields the condition $K \le \sqrt{p}$ (up to factors of $\log p$). After the initial appearance of this paper, the above argument was made into a rigorous upper bound: if $K \le \sqrt{p}/\mathrm{polylog}(p)$ and the $K$ signals are i.i.d.\ Gaussian with uniform mixing weights, then polynomial-time recovery is possible (see Chapter~5 of \cite{alex-thesis}).

\subsection{$S^2$ registration}\label{sec:ex-s2}

Recall that this is the case where the signal $\theta$ is a real-valued function defined on the unit sphere $S^2$ in $\RR^3$. The formal setup is as follows.

Let $G = \mathrm{SO}(3)$. For each $\ell = 0,1,2,\ldots$ there is an irreducible representation $V_\ell$ of $\mathrm{SO}(3)$ of dimension $2\ell+1$. These representations are of real type, i.e.\ they can be defined over the real numbers so that $V_\ell = \RR^{2\ell+1}$. Let $\mathcal{F}$ be a finite subset of $\{0,1,2,\ldots\}$ and consider the orbit recovery problem in which $G$ acts on $V = \oplus_{\ell \in \mathcal{F}} V_\ell$.

As intuition for the above setup, $V_\ell$ is the space of degree-$\ell$ \emph{spherical harmonic} functions $S^2 \to \RR$ defined on the surface of the unit sphere $S^2 \subseteq \RR^3$. The spherical harmonics are a complete set of orthogonal functions on the sphere and can be used (like a ``Fourier series'') to represent a function $S^2 \to \RR$. Thus the signal $\theta \in V$ can be thought of as a function on the sphere, with $\mathrm{SO}(3)$ acting on it by rotating the sphere. See Appendix~\ref{app:so3} for details on spherical harmonics.

The primary case of interest is $\mathcal{F} = \{1,\ldots,F\}$ for some $F$ (the number of ``frequencies''). We will see that generic list recovery is possible at degree $3$ for $10 \le F \le 16$ (and conjecturally for all $F \ge 10$).
We will see that it is convenient to not include $0 \in \mathcal{F}$, but we now justify why this is without loss of generality. $V_0$ is the trivial representation, i.e.\ the 1-dimensional representation on which every group element acts as the identity. In the interpretation of spherical harmonics, the $V_0$-component is the mean value of the function over the sphere. We claim that the $S^2$ registration problem with $0 \in \mathcal{F}$ can be easily reduced to the problem with $\mathcal{F}' = \mathcal{F} \setminus \{0\}$. This is because the $V_0$-component is itself a degree-1 invariant; given the value of this invariant, one can subtract it off and reduce to the case without a $V_0$-component (i.e.\ the case where the function on the sphere is zero-mean). Thus we have that e.g.\ generic list recovery is possible (at a given degree) for $\mathcal{F}$ if and only if it is possible for $\mathcal{F}'$.

Using Proposition~\ref{prop:trdeg-general} we compute that $\trdeg(\RxG) = p - p'$, where
$$p = \dim(V) = \sum_{\ell \in \mathcal{F}} (2\ell+1)$$
and
$$p' = \left\{\begin{array}{ll} 0 & \ell_\mathrm{max} = 0 \\ 2 & \ell_\mathrm{max} = 1 \\ 3 & \ell_\mathrm{max} \ge 2 \end{array}\right. \qquad\text{where } \ell_\mathrm{max} = \max_{\ell \in \mathcal{F}} \ell.$$

\noindent After all, $V_0$ is the trivial representation on the 1-dimensional vector space, with 3-dimensional stabilizer $\mathrm{SO}(3)$, and $V_1$ is the standard 3-dimension representation of $\mathrm{SO}(3)$ on $\RR^3$ by rotations, which yields a one-dimensional $\mathrm{SO}(2)$ stabilizer at each nonzero point. When $\ell_\mathrm{max} \ge 2$, the representation $V$ is known to have zero-dimensional stabilizer at some points (see e.g.\ \cite{so3-stabilizer}).

In the following we restrict to the case $0 \notin \mathcal{F}$ for simplicity (but recall that this is without loss of generality). There are therefore no degree-1 invariants, i.e.\ $\RxG_1$ is empty. By Theorem~\ref{thm:generic-list}, if $\dim(\RxG_2) + \dim(\RxG_3) < \trdeg(\RxG)$ then generic list recovery is impossible at degree 3; this rules out generic list recovery for $\mathcal{F} = \{1,2,\ldots,F\}$ when $F \leq 9$. (We will see below how to compute $\dim(\RxG_d)$.) Beyond this threshold, the situation is more hopeful:

\begin{theorem} If $\mathcal{F} = \{1,2,\ldots,F\}$ and $10 \leq F \leq 16$ then the degree-$3$ method of moments achieves generic list recovery. \end{theorem}
\noindent This theorem is based on computer verification of the Jacobian criterion for $10 \leq F \leq 16$ using exact arithmetic in a finite extension of $\QQ$. This result lends credence to the following conjecture.

\begin{conjecture}
Consider the $S^2$ registration problem with $0 \notin \mathcal{F}$. We conjecture the following.
\begin{itemize}
    \item Generic list recovery is possible at degree 3 if and only if $\dim(\RxG_2) + \dim(\RxG_3) \ge \trdeg(\RxG)$ (where $\trdeg(\RxG)$ is computed above and $\dim(\RxG_d)$ can be computed from Proposition~\ref{prop:S2-dim} below).
    \item In particular, if $\mathcal{F} = \{1,2,\ldots,F\}$ then generic list recovery is possible at degree 3 if and only if $F \ge 10$.
\end{itemize}
\end{conjecture}

The reason it is convenient to exclude the trivial representation is because it simplifies the parameter-counting: if we use the trivial representation then we have a degree-1 invariant $f$ and so there is an algebraic relation between the degree-2 invariant $f^2$ and the degree-3 invariant $f^3$.

We now discuss how to compute $\dim(\RxG_d)$. Using the methods in Section~4.6 of \cite{dk-book}, we can give a formula for the Hilbert series of $\RxG$; see Section~\ref{sec:S2-molien}. However, if one wants to extract a specific coefficient $\dim(\RxG_d)$ of the Hilbert series, we give an alternative (and somewhat simpler) formula:
\begin{proposition}\label{prop:S2-dim}
Consider $S^2$ registration with frequencies $\mathcal{F}$. Let $\chi_d(\phi): \RR \to \RR$ be defined recursively by
$$\chi_0(\phi) = 1,$$
$$\chi_1(\phi) = \sum_{\ell \in \mathcal{F}} \left[1 + 2 \sum_{m=1}^\ell \cos(m \,\phi)\right],\text{ and}$$
$$\chi_d(\phi) = \frac{1}{d} \sum_{i=1}^d \chi_1(i \phi) \chi_{d-i}(\phi).$$
Then we have
$$\dim(\Rx^G_d) = \frac{1}{\pi} \int_0^\pi (1-\cos \phi) \chi_d(\phi)\,\dee \phi.$$
\end{proposition}

\noindent We give the proof in Section~\ref{sec:S2-dim}. Additionally, in Appendix~\ref{app:so3-count} we give explicit formulas for the invariants (up to degree 3), which yields a combinatorial analogue of Proposition~\ref{prop:S2-dim} (up to degree 3).

\subsubsection{Formula for Hilbert series of $\RxG$}
\label{sec:S2-molien}

We can derive the Hilbert series of $\RxG$ for $S^2$ registration using the methods in Section~4.6 of \cite{dk-book}.

\begin{proposition}
Consider $S^2$ registration with frequencies $\mathcal{F}$. For $|t| < 1$, the Hilbert series of $\RxG$ is given by
$$H(t) = \sum_{z \in \mathcal{P}} \mathrm{Res}(f,z)$$
where
$$f(z) = \frac{1-\frac{1}{2}(z+1/z)}{z \prod_{\ell \in \mathcal{F}} \prod_{m=-\ell}^\ell (1-tz^m)}
= \frac{-z^{N-2}(1-z)^2}{2 \prod_{\ell \in \mathcal{F}} \left[\prod_{m=1}^\ell (z^m-t) \prod_{m=0}^\ell (1-tz^m) \right]}$$
with $N = \frac{1}{2} \sum_{\ell \in \mathcal{F}} \ell(\ell+1)$. Here $\mathrm{Res}(f,z)$ denotes the residue (from complex analysis) of the function $f$ at the point $z$, and $\mathcal{P}$ is the set of poles of $f(z)$ inside the unit circle (in $\CC$). Namely, $\mathcal{P}$ contains $t^{1/m} e^{2\pi i k/m}$ for all $m \in \{1,2,\ldots,\max_{\ell \in \mathcal{F}}\ell\}$ and for all $k \in \{0,1,\ldots,m-1\}$. If $N \le 1$, $\mathcal{P}$ also contains $0$.
\end{proposition}

\begin{proof}
Recall Molien's formula (Proposition~\ref{prop:molien}):
$$H(t) = \Ex_{g \sim \Haar(G)} \det(I-t \,\rho(g))^{-1}.$$
Note that $\det(I-t \,\rho(g))$ depends only on the conjugacy class of $g$. In $\mathrm{SO}(3)$, two elements are conjugate if and only if they rotate by the same angle $\phi$. When $g \sim \Haar(\mathrm{SO}(3))$, the angle $\phi = \phi(g)$ is distributed with density function $\frac{1}{\pi}(1-\cos \phi)$ on $[0,\pi]$ (see e.g.\ \cite{quat-rot}). If $g$ has angle $\phi$, the matrix $\rho_\ell(g)$ by which it acts on the irreducible representation $V_\ell$ has eigenvalues $e^{-i \ell \phi},e^{-i(\ell-1) \phi},\ldots,e^{i\ell \phi}$ (see e.g.\ \cite{vve-notes}). The matrix $\rho(g)$ by which $g$ acts on $V = \oplus_{\ell \in \mathcal{F}} V_\ell$ is block diagonal with blocks $\rho_\ell(g)$. Using the above we write an expression for the Hilbert series:
$$H(t) = \frac{1}{\pi} \int_0^\pi \frac{1-\cos \phi}{\prod_{\ell \in \mathcal{F}} \prod_{m = -\ell}^\ell (1-t e^{im\phi})}\, \dee \phi = \frac{1}{2\pi} \int_0^{2\pi} \frac{1-\frac{1}{2}(e^{i\phi}+e^{-i\phi})}{\prod_{\ell \in \mathcal{F}} \prod_{m = -\ell}^\ell (1-t e^{im\phi})}\, \dee \phi.$$
Now write this as a complex contour integral around the unit circle in $\CC$ and apply the residue theorem from complex analysis to arrive at the result.
\end{proof}

\subsubsection{Formula for dimension of $\RxG_d$}\label{sec:S2-dim}

The dimension of $\RxG$ can be extracted as the coefficient of $t^d$ in the Hilbert series from the previous section, but here we give a different formula based on character theory from representation theory. The \emph{character} of a representation $\rho: G \to \mathrm{GL}(V)$ (where $V$ is a finite-dimensional real vector space) is the function $\chi_V: G \to \RR$ defined by $\chi_V(g) = \mathrm{tr}(\rho(g))$.

In our case, using the eigenvalues of $\rho_\ell(g)$ from the previous section, we have
$$\chi_{V_\ell}(g) = 1 + 2 \sum_{m=1}^\ell \cos(m \,\phi(g))$$
where $\phi(g)$ is the angle of rotation of $g$. For $V = \oplus_{\ell \in \mathcal{F}} V_\ell$ we then have $\chi_V(g) = \sum_{\ell \in \mathcal{F}} \chi_{V_\ell}(g)$.

As a representation of $G = \mathrm{SO}(3)$, $\Rx_d$ is (isomorphic to) the $d$th symmetric power of $V$, denoted $S^d(V)$. (This is using the fact that a real representation is isomorphic to its dual.) There is a recursive formula for the character of $S^d(V)$:
$$\chi_{S^d(V)}(g) = \frac{1}{d} \sum_{i=1}^d \chi_V(g^i) \chi_{S^{d-i}(V)}(g).$$
\noindent This comes from the Newton--Girard formula for expressing complete homogeneous symmetric polynomials in terms of power sum polynomials.

The representation $\Rx_d = S^d(V)$ decomposes as the direct sum of irreducible representations $V_\ell$. The subspace of $\Rx_d$ consisting of all copies of the trivial representation $V_0$ (the 1-dimensional representation on which every group element acts as the identity) is precisely $\RxG_d$. Thus, $\dim(\RxG_d)$ is the number of copies of the trivial representation in the decomposition of $\Rx_d$. This can be computed using characters: $\dim(\RxG_d) = \langle \chi_{S^d(V)}, \chi_{V_0} \rangle = \langle \chi_{S^d(V)},1 \rangle$ where $\langle f_1,f_2 \rangle \defeq \EE_{g \sim \Haar(G)} [f_1(g) f_2(g)]$. Since characters are \emph{class functions} (i.e.\ they are constant on conjugacy classes), we can compute this inner product by integrating over the angle $\phi$ (as in the previous section). This yields the formula stated in Proposition~\ref{prop:S2-dim}.

\subsection{Cryo-EM}\label{sec:ex-cryo}

We now define a simple model for the cryo-EM reconstruction problem. We will use properties of the 3-dimensional Fourier transform, including the projection-slice theorem; see e.g.\ \cite{fourier-notes} for a reference.

The signal is a 3-dimensional molecule, which we can think of as encoded by a density function $f: \RR^3 \to \RR$. The 3-dimensional Fourier transform of $f$ is $\hat f:\RR^3 \to \CC$ given by
\begin{equation}\label{eq:fourier}
\hat f(k_x,k_y,k_z) = \int_{-\infty}^\infty \int_{-\infty}^\infty \int_{-\infty}^\infty e^{-2 \pi i (x k_x + y k_y + z k_z)} f(x,y,z) \,\dee x \,\dee y \,\dee z.
\end{equation}
It is sufficient to learn $\hat f$ because we can then recover $f$ using the inverse Fourier transform. $\mathrm{SO}(3)$ acts on the molecule by rotating it in 3-dimensional space (keeping the origin fixed). When $f$ is rotated in $(x,y,z)$ coordinates, $\hat f$ is also rotated in $(k_x,k_y,k_z)$-coordinates by the same rotation. Each observation is a 2-dimensional image obtained by first rotating $f$ by a random element of $\mathrm{SO}(3)$ and then projecting $f$ parallel to the $z$ axis. Specifically, the projection of $f$ is $f_\mathrm{proj}: \RR^2 \to \RR$ given by
\begin{equation} \label{eq:cryo-projection}
f_\mathrm{proj}(x,y) = \int_{-\infty}^\infty f(x,y,z) \, \dee z.
\end{equation}
By the \emph{projection-slice theorem}, the 2-dimensional Fourier transform of $f_\mathrm{proj}$ is equal to the slice $\hat f_\mathrm{slice}: \RR^2 \to \CC$ given by
$$\hat f_\mathrm{slice}(k_x,k_y) = \hat f(k_x,k_y,0).$$
Thus we think of $\hat f$ as our unknown signal with $\mathrm{SO}(3)$ acting by rotation, and with post-projection which reveals only the slice of $\hat f$ lying in the plane $k_z = 0$.

This does not yet conform to our definition of a (generalized) orbit recovery problem because the signal needs to lie in a finite-dimensional real vector space. Instead of thinking of $\hat f$ as a function on $\RR^3$, we fix a finite number $S$ of nested spherical shells in $\RR^3$, each of different radius and all centered at the origin. We consider only the restriction of $\hat f$ to these shells. We fix a finite number $F$ of frequencies and on each shell we expand $\hat f$ (restricted to that shell) in the basis of spherical harmonics, truncated to $1 \le \ell \le F$. (As in $S^2$ registration, we can discard the trivial representation $\ell=0$ without loss of generality, and it is convenient to do so.) Being the Fourier transform of a real-valued function, $\hat f$ satisfies
\begin{equation}\label{eq:real-fourier}
\hat f(-k_x,-k_y,-k_z) = \overline{\hat f(k_x,k_y,k_z)}
\end{equation}
(see (\ref{eq:fourier})) and so we can use a particular basis $H_{\ell m}$ of spherical harmonics for which the expansion coefficients are real; see Appendix~\ref{app:so3}. We have now parametrized our signal by a finite number of real values $\theta_{s \ell m}$ with $1 \le s \le S$, $1 \le \ell \le F$, and $-\ell \le m \le \ell$. In particular, the restriction of $\hat f$ to shell $s$ has expansion
$$\sum_{1 \le \ell \le F} \sum_{-\ell \le m \le \ell} \theta_{s \ell m} H_{\ell m}.$$

\noindent $\mathrm{SO}(3)$ acts on each shell by 3-dimensional rotation; see Appendix~\ref{app:so3} for the details of how $\mathrm{SO}(3)$ acts on spherical harmonics. The projection $\Pi$ reveals only the values on the equator $z = 0$ (or in spherical coordinates, $\theta = \pi/2$) of each shell. Using again the property (\ref{eq:real-fourier}), the output of $\Pi$ on each shell has an expansion with real coefficients in a particular finite basis $h_m$; see Appendix~\ref{app:so3-proj}.

\begin{remark}
There are various other choices one could make for the basis in which to represent the (Fourier transform of the) molecule. Each of our basis functions is the product of a spherical harmonic and a radial delta function (i.e.\ a delta function applied to the radius, resulting in a spherical shell). Another common basis is the Fourier--Bessel basis (used in e.g.\ \cite{cryo-2clean}) where each basis function is the product of a spherical harmonic and a radial Bessel function. More generally we can take the product of spherical harmonics with any set of radial basis function. It turns out that the choice of radial basis is unimportant because the resulting problem will be isomorphic to our case (spherical shells) and so the same results hold.
\end{remark}

We now present our results on the above cryo-EM model. We focus on identifying the regime of parameters for which generic list recovery is possible at degree 3.

We first prove that adding more shells can only make the problem easier.

\begin{proposition}\label{prop:more-shells}
Suppose generic list recovery is possible at degree $d$ for cryo-EM with $S \ge 2$ shells and $F$ frequencies. Then for any $S' \ge S$, generic list recovery is possible at degree $d$ with $S'$ shells and $F$ frequencies. The same result also holds if ``list recovery'' is replaced by ``unique recovery'' everywhere. The same result still holds for heterogeneous cryo-EM (with some fixed $K$).
\end{proposition}
\begin{proof}
To solve the problem with $S'$ shells, solve the sub-problem on each subset of $S$ shells and then patch these partial solutions together, i.e.\ rotate each partial solution by an element of $\textup{SO}(3)$ to make them all consistent. (It is not necessarily to use every subset of $S$ shells.) To show that the solutions can be patched uniquely, it is sufficient to show that for a single-shell homogeneous problem (with $F \ge 2$ frequencies), a generic signal $\theta$ is not stabilized by any element of $\textup{SO}(3)$ except the identity. To see this, first note that the restriction of $\theta$ to frequency $\ell=1$ (the standard 3-dimensional representation of $\textup{SO}(3)$) has an $\textup{SO}(2)$ stabilizer. Thus we can restrict to rotations about this fixed axis. Such rotations act diagonally (in a particular basis) on the $\ell=2$ frequency: rotation by angle $\phi$ has eigenvalues $\exp(im\phi)$ for $m \in \{-2,1,\ldots,2\}$ (see \cite{vve-notes}). It is clear that for generic $\theta$, no such rotation will stabilize $\theta$ unless $\phi = 0$. This completes the proof in the homogeneous case. The heterogeneous case is similar: now there are multiple options for which heterogeneous component should match to which in the patching process, but it will be impossible to patch them together if matched incorrectly.
\end{proof}

As in $S^2$ registration, by Proposition~\ref{prop:trdeg-general} we have for $F \geq 2$:
\begin{equation}\label{eq:cryo-trdeg}
\trdeg(\RxG) = \dim(V) - 3 = S \sum_{\ell = 1}^F (2\ell + 1) - 3 = S(F^2+2F)-3
\end{equation}
where again we have a zero-dimensional stabilizer.

In Appendix~\ref{app:so3} we give an explicit construction of the invariant polynomials in $U^T_{\le 3}$. By testing the Jacobian criterion in exact arithmetic on small examples, we arrive at the following theorem:

\begin{theorem}
Consider the homogeneous ($K=1$) cryo-EM problem with $S$ shells and $F$ frequencies.
\begin{itemize}
\item If $S=1$ then for any $F \geq 2$, generic list recovery is impossible at degree $3$.
\item If $S \ge 2$ and $2 \leq F \leq 6$, the degree-$3$ method of moments achieves generic list recovery.
\end{itemize}
\end{theorem}
The first assertion results from a simple counting argument: there are fewer invariants at degree $\leq 3$ than degrees of freedom. The second part is by Proposition~\ref{prop:more-shells} and by confirming that the Jacobian of the invariants has rank equal to $\trdeg(\RxG)$, through computer-assisted exact arithmetic over an appropriate finite extension of $\QQ$.

In floating-point arithmetic, we have further verified that the Jacobian appears to have appropriate rank for $S = 2$ and $2 \leq F \leq 10$, leading us to conjecture the following:
\begin{conjecture} If $S \ge 2$ and $F \ge 2$ then generic list recovery is possible at degree 3 (but not at degree 2).
\end{conjecture}

Intuitively, when there is a single shell ($S=1$) there are simply not enough invariants in $U^T_{\le 3}$. However, when $S \ge 2$, the number of invariants increases dramatically due to cross-terms that involve multiple shells.

When $F = 1$, generic list recovery is possible at degree 2. To see this, first note that at degree 2, no information is lost in the projection (see Appendix~\ref{app:cryo-inv}). Next note that without projection, the problem is equivalent to ``learning a rigid body'' with $S$ vectors in 3 dimensions, so generic list recovery is possible at degree 2 (see Section~\ref{sec:rigid-body}).

\subsubsection{Unprojected cryo-EM} \label{sec:cryoET}

At present, we have not bounded the list size for cryo-EM with $d=3$; the algorithm in Theorem \ref{thm:alg-ext-deg} using Gr\"obner bases is too slow for this problem.  One possible remedy would be a specially-designed algorithm for cryo-EM. To this direction, in this section, we consider \textit{unprojected} cryo-EM.  Exploiting particular structure, we show how to rapidly check that generic unique recovery holds at $d=3$ for unprojected cryo-EM.

In terms of generalized orbit recovery, the group $G$ is $\textup{SO}(3)$, the projection $\Pi$ is the identity, and the signal is the 3-dimensional Fourier transform $\hat{f}: \RR^{3} \rightarrow \CC$ of a function $f: \RR^{3} \rightarrow \RR$, restricted to $S$ shells and expanded in spherical harmonics of frequencies $\{1, \ldots, F\}$.  We shall assume $S \geq 3$ (when $S=1$, the problem is $S^2$ registration).  Our result is an \textit{efficient, constructive algorithm} showing unprojected cryo-EM
has generic unique recovery at $d=3$.   The method relies on particular ``triangular" structure in the unprojected invariants, which unfortunately disappears when the projection $\Pi$ is included.  So, the result does not easily extend to projected cryo-EM. After the initial appearance of this paper, a more refined analysis of a similar algorithm was given by~\cite{LM-so3}.

\begin{remark}
Unprojected cryo-EM is indeed of  scientific interest.  For structural biologists, unprojected cryo-EM represents a simplified version of \textit{cryo-electron tomography} (cryo-ET), with no ``missing wedge problem".  
Cryo-ET is a popular molecular imaging technique delivering noisy 3D tomograms from tilt series  (see \cite{Fra} for more).  To emphasize this connection, we freely interchange the names unprojected cryo-EM and cryo-ET here, with the tacit understanding that technically we mean a simplification of cryo-ET.
\end{remark}

In Appendix \ref{app:cryo-inv}, the degree-2 and degree-3 invariants for cryo-ET are determined, denoted there by $\mathcal{I}_{2}(s_{1}, s_{2}, \ell)$ and $\mathcal{I}_{3}(s_{1}, \ell_{1}, s_{2}, \ell_{2}, s_{3}, \ell_{3})$.  Here we will show that given the values of these invariants, we may solve for the unknowns $x_{s \ell m}$ with $1 \leq s \leq S, 1 \leq \ell \leq F,$ and $-\ell \leq m \leq \ell$, up to the group action by $\textup{SO}(3)$.  The strategy is \textit{frequency-marching}: solving first for the $\ell=1$ coefficients, then $\ell=2$ coefficients, and so on.

We first consider the degree-2 invariants with fixed frequency $\ell$ arranged into an $S \times S$ matrix,
\begin{equation*}
    \mathcal{I}_{2}(\ell) \defeq \left(  \mathcal{I}_{2}(s_1, s_2, \ell) \right)_{s_1=1, \ldots, S, s_2 = 1, \ldots, S} \in \mathbb{R}^{S \times S}.
\end{equation*}
Assembling the coefficients for the $\ell$-isotypic component  into a matrix
\begin{equation*}
    X_{\ell} \defeq \left( x_{s \ell m} \right)_{s=1, \ldots, S, m= -\ell, \ldots, \ell} \in \mathbb{R}^{S \times (2\ell+1)},
\end{equation*}
and letting 
\begin{equation*}
    \mathcal{Q}_{\ell} \defeq \left( (-1)^m \one_{m=-m'} \right)_{m=-\ell, \ldots, \ell, m' = -\ell, \ldots, \ell} \in \mathbb{R}^{(2\ell+1) \times (2\ell+1)},
\end{equation*}
(that is, $\mathcal{Q}_{\ell}$ is an anti-diagonal matrix with alternating -1's and 1's), from Appendix \ref{app:deg2invs} it holds
\begin{equation} \label{eq:cholesky}
    (2 \ell + 1)\, \mathcal{I}_2(\ell) \,= \, X_{\ell} \, \mathcal{Q}_{\ell} \, X_{\ell}^{\top}.
\end{equation}
We can interpret \eqref{eq:cholesky} as follows.  Note $\mathcal{Q}_{\ell}$ defines a non-degenerate, indefinite inner product on $\mathbb{R}^{2 \ell +1}$, since $\mathcal{Q}_\ell$ is a full-rank symmetric matrix. 
Then, the right-hand side of \eqref{eq:cholesky} is the Gram matrix for the set of $S$ vectors in $\mathbb{R}^{2\ell+1}$ given by the rows of $X_{\ell}$.
As in \cite{Kam80}, the Gram matrix determines the vectors, up to the action of an element of the indefinite orthogonal group $\textup{O}(\mathcal{Q}_{\ell}) = \{A \in \RR^{(2\ell+1)\times(2\ell+1)} \, | \, A\mathcal{Q}_{\ell} A^{\top} = \mathcal{Q} \}$.  So, degree-2 invariants  tell us that
    $X_{\ell} = Y_{\ell} A_{\ell}$,
for some unknown  $A_{\ell} \in O(\mathcal{Q}_{\ell})$ and known $Y_{\ell} \in \mathbb{R}^{S \times (2 \ell + 1)}$.  Further, $Y_{\ell}$ may be computed from $\mathcal{I}_2(\ell)$ using an indefinite variant of Cholesky factorization (see \cite{Kam80}).

Now let us start frequency-marching, by setting $\ell =1$.
We compute a $\mathcal{Q}_1$-Cholesky factorization of  $I_{2}(1)$.  This
retrieves $X_1$, up to an unknown $3 \times 3$ orthogonal matrix.  However, fibers of the map $x \mapsto (\mathcal{I}_{2}(x), \mathcal{I}_{3}(x))$ are $\textup{SO}(3)$-sets (by construction).  Thus, we are free to pick out a unique representative per generic orbit by assuming the missing $3 \times 3$ orthogonal matrix
to be $\pm I_{3}$.  Then, the missing sign is determined by comparing with the degree-3 invariants $\left( \mathcal{I}_{3}(s_{1}, 1, s_{2}, 1, s_{3}, 1)\right)_{s_1, s_2, s_3}$.  This $S \times S \times S$ tensor depends only on $X_1$, and cubically so, therefore $X_1$ is fixed.
At this point, we have solved for the $\ell=1$ isotypic component of $x$ and killed off $\textup{SO}(3)$ ambiguities, using the invariants $\mathcal{I}_{2}(s_{1},s_{1},1)$ and $\mathcal{I}_{3}(s_{1}, 1, s_{2}, 1, s_{3}, 1)$. The rest of frequency-marching just uses degree-3 invariants.

Let $\ell > 1$.  Assume for $\ell' < \ell$ the coefficients $x_{s\ell'm}$ are (generically) uniquely determined ($x_{s1m}$ fixed as above), and that we have all of these coefficients.  Now fix $s$.  Upon forward-substituting the known coefficients $x_{s\ell'm}$ into $\mathcal{I}_3$  (see Appendix \ref{app:deg3invs}), the invariants 
\begin{equation} \label{eq:deg3system}
\big{\{} \mathcal{I}_{3}(s_1, \ell_1, s_2, \ell_2, s, \ell) \,\, | \, s_{1},s_{2} = 1, \ldots, S \textup{ and } \ell_{1}, \ell_{2} < \ell \big{\}}
\end{equation}
provide a system of affine-linear equations in the unknowns $(x_{s \ell m})_{m=-\ell, \ldots, \ell} \in \mathbb{R}^{2\ell+1}$.  If this  system has full column rank, its solution is unique, we can find it, e.g., by Gaussian elimination and then march to $\ell +1$. 

Counting equations and unknowns, the invariants \eqref{eq:deg3system} determine a linear system with $2\ell + 1$ unknowns $(x_{s \ell m})_{m=-\ell, \ldots, \ell}$ and (after removing redundancies in $\mathcal{I}_{3}$) the following number of linear equations:

\vspace{-0.2em}

$$\begin{cases} \big{(} \frac{\ell}{2} - 1 \big{)}\big{(} \frac{\ell}{2}\big{)}S^{2} + \big{(}\frac{\ell}{2}\big{)} \binom{S+1}{2}, \,\,\, \textup{ if } \ell \textup{ is even} \\ \big{(} \frac{\ell -1}{2}\big{)}^2 S^2 + \big{(} \frac{\ell -1}{2}\big{)} \binom{S}{2}, \,\,\,\,\,\,\,\,\,\,\,\,\,\,\, \textup{ if } \ell \textup{ is odd}. \end{cases}$$

\vspace{0.5em}

\noindent When $S \geq 3$, equations outnumber unknowns, so we might expect full column rank to hold generically.  It remains to prove this.  We note that it is sufficient to exhibit one instance $x$ where full rank holds: since the system's coefficient matrix depends polynomially on $\{x_{s\ell'm}\}_{\ell',m}$ (for $\ell' < \ell$), full rank is a Zariski-open condition.  To this end, for each $\ell > 1$ we evaluated on computer the system's coefficient matrix at randomly chosen $x$, and computed the resulting rank.  If we found the rank to be full for each $s$, then $\ell$ isotypic component for $x$ is generically uniquely determined by \eqref{eq:deg3system}, and frequency-marching can proceed.  Performing these computations up to $F=15$, we found that generic unique recovery is possible for cryo-ET at $d=3$ (in the tested range), and our frequency-marching is an efficient algorithm for the recovery.   

\begin{theorem}\label{thm:cryo-et}
Let $F \geq 1$ and $S \geq 3$.  If the linear system  from \eqref{eq:deg3system} is full-rank generically, then homogeneous cryo-ET achieves generic unique recovery at $d=3$.  Furthermore, there exists an efficient algorithm for generic  signal recovery.   
\end{theorem}

\noindent We numerically verified that the linear system arising from \eqref{eq:deg3system} indeed generically has full-rank for $1 \leq F \leq 15$.

We note that in our above discussion of frequency marching, we have assumed access to the exact values of the invariant polynomials, whereas in reality finitely many noisy samples yield only noisy estimates. After the initial appearance of this paper, the work of~\cite{LM-so3} gave a more refined analysis of a related algorithm that accounts for the noise and achieves quasipolynomial runtime when operating on samples.

\subsubsection{Heterogeneous cryo-EM} \label{subsec:het-cryo}

We now consider heterogeneous (projected) cryo-EM ($K \ge 2$). By combining (\ref{eq:cryo-trdeg}) with Proposition~\ref{prop:trdeg-het} we can compute $\trdeg(\RxG)$. Based on testing the Jacobian criterion on small examples, we conjecture that the degree-$3$ method of moments achieves generic list recovery if and only if $\dim(U^T_2) + \dim(U^T_3) \ge \trdeg(\RxG)$. In other words, we expect no unexpected algebraic dependencies among $U^T_{\le 3}$. (Recall that there are no degree-1 invariants since we are not using the trivial representation $\ell = 0$).  Additionally, the Hessian test is passed by heterogeneous cryo-EM on small examples whenever $\dim(U^T_2) + \dim(U^T_3) > \trdeg(\RxG)$.  We conjecture this pattern continues, so that heterogenous cryo-EM enjoys generic unique de-mixing.

In Appendix~\ref{app:so3-count} we give a conjectured formula for the exact value of $\dim(U^T_2) + \dim(U^T_3)$ for all $S \ge 1, F \ge 2$. As a result we can determine for any given $S \ge 1$ and $F \ge 2$, the exact condition on $K$ for which we believe generic list recovery is possible. For $S$ and $F$ large, this condition is approximately $K \le S^2/4$.

\subsection{Cryo-EM with a symmetric molecule}
\label{sec:symm}

Consider the projected cryo-EM problem but suppose the molecule has known symmetries. Specifically, let $H$ be a subgroup of $G = \textup{SO}(3)$ and suppose we know that the molecule belongs to the subspace $V^H = \{\theta \in V \;:\; h \cdot \theta = \theta \quad\forall \,h \in H\}$. Symmetries arise frequently in biology: nearly 40\% of molecules in the Protein Data Bank~\cite{PDB} have rotational symmetries. We might expect symmetries to make the cryo-EM problem easier because there are now fewer degrees of freedom to learn. On the other hand, we might worry that they make the cryo-EM problem harder because symmetric molecules might not be generic, i.e.\ $V^H$ might be contained in the measure-zero set of ``bad'' signals. In this section we will explain how to adapt our techniques to the setting of symmetric molecules. In particular, we show how to determine the number of moments required for list recovery of a generic element of $V^H$.

We focus on a particular case where the symmetry group is $H = \ZZ/L$, generated by a rotation of $2\pi/L$ radians about the $z$-axis. By testing small examples we arrive at the following conjecture.

\begin{conjecture}
Consider (homogeneous) cryo-EM with $\ZZ/L$ symmetry. For any $S \ge 2, F \ge 2, L \ge 2$, generic list recovery is possible at degree 3. Furthermore, it is possible at degree 2 if and only if $F < L$.
\end{conjecture}

In the case $F < L$, the molecule has not only $\ZZ/L$ symmetry, but $\textup{SO}(2)$ symmetry (about the same axis).
In Section~\ref{sec:sym-het} we also point out how our techniques can handle symmetric molecules on top of heterogeneity.

\subsubsection{The homogeneous case} \label{sec:sym-hom}

Consider the setting described above, namely the homogeneous (projected) cryo-EM problem with known symmetries. Although the ideas here apply more generally, we will illustrate them by fixing the example $H = \mathbb{Z}/L$. In particular, fix $L \ge 2$ and let $H$ be generated by a rotation of $2 \pi / L$ radians about the $z$-axis. In this case, the subspace $V^H$ takes a particularly simple form: it is spanned by the spherical harmonic basis functions $H_{\ell m}$ for which $m \equiv 0 \; (\mathrm{mod} \; L)$. This can be seem because a rotation of angle $\alpha$ about the $z$-axis acts on the degree-$\ell$ complex spherical harmonics via the (diagonal) Wigner D-matrix $(D^\ell(\alpha))_{mm'} = \one_{m=m'}\, e^{-i m \alpha}$.

The inclusion map $\iota: V^H \to V$ induces the map $\phi: \RR[V] \to \RR[V^H]$ given by $\phi(f) = f \circ \iota$. In other words, $\phi$ takes a polynomial and plugs in zero for the variables corresponding to basis elements that are not in $V^H$. Let $\tilde{\RR}[{\bf x}]^G = \phi(\RR[V]^G)$ be our new invariant ring, and let $\tilde U^T_{\le d} = \phi(U^T_{\le d})$ be the new set of invariants we have access to.

Our goal is to show $\trdeg(\tilde U^T_{\le d}) = \trdeg(\tilde{\RR}[{\bf x}]^G)$; it then follows (as in Section~\ref{sec:generic-proof}) that for generic $\theta \in V^H$, the values $\{f(\theta)\}_{f \in \tilde U^T_{\le d}}$ determine the values $\{f(\theta)\}_{f \in \tilde{\RR}[{\bf x}]^G}$ up to finite ambiguity. By Theorem~\ref{thm:RxG-separates}, this means the $\textup{SO}(3)$-orbit of $\theta$ is determined up to finite ambiguity, and so we have generic list recovery.

For intuition, note that we expect $\trdeg(\tilde{\RR}[{\bf x}]^G) = \dim(V)-1$; this is because if we are given access to all the $\textup{SO}(3)$ invariants then (because we've fixed the axis of $\ZZ/L$-symmetry) we expect to be able to determine $\theta \in V^H$ up to a single degree of freedom, namely rotation by $\textup{SO}(2)$ around the axis of $\ZZ/L$-symmetry.

Rigorously, this reasoning gives us $\trdeg(\tilde{\RR}[{\bf x}]^G) \le \dim(V)-1$. This is because otherwise, knowledge of all invariants would (generically) determine a finite number of possible signals $\theta$; but this is impossible due to the $\textup{SO}(2)$ ambiguity. (In order to ensure that $\textup{SO}(2)$ ambiguity exists, we are assuming that the number of frequencies $F$ is at least $L$. Otherwise, the basis for $V^H$ only contains spherical harmonics with $m=0$, which are invariant under $\textup{SO}(2)$ action about the symmetry axis. If $F < L$ then we expect $\trdeg(\tilde{\RR}[{\bf x}]^G) = \dim(V)$.)

Therefore, to show generic list recovery at degree $d$, it is sufficient to verify via the Jacobian criterion (see Section~\ref{sec:trdeg-U}) that $\trdeg(\tilde U^T_{\le d}) = \dim(V)-1$ (or $\dim(V)$ in the case $F < L$). As reported in Section~\ref{sec:symm}, we have carried this out on small examples and observed that cryo-EM with $\ZZ/L$ symmetry appears to require degree 3 when $F \ge L$, and degree 2 when $F < L$.  

We note that it is a simple matter to rigorously prove degree 2 achieves list recovery, for $F < L$.  Following \cite{Kam80}, the projected degree-2 invariants $\mathcal{P}_{2}(s_{1}, s_{2}, m)$ span the same subspace of polynomials as the unprojected degree-2 invariants $\mathcal{I}_{2}(s_{1},s_{2},\ell)$.   As mentioned above, $F < L$ forces $x_{s\ell m}=0$ unless $m=0$.  This implies the Gram matrices $\big{(}\mathcal{I}_{2}(s_{1},s_{2},\ell)\big{)}_{s_{1},s_{2}}$ (for $\ell$ fixed) are all rank-1, so that their Cholesky factorizations, recovering $(x_{s\ell 0})_{s}$, are unique up to sign.  Generic list recovery at $d=2$ follows, with list size $2^{F}$ (accounting for an ambiguous sign per $\ell \in \{1, \ldots, F\}$).

\subsubsection{The heterogeneous case} \label{sec:sym-het}

We finish this section by explaining how to analyze  symmetry copresent with heterogeneity in cryo-EM.  To fix ideas, stick with $H = \ZZ/L$ and only consider $F \geq L$.  Suppose signals $\theta_{1}, \ldots, \theta_{K}$ (as in Problem \ref{prob:gen-orbit}) are generic in $V^{H}$.  Physically, we are assuming that each molecular conformation $\theta_{i}$ has the same known symmetry group $H$.  This assumption lets us invoke Proposition \ref{prop:trdeg-het} and Theorem \ref{thm:hessian}, and it commonly holds in practice \cite{MWC}.  By Proposition \ref{prop:trdeg-het}, the transcendence degree for our problem's invariant ring is
$K \cdot \textup{trdeg}(\tilde \RR[{\bf x}]) + K - 1$.  As $\textup{trdeg}(\tilde \RR[{\bf x}]) \leq \dim(V) - 1$, to prove generic list recovery at degree $d$, it suffices to check with the Jacobian criterion that the subspace of heterogeneous moments corresponding to $\tilde U^{T}_{\le d}$ has transcendence degree $K \cdot \big{(} \dim(V) - 1 \big{)} + K - 1$.  Then by Theorem~\ref{thm:hessian}, to verify generic unique de-mixing holds, it suffices to check $K \cdot \big{(} \dim(V) - 1 \big{)} + K - 1 < \dim(\tilde U^T_{\leq d})$ and apply the Hessian test to $K$ and $\tilde U^T_{\leq d}$.  

Testing a range of parameters $F,S,L,K$ for heterogenous cryo-EM with cyclic symmetries and $F \geq L$, we found that results vary according to $F > L$ or $F = L$.  If $F > L$, we observed that generic list recovery and generic unique de-mixing are possible at $d=3$ up to the critical $K$.  Additionally, $\dim(\tilde U^{T}_{\leq 3}) = \dim(U^{T}_{\leq 3})$, so that plugging in Conjecture \ref{conj:het-cryo}, we can predict the critical value of $K$.  By contrast, if $F = L$, the behavior is erratic and generic list recovery appears to sometimes fail at $d=3$ and critical $K$.

\section{Open questions}
\label{sec:open-questions}

We leave the following as directions for future work.
\begin{enumerate}
    \item Our methods require testing the rank of the Jacobian on a computer for each problem size. It would be desirable to have analytic results for e.g.\ (variants of) MRA in any dimension $p$.
    \item We have given an efficient test for whether generic \emph{list} recovery is possible, but have not given a similarly efficient test for generic \emph{unique} recovery. For heterogeneous problems, we presented an efficient way to reduce the question of generic unique recovery to the homogeneous problem, but then no efficient approach for the homogeneous problem.  In cases where unique recovery is impossible, it would be nice to give a tight bound on the size of the list; for instance, for MRA with projection, we conjecture that the list has size exactly 2 (due to ``chirality''), but we lack a proof for this fact. Our general algorithms for testing generic unique recovery are based on Gr\"obner bases, the calculation of which is known to be computationally hard in the worst case~\cite{Huy86}. Unfortunately, the algorithms we have proposed are also extremely slow in practice, though a faster implementation may be possible.
    \item Our procedure for recovering $\theta$ from the samples involves solving a polynomial system of equations. While solving polynomial systems is NP-hard in general, the fact that the polynomials used in the orbit recovery problem have special structure leaves open the possibility of finding an efficient (polynomial time) method with rigorous guarantees for general orbit recovery problems. Possible methods include tensor decomposition \cite{PerWeBan17}, non-convex optimization \cite{BenBouMa17,mra-het}, or numerical path-tracking \cite{SW05} with \textit{a posteriori} certification \cite{Sma86}.

    \item We have addressed the statistical limits of orbit recovery problems. However, prior work has indicated the presence of statistical-to-computational gaps in related synchronization problems~\cite{pwbm-amp}, and we expect such gaps to appear in orbit recovery problems too. As discussed in Section~\ref{sec:mra-het}, the results of \cite{mra-het} suggest a possible gap of this kind for heterogeneous MRA.
\end{enumerate}

\section{Proofs for Section~\ref{sec:stat}: statistical results} \label{sec:proofs-stat}
\subsection{Upper bounds}

We first prove Theorem~\ref{thm:statistical-upper}.
To do so, we establish the following auxiliary result.
Denote by $\Mthetad^\eps$ the $\emph{$\eps$-fattening}$ of $\Mthetad$, i.e., the set of all $\phi \in \Theta$ such that $\min_{\tau \in \Mthetad} \|\phi - \tau\| \leq \eps$.

\begin{proposition}
\label{prop:statistical-upper}
For any positive integer $n$, noise level $\sigma \geq \max_{\theta \in \Theta} \|\theta\|$, and accuracy parameter $\delta > 0$, there exists an estimator $\hat{\mathcal{M}}_n = \hat{\mathcal{M}}_n(y_1, \dots, y_n) \subseteq V$ such that, for any positive constant $\eps$, if $y_1, \dots, y_n \sim \mathrm{P}_{\theta}$~i.i.d.\ and $n \geq c_{\theta, \eps, d} \log(1/\delta) \sigma^{2d}$, then with probability at least $1-\delta$,
$$
\Mthetad \subseteq \hat{\mathcal{M}}_n \subseteq \Mthetad^\eps\,.
$$
\end{proposition}

The constant $c_{\theta, \eps, d}$ in Proposition~\ref{prop:statistical-upper} can be replaced by $c_{\theta,d} \eps^{-2}$ in the unique recovery setting if $\theta$ is suitably generic, but the dependence on $\eps$ can be worse in general.
What is key is that $c_{\theta, \eps, d}$ does not depend on $\sigma$, so that Proposition~\ref{prop:statistical-upper} implies that, if $n = \omega(\sigma^{2d})$ as $\sigma \to \infty$, then there exists a sequence of estimators $\{\mathcal{M}_n\}$ such that $\mathcal{M}_n$ converges to $\Mthetad$ in Hausdorff distance with high probability.

Let us first show how this yields Theorem~\ref{thm:statistical-upper}.
\begin{proof}[Proof of Theorem~\ref{thm:statistical-upper}]
Since $G$ is a compact group acting continuously on $V$, the orbits are compact. By assumption $\Mthetad$, is a union of a finite number of orbits, so there exists an $\eps_\theta$ such that $d_G(\mathfrak{o}_i, \mathfrak{o}_j) \geq 4\eps_\theta$ for any $i \neq j$.
For any $\eps < \eps_\theta$, let $\cN$ be an $\eps/2$-net of $V/G$, and construct $\hat{\mathcal{M}}_n$ as in Proposition~\ref{prop:statistical-upper}. Proposition~\ref{prop:statistical-upper} implies the existence of a constant $c_{\theta, \eps, d}$ such that as long as $n \geq c_{\theta,\eps,d} \log(1/\delta) \sigma^{2d}$, with probability at least $1-\delta$, $\Mthetad \subseteq \hat{\mathcal{M}}_n \subseteq \Mthetad^{\eps/2}$.

Consider the set $\mathcal{C}$ consisting of $\mathfrak{o} \in \cN$ such that $d_G(\mathfrak{o}, \hat{\mathcal{M}}_n) \leq \eps/2$. 
With probability $1-\delta$, any element of $\mathcal{C}$ is within $\eps$ of $\mathfrak{o}_i$ for some $i \in [M]$, and for each $\mathfrak{o}_i$ there exists an $\mathfrak{o} \in \mathcal{C}$ that is at most $\eps$ away.
By assumption, distinct orbits in $\Mthetad$ are separated by more than $4 \eps$, so if any two elements of $\mathcal{C}$ are separated by at most $2 \eps$, then they are close to the same element of $\Mthetad$. As a result, it is possible to partition $\mathcal{C}$ into $M$ sets $\mathcal{C}_1$, \dots, $\mathcal{C}_M$ such that $d_G(\mathfrak{o}, \mathfrak{o}') \leq 2\eps$ if $\mathfrak{o},\mathfrak{o}'$ are in the same set, and $d_G(\mathfrak{o}, \mathfrak{o}') > 2\eps$ if $\mathfrak{o}$ and $\mathfrak{o}'$ are in different sets. For $i \in [M]$, let $\hat \theta_i$ be any element of $V$ such that the orbit of $\hat \theta_i$ lies in $\mathcal{C}_i$. The claim follows.
\end{proof}

We now show how to establish Proposition~\ref{prop:statistical-upper}.
This result in fact holds for more general mixture problems, not merely those arising from the orbit recovery problems defined in Problem~\ref{prob:gen-orbit}.
For convenience, we will state and prove the theorem in its general form.
\begin{problem}[mixture recovery]\label{prob:subg-mixture}
Let $V = \RR^p$, and let $\Theta \subset V$ be compact. For $\theta \in \Theta$, let $\mu_\theta$ be a measure on $\RR^p$ whose support is contained in the unit ball, and assume the map $\theta \mapsto \mu_\theta$ is continuous with respect to the weak topology. Let $\sigma \geq 1$.
For $i \in [n] = \{1, 2, \dots, n\}$, we observe
\begin{equation*}
    y_i = x_i + \sigma \xi_i\,,
\end{equation*}
where $x_i \sim \mu_\theta$ and $\xi_i~\cN(0, I_{p \times p})$ are independent. The goal is to estimate $\theta$.
\end{problem}
Write $\mathrm{P}_\theta$ for the distribution arising from the parameter $\theta$, and let $\EE_\theta$ be expectation with respect to this distribution. 
We denote by $\EE_\theta^n$ the expectation taken with respect to $n$ i.i.d.\ samples from $\mathrm{P}_\theta$. Where there is no confusion, we also write $\EE_\theta$ for expectation with respect to the distribution $\mu_\theta.$

We require that $\mu_\theta$ have bounded support; the requirement that it be supported in the unit ball is for normalization purposes only.
We assume throughout that $\sigma \geq 1$.
The following definiton gives the generalization of Definition~\ref{def:order-d-set} to the mixture recovery problem.
\begin{definition}\label{def:gen-order-d-set}
Given a positive integer $d$ and $\theta \in V$, the \emph{order-$d$ matching set $\mathcal{M}_{\theta, d}$} is the set consisting of all $\phi \in V$ such $\EE_{\theta}[x^{\otimes k}] = \EE_{\phi}[x^{\otimes k}]$ for $k = 1, \dots, d$, where $\EE_{\zeta}$ represents expectation with respect to $x \sim \mu_\zeta$.
\end{definition}

Problem~\ref{prob:subg-mixture} generalizes Problem~\ref{prob:gen-orbit} by allowing the random vector $x_i$ to arise from more general mixtures than those arising from group actions. Note that when the mixtures do arise from a generalized orbit recovery problem, i.e., when $x_i = \Pi(g_i \cdot \theta_{k_i})$, where $g_i$ and $k_i$ are distributed as in Problem~\ref{prob:gen-orbit}, then Definition~\ref{def:gen-order-d-set} reduces to Definition~\ref{def:order-d-set}.

Having made these definitions, our goal in this section is to show that Theorem~\ref{thm:statistical-upper} holds word-for-word in the setting of mixture recovery. To do so, we show that entries of the moment tensors $\EE_{\theta}[x^{\otimes k}]$ can be estimated on the basis of $O(\sigma^{2k})$ samples from $\mathrm{P}_\theta$ .

\subsubsection{Estimation of moments}

Our estimators will be based on the univariate Hermite polynomials, which are defined by
\begin{equation*}
    H_k(x) = (-1)^k e^{x^2/2} \frac{d^k}{dx^k} e^{-x^2/2}\,.
\end{equation*}
It is well known that $H_k$ is a monic polynomial of degree $k$.
Moreover, they satsify an orthogonality relation with respect to the Gaussian law; indeed, under the inner product given by $\langle f, g \rangle = \EE_{\xi \sim \cN(0,1)}[f(\xi)g(\xi)]$, the polynomials $H_0, \dots, H_k$ form an orthogonal basis for the space of polynomials of degree at most $k$.

We rely on the crucial fact (see, e.g.,~\cite{CaiLow11}) that these polynomials satisfy
\begin{equation}\label{eqn:mean-power}
\EE_{\xi \sim \cN(0, 1)}[H_k(x + \xi)] = x^k\,.
\end{equation}

We briefly review multi-index notation.
\begin{definition}
A $p$-dimensional \emph{multi-index} is a tuple $\alpha = (\alpha_1, \dots, \alpha_p)$ of nonnegative integers.
For $x \in \RR^p$, let $x^\alpha = \prod_{j=1}^p x_j^{\alpha_j}$.
\end{definition}
For any multi-index $\alpha$, we write $|\alpha| = \sum_{j=1}^p \alpha_j$. Given independent samples $y^{(1)}, \dots, y^{(n)}$ from $\mathrm{P}_\theta$, consider the estimate for $\EE_{x \sim \mu_\theta}[x^\alpha]$ given by
\begin{equation*}
\widetilde{x^\alpha} \defeq \frac 1 n \sum_{i=1}^n \prod_{j=1}^p \sigma^{\alpha_j}H_{\alpha_j}(\sigma^{-1}y^{(i)}_j)\,.
\end{equation*}
We first show that $\widetilde{x^\alpha}$ is unbiased.

\begin{lemma}
For all $\theta \in \Theta$, $\EE_{\theta}[\widetilde{x^\alpha}] = \EE_{\theta}[x^\alpha]$.
\end{lemma}
\begin{proof}
Since $\widetilde{x^\alpha}$ is a sum of i.i.d. terms, it suffices to prove the claim for a single sample.
By~\eqref{eqn:mean-power},
\begin{align*}
\EE_\theta\Big[\prod_{j=1}^p \sigma^{\alpha_j}H_{\alpha_j}(\sigma^{-1}y_j)\Big] & = \Ex_{x \sim \mu_\theta}\Big[\Ex_{\xi_1, \dots, \xi_p \sim \cN(0, 1)} \big[\prod_{j=1}^p \sigma^{\alpha_j} H_{\alpha_j}(\sigma^{-1}(x_j + \sigma \xi_j)) \big | x\big]\Big] \\
& = \Ex_{x \sim \mu_\theta}\Big[ \prod_{j=1}^p \Ex_{\xi_j \sim  \cN(0,1)} \big[\sigma^{\alpha_j}H_{\alpha_j}(\sigma^{-1}x_j + \xi_j) \big | x\big]\Big] \\
& = \Ex_{x \sim \mu_\theta}\Big[ \prod_{j=1}^p\sigma^{\alpha_j}(\sigma^{-1}x_j)^{\alpha_j}\Big] = \EE_{\theta}[x^\alpha]\,.
\end{align*}
\end{proof}

It remains to bound the variance.
\begin{proposition}\label{prop:moment-estimator}
For any multi-index $\alpha$, there exists a constant $c_\alpha$ such that for all $\theta \in \Theta$, $\Var_\theta[\widetilde{x^\alpha}] \leq c_\alpha n^{-1}\sigma^{2|\alpha|}$.
\end{proposition}
\begin{proof}
Since $\widetilde{x^\alpha}$ is a sum of i.i.d. terms, it suffices to prove the claim for $n = 1$.
Given a multi-index $\alpha$, let $c_\alpha = \sup_{x : \|x\| \leq 1} \prod_{j = 1}^p  \EE_{\xi_j \sim \cN(0, 1)}[H_{\alpha_j}(x_j + \xi_j)^2]$, and note that $c_\alpha$ is independent of $\sigma$.
We obtain
\begin{align*}
\Var_\theta[\widetilde{x^\alpha}] & \leq \Ex_{x \sim \mu_\theta}\Big[ \prod_{j=1}^p \Ex_{\xi_j \sim  \cN(0, 1)} \big[\sigma^{2\alpha_j}H_{\alpha_j}(\sigma^{-1}x_j + \xi_j)^2 \big | x\big]\Big] \\
& \leq \sup_{x: \|x\| \leq 1}  \prod_{j=1}^p \Ex_{\xi_j \sim  \cN(0, 1)} \big[\sigma^{2\alpha_j}H_{\alpha_j}(\sigma^{-1}x_j + \xi_j)^2 \big] \\
& \leq c_\alpha \sigma^{2|\alpha|}\,,
\end{align*}
where the final inequality uses $\sigma \geq 1$.
\end{proof}

Finally, we apply the ``median-of-means'' trick~\cite{NemYud83} to show that we can combine the estimators defined above to obtain estimates for the moment tensors $\EE_\theta[x^{\otimes k}]$ for $k \leq d$ which are close to their expectation with high probability.

\begin{proposition}\label{prop:high-probability}
Let $y_1, \dots, y_n$ be i.i.d. samples from $\mathrm{P}_\theta$.
For any degree $d$ and accuracy parameter $\delta$, there exist estimators $\widehat{x^\alpha} = \widehat{x^\alpha}(y_1, \dots, y_n)$ for all $\alpha$ with $|\alpha| \leq d$ such that with probability at least $1-\delta$,
\begin{equation*}
\max_{\alpha: |\alpha| \leq d} | \EE_{\theta}[x^\alpha] - \widehat{x^\alpha}| \leq c_d\sigma^d\sqrt{\frac{\log(p/\delta)}{n}}\,,
\end{equation*}
for some constant $c_d$.
\end{proposition}
\begin{proof}
Split the samples into $m$ subsamples of equal size, for some $m$ to be specified, and for each $\alpha$ construct the $m$ estimators $\widetilde{x^\alpha_1}, \dots, \widetilde{x^\alpha_m}$ on the basis of the $m$ subsamples. (We assume for convenience that $m$ divides $n$.)
Let $\widehat{x^\alpha}$ be the median of $\widetilde{x^\alpha_1}, \dots, \widetilde{x^\alpha_m}$.

Chebyshev's inequality together with Proposition~\ref{prop:moment-estimator} implies that there exists a constant $c_d$ such that, for each $j = 1, \dots, m$ and multi-index $\alpha$,
\begin{equation*}
\Pr\Big[|\widetilde{x^\alpha_j} - \EE_{\theta}[x^\alpha]| > c_d\sigma^d \sqrt{\frac{m}{n}}\Big] \leq 1/4\,,
\end{equation*}
and since the estimators $\widetilde{x^\alpha_1}, \dots, \widetilde{x^\alpha_m}$ are independent, a standard concentration argument shows that
\begin{equation*}
\Pr\Big[|\widehat{x^\alpha} - \EE_{\theta}[x^\alpha]| > c_d\sigma^d \sqrt{\frac{m}{n}}\Big] \leq e^{-m/4}\,.
\end{equation*}
There are $\binom{p + d}{d}$ multi-indices $\alpha$ satisfying $|\alpha| \leq d$, so taking a union bound over all choices of $\alpha$ yields
\begin{equation*}
\max_{\alpha: |\alpha| \leq d} | \EE_{\theta}[x^\alpha] - \widehat{x^\alpha}| \leq c_d\sigma^d \sqrt{\frac{m}{n}}
\end{equation*}
with probability at least $1 - \binom{p + d}{d} e^{-m/4}$.
Choosing $m = 4 \log(\binom{p + d}{d}/\delta)$ and taking $c_d$ sufficiently large in the statement of the theorem yields the claim.
\end{proof}

The constant $c_d$ in the statement of Proposition~\ref{prop:high-probability} can be made explicit; however, we do not pursue this direction here.

\subsubsection{Robust solutions to polynomial systems}
We now show that approximate knowledge of the moment tensors $\EE_\theta[x^{\otimes k}]$ for $k = 1, \dots, d$ suffices to approximately recover $\theta$.
\begin{lemma}\label{lem:moment-continuity}
For all $\theta \in \Theta$ and $\eps > 0$, there exists a $\eps' > 0$ such that, if $\phi \in \Theta$ satisfies $\max_{k \leq d} \|\EE_\theta[x^{\otimes k}] - \EE_\phi[x^{\otimes k}]\|_{\infty} < \eps'$, then there exists a $\tau \in \Mthetad$ such that $\|\phi - \tau\| < \eps$.
\end{lemma}
\begin{proof}
We employ a simple compactness argument. Consider the set $F = \{\phi \in \Theta: \forall\, \tau \in \Mthetad \, \, \|\phi - \tau\| \geq \eps\}$. Since $\Theta$ is compact, so is $F$. Set
\begin{equation*}
\eps' = \inf_{\phi \in F} \max_{k \leq d} \|\EE_\theta[x^{\otimes k}] - \EE_\phi[x^{\otimes k}]\|_{\infty}\,.
\end{equation*}
Clearly if $\max_{k \leq d} \|\EE_\theta[x^{\otimes k}] - \EE_\phi[x^{\otimes k}]\|_{\infty} < \eps'$ for some $\phi \in \Theta$, then there exists a $\tau \in \Mthetad$ such that $\|\phi - \tau\| < \eps$, so it remains to check that $\eps' > 0$.

Since $\theta \mapsto \mu_\theta$ is continuous with respect to the weak topology and $\mu_\theta$ is supported on a compact set for all $\theta \in \Theta$, the moment map $\theta \mapsto \EE_\theta[x^{\otimes k}]$ is also continuous for all $k \leq d$. If $\phi \in F$, then in particular $\phi \notin \Mthetad$, so there exists a $k \leq d$ for which $\EE_\theta[x^{\otimes k}] \neq \EE_\phi[x^{\otimes k}]$. Therefore $\eps' > 0$, as desired.
\end{proof}

Lemma~\ref{lem:moment-continuity} is simply stating that the function $\phi \mapsto \min_{\tau \in \Mthetad} \|\phi  - \tau\|$ is continuous at $\theta$ with respect to the topology induced by the moment maps.
Note that, for generic $\theta$ when $\mu_\theta$ arises from an orbit recovery problem, the moment map will be continuously differentiable with a nonsingular Jacobian, so the inverse function theorem implies $\eps'$ can be taken to be $\Omega(\eps)$. In general, however, the dependence could be worse.

\subsubsection{Proof of Proposition~\ref{prop:statistical-upper}}

Construct the estimators $\widehat{x^\alpha}$ as in  Proposition~\ref{prop:high-probability}, and let
\begin{equation*}
\hat{\mathcal{M}}_n = \left\{\phi \in \Theta : \max_{\alpha: |\alpha| \leq d} | \EE_{\phi}[x^\alpha] - \widehat{x^\alpha}| \leq c_d\sigma^d\sqrt{\frac{\log(p/\delta)}{n}}\right\}
\end{equation*}

Applying Proposiiton~\ref{prop:high-probability}, we have with probability at least $1-\delta$ that $\Mthetad \subseteq \hat{\mathcal{M}}_n$ and that, for all $\phi \in \hat{\mathcal{M}}_n$,
\begin{equation*}
\max_{k \leq d} \|\EE_{\phi}[x^{\otimes k}]- \EE_\theta[x^{\otimes k}]\|_\infty
 = \max_{\alpha: |\alpha| \leq d} |\EE_{\phi}[x^{\alpha}] - \EE_\theta[x^{\alpha}]| \leq 2c_d \sigma^d\sqrt{\frac{\log(p/\delta)}{n}}\,.
\end{equation*}

By Lemma~\ref{lem:moment-continuity}, there exists an $\eps'_{\theta, \eps}$ such that, as long as $2c_d \sigma^d\sqrt{\frac{\log(p/\delta)}{n}} < \eps'_{\theta, \eps}$, then with probability at least $1-\delta$, we have the desired inclusion $\hat{\mathcal{M}}_n \subseteq \Mthetad^\eps$.

Therefore taking $n > (2 c_d/\eps'_{\theta, \eps})^2 \log(p/\delta) \sigma^{2d} = c_{\theta, \eps, d} \log(1/\delta) \sigma^{2d}$ suffices.
\qed

\subsection{Lower bounds}
We now develop tools to prove the promised lower bound, Theorem~\ref{thm:statistical-lower}.
\subsubsection{Information geometry of gaussian mixtures}\label{sec:information-geometry}
In this section, we establish an upper bound on the Kullbeck-Leibler divergence between different gaussian mixtures, which we denote by $D(\cdot \,\|\, \cdot)$. We refer the reader to~\cite{Tsy09} for the definition of the Kullback-Leibler divergence and related notions.

The proof follows the outline used in~\cite{BanRigWee17}, based off a technique developed in~\cite{LepNemSpo99,CaiLow11}.
\begin{proposition}\label{prop:kl-bound}
Let $\theta, \phi \in \Theta$, let $\sigma \geq 1$, and let $d$ be any positive integer.

There exist universal constants $C$ and $c$ such that if $\EE_\theta[x^{\otimes k}] = \EE_\phi[x^{\otimes k}]$ for $k \leq d-1$, then
\begin{equation*}
D(\mathrm{P}_\theta \, \|\, \mathrm{P}_\phi) \leq C\frac{(c\sigma)^{-2d}}{d!}\,.
\end{equation*}
\end{proposition}
\begin{proof}
We first establish the claim when $d = 1$.
Note that the condition on the moment tensors is vacuous in this case.
By the convexity of the divergence,
\begin{equation*}
D(\mathrm{P}_\theta \, \|\, \mathrm{P}_\phi) \leq \Ex_{\substack{x \sim \mu_\theta \\ x' \sim \mu_\phi}} D(\cN(x, \sigma^2) \,\|\, \cN(x', \sigma^2)) = \Ex_{\substack{x \sim \mu_\theta \\ x' \sim \mu_\phi}} \frac{\|x - x'\|^2}{2\sigma^2} \leq 2 \sigma^{-2}\,,
\end{equation*}
where in the last step we used the fact that $x$ and $x'$ lie in the unit ball almost surely.

Now, assume $d > 1$, so in particular $\EE_\theta[x] = \EE_\phi[x]$. Denote their common mean by $v$.
For $\zeta \in \{\theta, \phi\}$ denote by $\bar{\mu}_\zeta$ the distribution of $x - v$ when $x \sim \mu_\zeta$, and let $\bar{\mathrm{P}}_\zeta$ denote distribution of $y$ when $y = x + \xi$ for $x \sim \bar{\mu}_\zeta$ and $\xi \sim \cN(0, \sigma^2I)$. Since this transformation is a deterministic bijection, the data processing inequality implies $D(\mathrm{P}_\theta \, \|\, \mathrm{P}_\phi) = D(\bar{\mathrm{P}}_\theta \, \|\, \bar{\mathrm{P}}_\phi)$.

Note that $\EE_{\bar \mu_\theta}[x] = \EE_{\bar \mu_\phi}[x] = 0$ and $\EE_{\bar \mu_\theta}[x^{\otimes k}] = \EE_{\bar \mu_\phi}[x^{\otimes k}]$ for $k \leq d-1$.
Hence without loss of generality we can reduce to the case where $\mu_\theta$ and $\mu_\phi$ are both centered and are supported in a ball of radius 2.

We bound the $\chi^2$-divergence between $\mathrm{P}_\theta$ and $\mathrm{P}_\phi$.
Let $f$ be the density of a standard $p$-dimensional Gaussian and for $\zeta \in \Theta$, let $f_\zeta$ be the density of $\mathrm{P}_\zeta$, which can be written explicitly as
\begin{equation*}
f_\zeta(y) = \EE_{x \sim \mu_\zeta} \sigma^{-p} f(\sigma^{-1}(y - x)) = \sigma^{-p} f(\sigma^{-1}y) \EE_{x \sim \mu_\zeta} e^{-\frac{1}{2\sigma^2}(x^2 - 2y^\top x)}\,.
\end{equation*}
Since $\|x\| \leq 2$ almost surely with respect $\mu_\zeta$, Jensen's inequality implies that
\begin{equation}\label{eqn:density-upper-bound}
f_\zeta(y) \geq \sigma^{-p} f(\sigma^{-1}y) e^{-\frac{1}{2\sigma^2}(4 - 2 y^\top \EE_{\zeta} x)} = \sigma^{-p} f(\sigma^{-1}y) e^{-\frac{2}{\sigma^2}}\,.
\end{equation}

Recall that the $\chi^2$ divergence is defined by
\begin{equation*}
\chi^2(\mathrm{P}_\theta, \mathrm{P}_\phi) = \int \frac{(f_\theta(y) - f_\phi(y))^2}{f_\theta(y)} \mathrm dy\,.
\end{equation*}
Applying~\eqref{eqn:density-upper-bound} to the denominator, expanding the definitions of $f_\theta$ and $f_\phi$, and applying a change of variables yields
\begin{align}
\nonumber \chi^2(\mathrm{P}_\theta, \mathrm{P}_\phi) & \leq e^{2/\sigma^2}\int (\EE_{\theta} e^{-\frac{1}{2\sigma^2}(x^2 - 2y^\top x)} - \EE_{\phi} e^{-\frac{1}{2\sigma^2}(x^2 - 2y^\top x)})^2 \sigma^{-p} f(\sigma^{-1} y) \,\mathrm dy\\
\nonumber & = e^{2/\sigma^2} \int (\EE_{\theta} e^{y^\top(\sigma^{-1}x) - \frac 12 \|\sigma^{-1}x\|^2} - \EE_{\phi} e^{y^\top(\sigma^{-1}x) - \frac 12 \|\sigma^{-1}x\|^2})^2 f(y) \,\mathrm dy \\ 
& = e^{2/\sigma^2} \EE (\EE_{\theta} e^{g^\top(\sigma^{-1}x) - \frac 12 \|\sigma^{-1}x\|^2} - \EE_{\phi} e^{g^\top(\sigma^{-1}x) - \frac 12 \|\sigma^{-1}x\|^2})^2, \quad \quad g \sim \cN(0, I)\,.
\label{eqn:chi-square-expansion}
\end{align}

Given $\zeta, \zeta' \in \Theta$, let $x \sim \mu_\zeta$ and $x' \sim \mu_{\zeta'}$ be independent.
Then applying Fubini's theorem and using the expression for the moment generating function of a standard Gaussian random variable, we obtain
\begin{equation*}
\EE_{\zeta, \zeta'} \EE_{g} e^{g^\top(\sigma^{-1}(x+x')) - \frac 12 (\|\sigma^{-1}x\|^2+\|\sigma^{-1}x'\|^2)}
 = \EE_{\zeta, \zeta'} e^{\frac{x^\top x'}{\sigma^2}}\,.
\end{equation*}
Applying this expression to~\eqref{eqn:chi-square-expansion} after expanding the square produces
\begin{equation*}
\chi^2(\mathrm{P}_\theta, \mathrm{P}_\phi) \leq e^{2/\sigma^2}\left(\Ex_{\substack{x \sim \mu_\theta \\x' \sim \mu_\theta}} e^{\frac{x^\top x'}{\sigma^2}} - 2 \Ex_{\substack{x \sim \mu_\theta \\x' \sim \mu_\phi}} e^{\frac{x^\top x'}{\sigma^2}} + \Ex_{\substack{x \sim \mu_\phi \\x' \sim \mu_\phi}} e^{\frac{x^\top x'}{\sigma^2}}\right)\,,
\end{equation*}
where in each expectation the random variables $x$ and $x'$ are independent.
Since $\mu_\theta$ and $\mu_\phi$ are compactly supported, Fubini's theorem implies we can expand each term as a power series and interchange expectation and summation to produce
\begin{align*}
\chi^2(\mathrm{P}_\theta, \mathrm{P}_\phi) & \leq e^{2/\sigma^2} \sum_{k=0}^\infty \frac{\sigma^{-2k}}{k!}\left(\Ex_{\substack{x \sim \mu_\theta \\x' \sim \mu_\theta}} (x^\top x')^k - 2 \Ex_{\substack{x \sim \mu_\theta \\x' \sim \mu_\phi}} (x^\top x')^k + \Ex_{\substack{x \sim \mu_\phi \\x' \sim \mu_\phi}} (x^\top x')^k\right) \\
& = e^{2/\sigma^2} \sum_{k=0}^\infty \frac{\sigma^{-2k}}{k!}\left( \langle \EE_{\theta} x^{\otimes k}, \EE_{\theta}x^{\otimes k}\rangle
- 2 \langle \EE_{\theta} x^{\otimes k}, \EE_{x\phi}x^{\otimes k}\rangle
+ \langle \EE_{\phi} x^{\otimes k}, \EE_{\phi}x^{\otimes k}\rangle
\right) \\
& = e^{2/\sigma^2} \sum_{k=0}^\infty \frac{\sigma^{-2k}}{k!} \|\EE_\theta[ x^{\otimes k}] - \EE_\phi[ x^{\otimes k}]\|_{\text{HS}}^2 \\
& = e^{2/\sigma^2} \sum_{k=d}^\infty \frac{\sigma^{-2k}}{k!} \|\EE_\theta[ x^{\otimes k}] - \EE_\phi[ x^{\otimes k}]\|_{\text{HS}}^2\,,
\end{align*}
where $\langle \cdot, \cdot \rangle$ denotes the Frobenius inner product between tensors and $\|\cdot \|_{\text{HS}}$ denotes the Hilbert-Schmidt norm.
Since under both $\mu_\theta$ and $\mu_\phi$, $\|x\|\leq 2$ almost surely, we have for all $k \geq 2$,
\begin{equation*}
\|\EE_\theta[ x^{\otimes k}] - \EE_\phi[ x^{\otimes k}]\|_{\text{HS}}^2 \leq 2 \|\EE_\theta[ x^{\otimes k}]\|_{\text{HS}}^2 + 2 \|\EE_\phi[ x^{\otimes k}]\|_{\text{HS}}^2 \leq 4^{k+1}\,.
\end{equation*}
Therefore
\begin{equation*}
\chi^2(\mathrm{P}_\theta, \mathrm{P}_\phi) \leq 4 e^{2/\sigma^2} \sum_{k=d}^\infty \frac{4^k\sigma^{-2k}}{k!} \leq 4e^{6/\sigma^2} \frac{4^d \sigma^{-2d}}{d!}\,,
\end{equation*}
and applying the inequality $D(\mathrm{P}_\theta \,\|\, \mathrm{P}_\phi) \leq \chi^2(\mathrm{P}_\theta, \mathrm{P}_\phi)$~\cite{Tsy09} proves the claim.
\end{proof}

\subsubsection{Proof of Theorem~\ref{thm:statistical-lower}}
By~\cite[Theorem 2.2]{Tsy09}, for any test $\psi$ which tries to distinguish between $\mathrm{P}_{\tau_1}$ and $\mathrm{P}_{\tau_2}$ on the basis of $n$ samples, we have
\begin{equation*}
\max_{j = 1, 2} \Pr_{\tau_j}(\psi \neq j)\geq \frac 12 \left(1 - \sqrt{\frac 12 D(\mathrm{P}^n_{\tau_1} \,\|\, \mathrm{P}^n_{\tau_2})} \right) = \frac 12 \left(1 - \sqrt{\frac n2 D(\mathrm{P}_{\tau_1} \,\|\, \mathrm{P_{\tau_2})}} \right)\,.
\end{equation*}
If $\tau_1$ and $\tau_2$ are both in $\mathcal{M}_{\theta, d-1}$, then by Proposition~\ref{prop:kl-bound}, the corresponding distributions $\mathrm{P}_{\tau_1}$ and $\mathrm{P}_{\tau_2}$ satisfy $D(\mathrm{P}_{\tau_1} \,\|\, \mathrm{P}_{\tau_2}) \leq \frac{C}{ (c\sigma)^{2d} d!}$. 
Therefore, to achieve an error probability of at most $1/3$, we must have $n \geq 2(c\sigma)^{2d}d!/(9C) = c_d \sigma^{2d}$, as claimed.
\qed

\section{Proofs for Section~\ref{sec:algebraic}: algebraic results}\label{sec:algebraic-proofs}

\subsection{Bounding the list size for generic signals without algebraic geometry} \label{sec:generic-proof}

In this section we give an algebraic geometry-free proof of Theorem~\ref{thm:bound-list}, i.e., list size is bounded by $D \defeq [\RR(\RxG) : \RR(U)]$, which implies the first part of Theorem~\ref{thm:generic-list}, which states that if $U$ has full transcendence degree in $\RxG$, then generic list recovery is possible.  (See Section~\ref{sec:generic-converse} for the second part, when the transcendence degree of $U$ is deficient.) We do make use of basic field theory; all needed definitions and facts are in Appendix~\ref{app:fields}, the most important for this proof being the primitive element theorem (Proposition~\ref{prop:primitive-element}).

If $U$ has full transcendence degree in $\RxG$, then the fraction field $\RR(\RxG)$ of $\RxG$ is an algebraic extension of the field $\RR(U)$. There is a finite set of generators for $\RxG$ as a ring, and these same elements generate $\RR(\RxG)$ over $\RR(U)$ as a field, so this extension is of finite degree (in view of Proposition~\ref{prop:field-finite}).  Therefore, Theorem~\ref{thm:bound-list} implies the first part of Theorem~\ref{thm:generic-list}, so it remains to prove Theorem~\ref{thm:bound-list}.

\begin{proof}[Proof of Theorem~\ref{thm:bound-list}]
\leavevmode \\

Recall that the hypothesis is that $[\RR(\RxG):\RR(U)]=D<\infty$, and the desired conclusion is that $U$ generically list-resolves $\theta\in V$ with a list of size at most $D$. In characteristic zero, every algebraic extension is separable, so by the primitive element theorem (Proposition~\ref{prop:primitive-element}), $\RR(\RxG) = \RR(U)(\alpha)$ for some $\alpha \in \RR(\RxG)$. Since $\alpha$ generates a degree-$D$ extension, $\alpha$ is the root of a degree-$D$ polynomial
\begin{equation}\label{eq:alpha-root}
X^D + b_{D-1} X^{D-1} + \cdots + b_1 X + b_0
\end{equation}
with coefficients $b_i \in \RR(U)$. Furthermore, every element of $\RR(\RxG)$ can be expressed as
$$c_0 + c_1 \alpha + \cdots + c_{D-1} \alpha^{D-1}$$
with $c_i \in \RR(U)$. In particular, let $g_1,\ldots,g_k$ be generators for $\RxG$ (as an $\RR$-algebra) and write
\begin{equation}\label{eq:expand-g}
g_i = c_0^{(i)} + c_1^{(i)} \alpha + \cdots + c_{D-1}^{(i)} \alpha^{D-1}, \; i=1,\dots,k.
\end{equation}
Let $S \subseteq V$ be the subset for which $\alpha$ and all the (finitely many) coefficients $b_i, c_j^{(i)}$ have nonzero denominators; $S$ is a non-empty Zariski-open set and thus has full measure. Now fix $\theta \in S$. Given the values $f(\theta)$ for all $f \in U$, each $b_i$ takes a well-defined value in $\RR$ and so from (\ref{eq:alpha-root}) there are at most $D$ possible values that $\alpha(\theta)$ can take. From (\ref{eq:expand-g}), each value of $\alpha(\theta)$ uniquely determines all the values $g_i(\theta)$ and thus uniquely determines all the values $f(\theta)$ for $f \in \RxG$. Since $\RxG$ resolves $\theta$ (by Setup~\ref{rmk:separation-in-terms-of-AG}; see Theorem~\ref{thm:RxG-separates} for more information), this completes the proof.
\end{proof}

\subsection{Generic list recovery converse without algebraic geometry}
\label{sec:generic-converse}

In this section we give an algebraic geometry-free proof of the second part of Theorem~\ref{thm:generic-list}, stating that if $U$ does not have full transcendence degree in $\RxG$, then it fails to list-resolve $\theta$ on a full-measure subset of $V$. The argument replaces the maps of varieties found in Section~\ref{sec:algebraic} with maps of smooth manifolds, so the details are different. Instead of the work being primarily about worrying whether we have enough $\RR$-points that an argument over $\CC$ also works over $\RR$ as in Section~\ref{sec:generic-list}, everything is over $\RR$ to begin with, so all that is needed is to make sure everything in sight is a manifold. This is done with the below lemma, an immediate consequence of the constant rank theorem from smooth manifold theory (which asserts that a smooth map of constant rank is locally linear). But the main idea is the same: the transcendence degree of a set of functions is the dimension of the space of points it is able to distinguish from each other, so if $U$'s transcendence degree is deficient, there will be a positive-dimensional (and therefore infinite) family of orbits of $G$ on $V$ that all look the same to $U$. The names chosen for sets and maps ($\pi,\varphi,X$) emphasize the parallels with the arguments found in Section~\ref{sec:generic-list}.

\begin{lemma}\label{lem:hershkovits}
Let $M,N$ be smooth manifolds of dimensions $m$ and $n$, and let $\varphi:M\rightarrow N$ be a smooth map. Suppose that for some $p_0\in M$, there is an open neighborhood $U$ of $p_0$ in $M$ such that for all $p\in U$, the derivative $D\varphi_p$ has constant rank $d$. Then there exists an open neighborhood $U'\subset U$ containing $p_0$ such that
\begin{itemize}
\item$\varphi^{-1}(\varphi(p_0))\cap U'$ is a submanifold of $M$, of dimension $m-d$, and
\item  $\varphi(U')$ is a submanifold of $N$, of dimension $d$, and $\varphi$'s restriction to $U'$ is a submersion onto this submanifold.
\end{itemize} 
\end{lemma}

\begin{proof}
    This is an immediate consequence of the constant rank theorem \cite[Theorem~4.12]{lee2013}, which asserts that there are smooth coordinate charts $U'$ centered at $p_0$ and $V$ centered at $\varphi(p_0)$ such that $\varphi(U')\subset V$, in which $\varphi$ takes the form
    \[
    (x_1,\dots,x_d,x_{d+1},\dots,x_m) \mapsto (x_1,\dots,x_d,0,\dots,0).
    \]
    Then $\varphi^{-1}(\varphi(\pi_0))\cap U'$ is the $(m-d)$-dimensional submanifold of $U'$ defined by the equations $x_1=\dots=x_d=0$. And $\varphi(U')$ is the $d$-dimensional submanifold of $V$ obtained by freely varying the first $d$ coordinates with the others set to zero, and $\varphi|_{U'}$ is a submersion onto it.
\end{proof}

\begin{proof}[Proof of the second part of Theorem~\ref{thm:generic-list}]
Let $p = \dim(V)$, $\trdeg(\RxG) = q$, and $\trdeg(U) = r$ so that $r < q \le p$. Let ${\bf f} = \{f_1,\ldots,f_m\}$ be any generating set for $\RxG$, and let ${\bf g} = \{g_1,\ldots,g_n\}$ be any generating set for $\RR[U]$ (e.g., a basis for $U$). Let $S \subseteq V$ be the set of points $\theta$ for which the Jacobian $J_{\bf x}({\bf f})|_{\bf x = \theta}$ has row rank $q$ and the Jacobian $J_{\bf x}({\bf g})|_{\bf x = \theta}$ has row rank $r$; by the Jacobian criterion (see Corollary~\ref{cor:subs}), $S$ is a non-empty Zariski-open set and thus has full measure. 

Define a smooth map 
\begin{align*}
\pi:V&\rightarrow \RR^m \\
\theta &\mapsto (f_1(\theta),\dots,f_m(\theta)).
\end{align*}
Because the $f_i$'s generate $\RxG$, which contains $U$, there exist polynomials $G_1,\dots,G_n\in \RR[X_1,\dots,X_m]$ such that
\[
g_i = G_i(f_1,\dots,f_m)
\]
for each $i$. Define a second smooth map
\begin{align*}
\varphi:\RR^m &\rightarrow \RR^n\\
\mathbf{x} &\mapsto (G_1(\mathbf{x}),\dots,G_n(\mathbf{x})).
\end{align*}
Note that the composed map $\varphi\circ \pi$ is precisely $\theta\mapsto (g_1(\theta),\dots,g_n(\theta))$. Furthermore, $J_{\bf x}({\bf f})|_{\bf x = \theta}$ is the matrix of the derivative $D\pi_\theta$ written with respect to the standard bases for the tangent spaces $T_\theta V\cong V$ and $T_{\pi(\theta)}\RR^m\cong \RR^m$. Similarly, $J_{\bf x}({\bf g})|_{\bf x = \theta}$ is the matrix of the derivative $D(\varphi\circ\pi)_\theta$.

Fix $\theta_0 \in S$. Since $S$ is open and $D\pi_\theta$ has rank $q$ for every $\theta$ in $S$, the lemma applies, so there is an open neighborhood $O$ containing $\theta_0$ such that $\pi(O)$ is a submanifold of $\RR^m$, of dimension $q$, and $\pi$ restricted to $O$ is a submersion onto it. Call it $X$.

Since $X$ is a smooth submanifold of $\RR^m$, the restriction $\varphi|_X$ of $\varphi$ to $X$ is a smooth map. Furthermore, given any $x\in X$, the rank of $D(\varphi|_X)_x$ is $r$, because by the chain rule we have
\[
D(\varphi\circ\pi)_\theta = D(\varphi|_X)_x\circ D\pi_\theta,
\]
where $\theta\in O$ is such that $\pi(\theta)=x$, and the left side has rank $r$ by the construction of $S\supset O$, while $D\pi_\theta$ is surjective onto the tangent space $T_xX$ of $X$ because $\pi|_O$ is a submersion. In particular, the derivative of $\varphi|_X:X\rightarrow\RR^n$ has constant rank $r$ in a neighborhood of $x_0=\pi(\theta_0)$. Let $y_0=\varphi(x_0)$. By a second application of the lemma, there is an open neighborhood of $x_0$ in $X$, call it $P$, such that
\[
\varphi|_X^{-1}(y_0)\cap P = \varphi|_P^{-1}(y_0)
\]
is a submanifold of $X$ of dimension $q-r$. By hypothesis, $q-r>0$. Therefore, $\varphi|_P^{-1}(y_0)$ contains infinitely many points.

Because the elements of $\RxG$ (and in particular the $f_i$) are constant along orbits of $G$ on $V$, the nonempty fibers of $\pi$ are (nonempty) unions of orbits. (In fact, they are single orbits by Theorem~\ref{thm:RxG-separates}, but we do not need this here.) The set $\varphi|_P^{-1}(\varphi(x_0))$ lies within $X=\pi(O)$, thus its preimage under $\pi$ consists of infinitely many orbits, one of which contains $\theta_0$. But meanwhile, its image under $\varphi$ is exactly $y_0$. Thus all these infinitely many orbits have the same image under $\varphi\circ\pi$. 

But $\varphi\circ\pi$ is the map $\theta \mapsto (g_1(\theta),\dots,g_n(\theta))$, so this is the statement that a generating set for the ring $\RR[U]$ evaluates the same way on the orbit of $\theta_0$ and on infinitely many other orbits. In other words, $U$ fails to list-resolve $\theta_0$. Since $\theta_0\in S$ was arbitrary, $S$ is the claimed full-measure set on which $U$ fails to list-resolve orbits.
\end{proof}

\subsection{Algorithm to find $[\RR(X):\RR(Y)]$}\label{sec:compute-degree}

In this section we prove Theorem~\ref{thm:alg-ext-deg}, giving an algorithm to find the extension degree of the field extension $\RR(X)/\RR(Y)$ (and therefore bound the list size for generic list recovery, by Theorem~\ref{thm:bound-list}) when $U$ has full transcendence degree. For background on Gr\"obner bases, see Appendix~\ref{app:grobner} (and for background on fields, see Appendix~\ref{app:fields}).

\begin{lemma}
Let $d$ be such that $\RxG_{\leq d}$ generates $\RR(X)$. For all but a measure-zero set of $\alpha \in \RxG_{\le d}$, $\RR(X) = \RR(Y)(\alpha)$.
\end{lemma}

\noindent This is closely related to the primitive element theorem (\ref{prop:primitive-element}). We include a proof for completeness.

\begin{proof}
The field extension $\RR(X)/\RR(Y)$ is finite and separable (since we are in characteristic zero), so by the fundamental theorem of Galois theory (see \ref{prop:fund-thm-galois-thy}), there are only finitely many intermediate fields. (Take the normal closure of $\RR(X)/\RR(Y)$; then the intermediate fields are in bijection with subgroups of a finite group, and a subset of these lie inside $\RR(X)$.) Let $\mathcal{L}$ be the collection of intermediate fields of $\RR(X)/\RR(Y)$ that are proper subfields of $\RR(X)$. We know $\RxG_{\le d}$ is a subspace of $\RR(X)$ that generates $\RR(X)$, and therefore is not contained within any field in $\mathcal{L}$. This means each field $L \in \mathcal{L}$, being an $\RR$-vector subspace of $\RR(X)$, intersects $\RxG_{\le d}$ in a proper subspace $V_L$ of $\RxG_{\le d}$. The finite union $\cup_{L \in \mathcal{L}} V_L$ is a measure-zero subset of $\RxG_{\le d}$, and any $\alpha$ outside of it satisfies $\RR(X) = \RR(Y)(\alpha)$.
\end{proof}

\begin{proof}[Proof of Theorem~\ref{thm:alg-ext-deg}]
Consider the setup from Theorem~\ref{thm:bound-list}: given a finite set $U = \{f_1,\ldots,f_m\} \subseteq \RxG$, we want to compute $[\RR(X) : \RR(Y)]$ where $\RR(Y)$ ($= \RR(U)$) is the field generated by $U$, and $\RR(X)$ is the field of fractions of $\RxG$. We can assume $[\RR(X) : \RR(Y)]$ is finite (since we can efficiently test whether this is the case by comparing transcendence degrees using the methods of Section~\ref{sec:generic-list}). Let $d$ be such that $\RxG_{\le d}$ generates $\RR(X)$ as a field (over $\RR$). (It is sufficient for $\RxG_{\le d}$ to generate $\RxG$ as an $\RR$-algebra; such a $d$ can be computed via the algorithm in \ref{app:generator-alg} below. If $G$ is finite then $d = |G|$ is sufficient; see Section~\ref{sec:worst-unique}.) 

Let $\alpha$ be a generic element of $\RxG_{\le d}$. By the lemma, it generates the field extension. In light of this, $[\RR(X):\RR(Y)]$ is equal to the smallest positive integer $D$ for which there exists a relation
$$Q_D(f_1,\ldots,f_m) \alpha^D + \cdots + Q_1(f_1,\ldots,f_m) \alpha + Q_0(f_1,\ldots,f_m) \equiv 0$$
for polynomials $Q_i$ with $Q_D(f_1,\ldots,f_m) \not\equiv 0$. This number can be identified with a method similar to the Gr\"obner basis test for field membership (see Appendix~\ref{app:grobner}, which uses standard ideas from~\cite[Chapter~3]{iva}): Compute a reduced Gr\"obner basis $B$ for the elimination ideal $J \subseteq \RR[t_1,\ldots,t_m,\tau]$ consisting of the relations among $f_1, \ldots, f_m, \alpha$; use a monomial order that favors $\tau$. Then $[\RR(X) : \RR(Y)]$ is equal to the smallest positive integer $D$ for which $B$ contains an element of degree $D$ in $\tau$.
\end{proof}

\subsection{Hessian test relating heterogeneity to homogeneity}\label{sec:hessian}

In this section we prove Theorem \ref{thm:hessian}.  This furnishes a sufficient and efficiently-testable condition for when moments in a heterogeneous problem may be ``pulled apart uniquely" into $K$ moments for homogeneous problems.  
The ideas come from computational algebraic geometry and literature on tensor rank decompositions.  See \cite{BCO, CC1, CC2, COV1, COV2} for this line of work.

To set up, fix a heterogeneous (possibly projected)
orbit recovery problem $(\tilde{G}^{K} \rtimes S_{K}, \tilde{V}^{\oplus K} \oplus \overline{\Delta}_K, \Pi, d)$.  For simplicity of notation, we stick to the case of using all degree $\leq d$ moments, so $U = U^{T}_{\leq d} \subseteq \RR[{\bf x}]^{G}$.  The proof goes through when $U$ equals any subspace of invariant polynomials of degree $\leq d$, by replacing  $\tilde{T}_{\leq d}$ below by a  fixed basis for $\tilde{U} \subseteq \RR[{\bf \tilde{x}}]^{\tilde{G}}$. We also set $\dim(\tilde{V}) = p$ and $\dim(U^{T}_{\leq d}) = N$. 

The first goal is to reformulate heterogeneity using a classic construction from algebraic geometry.
The ``home" for heterogeneous moments is the $K^{\textup{th}}$ \textit{secant variety} of the set of homogeneous moments.  Definitions and explanation are below.  We start with the following  calculation:

$$ \mathbb{E}_{\,k \stackrel{w}{\sim} [K], \,\, g \sim \textup{Haar}(\tilde{G})} \big{[}(\Pi (g \cdot \theta_{k}))^{\otimes d}\big{]} \,\,=\,\,\, \sum_{k=1}^{K} \, w_{k} \,\, \mathbb{E}_{\, g \sim \textup{Haar}(\tilde{G})} \big{[}(\Pi (g \cdot \theta_{k}))^{\otimes d}\big{]}, $$ 

$$\,\,\, \textup{i.e.,} \,\,\,\,\,T_{d}(\Theta, w) \,\,=\,\, \sum_{k=1}^{K} \,w_{k} \,T_{d}(\theta_k).\,\,\,\,\,\,\,\,$$

\noindent Since a similar relation holds in lower degree, we may concatenate moments and write

\vspace{-0.1em}

$$ T_{\leq d}(\Theta, w) \,\,=\,\, \sum_{k=1}^{K} \,w_{k} \,T_{\leq d}(\theta_k),$$

\noindent where $T_{\leq d}(\Theta, w)$ denotes the list $\big{(}T_{1}(\theta, w) \oplus \ldots \oplus T_{d}(\theta, w)\big{)}$, and likewise for $T_{\leq d}(\theta_{k})$.  Given a basis for $\tilde{U}$, let $\tilde{T}_{\leq d}$ be the list of basis elements.  For example, for $d' = 1, \ldots, d$, we could fix (once and for all) a subset of entries in $T_{d'}$ corresponding to a basis for $U^{T}_{d'}$, and let $\tilde{T_{d'}}$ be the restriction of $T_{d'}$ to those entries, and likewise for $\tilde{T}_{\leq d}$.  With a basis chosen, 

\vspace{-0.2em}
\begin{equation}\label{eqn:mix}
\tilde{T}_{\leq d}(\Theta, w) \,\,=\,\, \sum_{k=1}^{K} \,w_{k} \,\tilde{T}_{\leq d}(\theta_k).
\end{equation}

 Define $M_{d} \subset \RR^{N}$ to be the set of all (restricted, concatenated) degree $\leq d$ moments for the homogeneous problem $(\tilde{G}, \tilde{V}, \Pi, d)$, i.e. 
 $$M_{d} \defeq \{ \tilde{T}_{\leq d}(\theta)  \,\, : \,\, \theta \in \tilde{V}\}.$$ 
Given (\ref{eqn:mix}), it is natural to recall the convex hull construction.

\begin{definition}\label{def:conv}
Given a subset $A \subset \RR^{m}$ and an integer $K \geq 1$, define $\textup{conv}_{K}(A)$ to be the set in $\RR^m$ of all convex combinations of $K$ points in $A$, i.e.:

$$\textup{conv}_{K}(A) \defeq \Big{\{} \sum_{1 \leq k \leq K} w_{k}a_{k} \,\,\, : \,\,\, w \in \Delta_{K-1} \textup{ and } a_{k} \in A \textup{ for each } k \Big{\}} \subset \RR^{m}.$$

\end{definition}

\noindent Thus, $\textup{conv}_{K}(M_{d}) \subset \RR^{N}$ is the set of all degree $\leq d$ moments for the heterogeneous problem $(G, V, \Pi, d)$, i.e.:

$$\textup{conv}_{K}(M_{d}) = \Big{\{} \, \tilde{T}_{\leq d}(\Theta, w) \,\,\, : \,\, (\Theta, \overline{w}) \in \tilde{V}^{\oplus K} \oplus \overline{\Delta}_{K-1} \Big{\}}.$$

\vspace{0.3em}

We will rely on the analog from algebraic geometry.

\begin{definition}\label{def:sec}
Given a complex affine algebraic cone\footnote{We work with affine cones, instead of projective varieties.  Also, we pass to $\CC$ for technical advantages.} $X \subset \CC^{m}$ and an integer $K \geq 1$, define the $K^{\textup{th}}$ \textit{secant variety} of $X$, denoted $\sigma_{K}(X)$, to be the Zariski closure in $\CC^m$ of all sums of $K$ points in $X$, i.e.:

$$\sigma_{K}(X) \defeq \overline{\Big{\{} \sum_{1 \leq k \leq K} p_{k} \,\,\, : \,\,\,  p_{k} \in X \,\,\, \forall \, k \Big{\}}} \subset \CC^{m}.$$

\vspace{0.05em}

\end{definition}

\vspace{-0,3em}

In the definition $X$ being a cone means that if $p \in X$ and $\lambda \in \CC$ then $\lambda p \in X$. Algebraicity of $X$ means that $X$ is the solution set of a system of polynomial equations on $\CC^{m}$. Background information on secant varieties is in   \cite{Zak}.  The following definition runs in parallel to Definition~\ref{def:demix}.

\begin{definition}
We say a secant variety $\sigma_{K}(X) \subset \CC^{m}$ has \textit{generic unique de-mixing} if for a generic point $p \in \sigma_{K}(X)$ we have a unique (up to permutation) decomposition $p = p_{1} + \ldots + p_{K}$ where $p_{k} \in X$.
\end{definition}

We now relate the convex hull and secant variety constructions.

\begin{proposition}\label{prop:reduce}
Consider the complex Zariski closure in $\CC^{N+1}$ of the cone over $M_{d}$, \textit{i.e.}:

\vspace{-0.4em}

$$\mathcal{M}_{d} \, \defeq \, \overline{\RR \cdot (M_{d} \oplus 1)} \, = \, \overline{\big{\{} \lambda \, (\tilde{T} \oplus 1) \,\,\, : \,\, \lambda \in \RR \textup{ and } \tilde{T} \in M_{d} \big{\}}} \, \subset \, \CC^{N+1}$$

\vspace{0.2em}

\noindent (where $\tilde{T} \oplus 1$ denotes concatenation and the scalar multiple $\lambda \, (\tilde{T} \oplus 1)$ lies in the cone over $M_{d}$).
Assume $\sigma_{K}(\mathcal{M}_{d})$ has generic unique de-mixing. For generic $\Theta = (\theta_{1}, \ldots, \theta_{K}) \in V = (\RR^{p})^{\oplus K}$ and $\overline{w} \in \overline{\Delta}_{K} = \RR^{K-1}$, if $\Theta' =(\theta'_{1}, \ldots, \theta'_{K}) \in V \otimes_{\RR} \CC = (\CC^{p})^{\oplus K}$ and $\overline{w}' \in \overline{\Delta}_{K} \otimes_{\RR} \CC = \CC^{K-1}$ satisfy $\tilde{T}_{\leq d}(\Theta', \overline{w}') = \tilde{T}_{\leq d}(\Theta, \overline{w})$, then $\tilde{T}_{\leq d}(\theta'_{1}) = \tilde{T}_{\leq d}(\theta_{1}), \ldots, \tilde{T}_{\leq d}(\theta_{K}) = \tilde{T}_{\leq d}(\theta'_{K})$ and $w'_{1} = w_{1}, \ldots, w'_{K} = w_{K}$ (up to permutation).
\end{proposition}

\begin{proof}
We may rewrite $\tilde{T}_{\leq d}(\Theta', \overline{w}') = \tilde{T}_{\leq d}(\Theta, \overline{w})$ and $\sum_{k=1}^{K} w'_{k} = \sum_{k=1}^{K} w_{k} = 1$  as

\vspace{-0.3em}

$$ w'_{1}\big{(}\tilde{T}_{\leq d}(\theta'_{1}) \oplus 1 \big{)} + \ldots + w'_{K}\big{(}\tilde{T}_{\leq d}(\theta'_{K}) \oplus 1 \big{)} = w_{1}\big{(}\tilde{T}_{\leq d}(\theta_{1}) \oplus 1 \big{)} + \ldots + w_{K}\big{(}\tilde{T}_{\leq d}(\theta_{K}) \oplus 1 \big{)}.$$

\vspace{0.3em}

\noindent We need to show that the right-hand side is sufficiently generic in  $\sigma_{K}(\mathcal{M}_{d})$.  More precisely, after scaling and taking complex Zariski closure, we require 

\vspace{-0.4em} 

$$\overline{\RR_{\leq 0} \cdot\big{(} \textup{conv}_{K}(M_{d}) \oplus 1\big{)}} = \sigma_{K}(\mathcal{M}_{d}).$$

\vspace{0.2em}

\noindent However this is true by the identity $\overline{f(\overline{A})} = f(\overline{A})$ when $f: X \rightarrow Y$ is a continuous map of topological spaces and $A \subset X$. Apply the identity to $X=(\CC^{N+1})^{\oplus K}$ with Zariski topology, $Y = \CC^{N+1}$ with Zariski topology, $f = K$-fold addition, and $A = \RR_{\geq 0} \cdot \big{(} M_{d} \oplus 1 \big{)}^{\times K}$.
\end{proof}

 The proposition reduces us to providing an efficient test for generic unique de-mixing of the secant variety $\sigma_{K}(\mathcal{M}_{d}) \subset \CC^{N+1}$.
For this we bring tools from the tensor decomposition literature to bear.
Specifically, we will adapt
the Hessian test for deciding tensor identifiability developed by Chiantini, Ottaviani and Vannieuwenhoven \cite{COV1, COV2} to our (more general) setting.

Background on Zariski tangent spaces to algebraic varieties \cite[Chapter~9]{iva} is required.  
Let $X \subset \CC^{m}$ be an irreducible complex affine algebraic variety, with corresponding prime ideal $I = I(X) \subset \CC[x_{1}, \ldots, x_{m}]$.  This entails that $X = \{p \in \CC^{m} \, : \, f(p) = 0 \textup{ for all } f \in I \}$ and $I = \{f \in \CC[x_{1}, \ldots, x_{m}] \,\, : \, f(p) = 0 \textup{ for all } p \in X \}$; see Hilbert's Nullstellensatz in \cite[Chapter~4]{iva}.   Now write $I = \langle f_{1}, \ldots, f_{t} \rangle$ for some ideal generators $f_{i} \in \CC[x_{1}, \ldots, x_{m}]$, which is possible by Hilbert's basis theorem \cite[Theorem~4, Chapter~2.5]{iva}.  For $p \in X$, the Zariski tangent space to $X$ at $p$, denoted $T_{p}X$, is defined as the following linear subspace of $\CC^{m}$:

\vspace{-0.8em}

$$T_{p}X \defeq \textup{ker}\, \big{(}(\partial f_{i} / \partial x_{j}) |_{x=p} \big{)}_{\substack{i = 1, \ldots, t \\ j = 1, \ldots, m}}$$

\vspace{-0.1em}

\noindent It may be shown that this subspace is independent of the choice of generators $f_{i}$.  Also, the dimension of $T_{p}X$ is constant in a nonempty Zariski-open subset of $X$. This common dimension is defined to be the dimension of $X$, written $\dim(X)$.  We have that $\dim(T_{p}X) \geq \dim(X)$ for all $p \in X$.  Additionally, if $\{\tilde{f}_{1}, \ldots, \tilde{f}_{\tilde{t}} \} \subset I$ is a subset not necessarily generating $I$, then it is clear $\dim(X) \leq \dim \textup{ker}\,(\partial \tilde{f}_{i} / \partial x_{j}) |_{x=p}$ for all $p \in X$.  Lastly, suppose $Y \subset \CC^{n}$ is another irreducible complex affine algebraic variety, and $G = (g_{1}, \ldots, g_{n}): Y \rightarrow X$ is a dominant algebraic map.  Thus $g_{i} \in \CC[y_{1}, \ldots, y_{n}]$ are functions on $Y$ and $\overline{G(Y)} = X$.  Then, by generic submersiveness we shall refer to the following two facts  (see e.g. \cite[Corollary~16.23]{Eis}).  Firstly, the derivative $DG |_{y=q} \,=\, (\partial g_{i} / \partial y_{j})|_{y=q}:  T_{q}Y \rightarrow T_{G(q)}X$ is a surjection of tangent spaces for generic $q \in Y$.  Secondly, this derivative has rank $\leq \dim(X)$ for all $q \in Y$.

Returning to the Hessian test, the following definition and proposition are the key ingredients.

\begin{definition}
Let $X \subset \CC^{m}$ be an irreducible complex affine algebraic cone.  For $p_{1}, \ldots, p_{K} \in X$, define the \textbf{\textup{contact locus}} to $X$ at $p_{1}, \ldots, p_{K}$, written $\textup{Con}_{p_{1}, \ldots, p_{K}}(X) \subset X$, as 

\vspace{-0.1em}

$$\textup{Con}_{p_{1}, \ldots, p_{K}}(X) \defeq \overline{\big{\{}p \in X \,\, : \,\,\, T_{p}X \,\subset \, T_{p_{1}}X + \ldots + T_{p_{K}}X \big{\}}},$$

\vspace{0.4em}

\noindent which is a cone in $X$ containing $p_{1}, \ldots, p_{K}$.
\end{definition}

\vspace{0.45em}

\begin{proposition}\label{prop:contact}
Let $X \subset \CC^{m}$ be an irreducible complex affine algebraic cone.  Suppose that $\dim \sigma_{K}(X) = K \dim(X)$.  If $\sigma_{K}(X)$ does not have generic unique de-mixing, then for generic points $p_{1}, \ldots, p_{K} \in X$, each irreducible component of $\textup{Con}_{p_{1}, \ldots, p_{K}}(X)$ has dimension $\geq 2$.   
\end{proposition}

\begin{proof}
The proof of \cite[Prop.~2.3]{COV1} remains valid when we replace ``Segre variety $\mathcal{S} \subset \CC^{n_{1} \times \ldots \times n_{d}}$ (i.e., rank-$1$ tensors)" by ``irreducible cone $X \subset \CC^{m}$\hspace{0.04em}".
\end{proof}

Thus the contact locus upgrades a finite-list statement to a conclusion about the dimension of a variety.

\begin{proof}[Proof of Theorem \ref{thm:hessian}]
The steps making up the Hessian test are highlighted in bold below.
 Recall we have assumed $U = U^{T}_{\leq d}$. 
 
First we argue $[F_{G} : \RR(U)] = [F_{\tilde{G}} : \RR(\tilde{U})]^{K}$  if $\sigma_{K}(\mathcal{M}_{d}) \subset \CC^{N+1}$ has generic unique de-mixing.
Consider
 
\vspace{-1em}
 
\begin{alignat*}{2}
 \big{(}\theta_{1}, \ldots, \theta_{K}, \overline{w}\big{)}  \,\,\, &\overset{f}{\longmapsto} \,\,\, \big{(}\tilde{T}_{\leq d}(\theta_{1}), \ldots, \tilde{T}_{\leq d}(\theta_{k}), \overline{w}\big{)} \,\,\, &\overset{g}{\longmapsto}& \,\,\, \sum_{k=1}^{K} w_{k} \tilde{T}_{\leq d}(\theta_{k}).
\end{alignat*}.
\vspace{-1em}
 
\noindent We regard the left and middle points as in the relevant quotient spaces, i.e.
 
\vspace{-1em}
 
\begin{alignat*}{2}
 \textup{Spec}\big{(} \RR\big{[}{\bf x}\big{]}^{G}\big{)} \,\,\, &\overset{f}{\longrightarrow} \,\,\, \textup{Spec}\big{(} \RR\big{[}\tilde{U}({\bf x}^{(1)}), \ldots, \tilde{U}({\bf x}^{(K)}), w_{1}, \ldots, w_{K}\big{]}^{S_{K}} \big{)}\,\,\, &\overset{g}{\longrightarrow}\,\,\, \textup{Spec}\big{(} \RR \big{[} U({\bf x}) \big{]} \big{)},
\end{alignat*}
 
\vspace{0.3em}

\noindent Due to Theorem \ref{thm:finite-gen}, the invariant rings are finitely generated, so the spaces  are algebraic varieties. Further, from the assumption $\textup{trdeg}(U) = \trdeg(\RR[{\bf x}]^{G})$, the source and target of $g \circ f$ share the same dimension. Since $g \circ f$ is dominant, it is a finite map of algebraic varieties with a well-defined degree, and likewise for $f$ and $g$.  The degree equals the degree of the corresponding field extension, or equivalently the cardinality of a general fiber after base change to $\CC$ (see e.g. \cite{Liu}).  
By the first description, $[F_{G} : \RR(U)] = \textup{deg}(g \circ f) = \textup{deg}(f) \cdot \textup{deg}(g)$.  By the second, $\textup{deg}(f) = [F_{\tilde{G}} : \RR(\tilde{U})]^{K}$ and $\textup{deg}(g) = 1$ using Proposition \ref{prop:reduce}, if $\sigma_{K}(\mathcal{M}_{d})$ has generic unique de-mixing.  The desired implication follows.

Next we give an efficient algorithm to verify the secant variety $\sigma_{K}(\mathcal{M}_{d})$ has generic unique de-mixing.  It is based on the contact locus defined above, and a tangent space computation.   

Clearly $\CC^{p+1} = (\tilde{V} \otimes \CC) \oplus \CC$ maps dominantly to $\mathcal{M}_{d} \subset \CC^{N+1}$ via $(\theta, \lambda) \mapsto (\lambda \tilde{T}_{\leq d}(\theta) \oplus \lambda)$.  Denoting this by $\pi$, we shall pull all contact locus computations back to $\CC^{p+1}$ via $\pi$.  This obviates computing $I(\mathcal{M}_{d})$.  

\textbf{Fix} $K$ \textbf{particular points} $(\theta_{1} \oplus \lambda_{1}), \ldots, (\theta_{K} \oplus \lambda_{K}) \in \CC^{p+1}$.  Set
$p_{1} = \pi(\theta_{1} \oplus \lambda_{1}), \ldots, p_{K} = \pi(\theta_{K} \oplus \lambda_{K})$.  To check $\dim \sigma_{K}(\mathcal{M}_{d}) = K \dim(\mathcal{M}_{d})$ (so generic list recovery holds in the heterogeneous problem), observe
$$\textup{Im}(D\pi|_{\theta_{1} \oplus \lambda_{1}}) + \ldots + \textup{Im}(D\pi|_{\theta_{K} \oplus \lambda_{K}}) \,\subset\, T_{p_{1}}(\mathcal{M}_{d}) + \ldots + T_{p_{K}}(\mathcal{M}_{d}) \,\subset\, T_{p_{1}+\ldots+p_{K}}\sigma_{K}(\mathcal{M}_{d}).$$ 
\textbf{Evaluate $K$ Jacobian matrices $D\pi$, concatenate them into} 

\vspace{-0.6em}

$$J  \defeq \begin{pmatrix} D\pi|_{\theta_{1} \oplus \lambda_{1}} & & \ldots & & D\pi|_{\theta_{K} \oplus \lambda_{K}}\end{pmatrix}^\top \in \CC^{K(p+1) \times (N+1)}$$

\vspace{0.4em}

\noindent \textbf{and compute $\textup{rank}(J)$}.  
If $\textup{rank}(J) = K \dim(\mathcal{M}_{d})$, then this verifies $\dim \sigma_{K}(\mathcal{M}_{d}) = K \dim(\mathcal{M}_{d})$ (by generic submersiveness and the chain of inclusions above).
Else, the algorithm does not proceed with $\theta_{1} \oplus \lambda_{1}, \ldots, \theta_{K} \oplus \lambda_{K}$, either choosing points again or terminating with the output ``fail".

Continuing, set $\mathcal{C} = \textup{Con}_{p_{1}, \ldots, p_{K}}(\mathcal{M}_{d})$ for the contact locus to $\mathcal{M}_{d} \subset \CC^{N+1}$ at $p_{1}, \ldots, p_{k}$ and consider
$$ \tilde{\mathcal{C}} \defeq \{\theta \oplus \lambda \in \CC^{p+1} \,\, : \,\, \textup{Im}(D\pi|_{\theta \oplus \lambda}) \subset \textup{Im}(D\pi|_{\theta_{1} \oplus \lambda_{1}}) + \ldots + \textup{Im}(D\pi|_{\theta_{K} \oplus \lambda_{K}})\} \subset \CC^{p+1}.$$
We regard $\tilde{\mathcal{C}}$ as an approximation to $\pi^{-1}\mathcal{C}$. 
We seek to bound the dimension of $\tilde{\mathcal{C}}$, and deduce the dimensions of $\pi^{-1} \mathcal{C}$ and $\mathcal{C}$ (particularly, that $\mathcal{C}$ has a 1-dimensional component).

\textbf{Compute $\textup{ker}(J) \subset \CC^{N+1}$ and fix a particular element\footnote{One may instead use a basis for $\textup{ker}(J)$, and obtain a stacked version of the test presented here. This is what is done in \cite{COV1} for tensors. In our experience, the stacked variant is often not needed and noticeably slower.} $\ell \in \textup{ker}(J)$}.  By the assumption $\textup{trdeg}(U) < \dim(\tilde{U})$, we may take nonzero $\ell$. Thus, $\ell$ is an equation for the RHS subspace in the definition of $\tilde{\mathcal{C}}$.  This imposes $p+1$ equations on $\tilde{\mathcal{C}}$, namely, each column of $D\pi|_{(\theta \oplus \lambda)}$ should have vanishing inner product with $\ell$.  Now, obviously $\theta_{1} \oplus \lambda_{1} \in \tilde{\mathcal{C}}$.  By the preceding discussion on Zariski tangent spaces, it follows that

\vspace{-0.6em}
\begin{equation}
\CC^{p+1} \,\,\,\,\, \supset \,\,\,\,\, \textup{ker}  \,\, \textup{Jac}\big{(} \ell^{T} \cdot D\pi|_{\theta \oplus \lambda}\big{)}\big{|}_{\theta_{1} \oplus \lambda_{1}} \,\,\, \supset \,\,\, T_{\theta_{1} \oplus \lambda_{1}}  (\tilde{\mathcal{C}}).
\label{eq:incl}
\end{equation}
 
\noindent In the middle term, we take Jacobian matrix of the $p+1$ equations of degree $\leq d+1$ for $\tilde{\mathcal{C}} \subset \CC^{p+1}$, and then evaluate at $\theta_{1} \oplus \lambda_{1}$.  This produces a $(p+1)\times(p+1)$ matrix.  In fact, it is easy to see we get the Hessian matrix for the pull back $\pi^{*} \ell \in \CC[\theta, \lambda]$, i.e.

\vspace{-0.4em}

$$\textup{Hess}(\pi^{*} \ell) = \begin{pmatrix} \frac{\partial^{2}(\pi^{*}\ell)}{\partial \theta^{1} \partial \theta^{1} }& \frac{\partial^{2}(\pi^{*}\ell)}{\partial \theta^{1} \partial \theta^{2}} & \ldots & \frac{\partial^{2}(\pi^{*}\ell)}{\partial \theta^{1} \partial \lambda} \\[0.3em] \vdots & \vdots & \ddots & \vdots \\[0.3em] \frac{\partial^{2}(\pi^{*}\ell)}{\partial \lambda \partial \theta^{1} }& \frac{\partial^{2}(\pi^{*}\ell)}{\partial \lambda \partial \theta^{2}} & \ldots & \frac{\partial^{2}(\pi^{*}\ell)}{\partial \lambda \partial \lambda}  \end{pmatrix}$$ 

\vspace{0.4em}

\noindent evaluated at $\theta \oplus \lambda = \theta_{1} \oplus \lambda_{1}$.  

\textbf{Compute the Hessian matrix $\textup{Hess}(\pi^{*} \ell)|_{\theta_{1} \oplus \lambda_{1}}$ and find its rank}.  If the Hessian rank equals $p-\max_{\tilde{v}}\dim(\tilde{G}\cdot \tilde{v})$, then, this proves $\sigma_{K}(\mathcal{M}_{d})$ has generic unique de-mixing (justification below), so the Hessian test terminates with the output ``pass".  Else, the algorithm re-chooses $\theta_{1} \oplus \lambda_{1}, \ldots, \theta_{K} \oplus \lambda_{K} \in \CC^{p+1}$ and starts anew, or it terminates with the output ``fail".

It remains to prove correctness of the algorithm, i.e. $\sigma_{K}(\mathcal{M}_{d})$ has generic unique de-mixing when ``pass" is outputted.  We first claim the following property
must hold for generic choices of $\theta_{1} \oplus \lambda_{1}, \ldots, \theta_{K} \oplus \lambda_{K}$ and $\ell$ since it holds for a particular choice: 

\vspace{-0.4em}

\begin{itemize}
\item \textup{the Hessian rank is maximal}: $\textup{rank}\big{(}\textup{Hess}(\pi^{*}\ell)|_{\theta_{1} \oplus \lambda_{1}}\big{)} \geq p-\max_{\tilde{v}}\dim(\tilde{G}\cdot \tilde{v}). $
\end{itemize}

\noindent Note $\textup{rank}(J) = K \dim(\mathcal{M}_{d})$ holds generically by generic submersiveness, since it holds for an instance.  Consider 

\vspace{-1.2em}

$$ \mathcal{X} \defeq \big{\{} (\theta_{1} \oplus \lambda_{1}, \ldots, \theta_{K} \oplus \lambda_{K}, \ell) \,\,\, : \,\,\, \textup{rank}(J) = K \dim(\mathcal{M}_{d}), \, \ell \in \textup{ker}(J) \big{\}} \subset (\CC^{p+1})^{\times K} \times \CC^{N+1}.$$

\vspace{0.5em}

\noindent Then $\mathcal{X}$ is an irreducible variety by \cite[Exercise~14.3]{Eis}, as the projection $\mathcal{X} \rightarrow (\CC^{p+1})^{\times K}$ has irreducible fibers of constant dimension (namely isomorphic Grassmannians).
At the same time,

\vspace{-0.6em}

$$\big{\{} (\theta_{1} \oplus \lambda_{1}, \ldots, \theta_{K} \oplus \lambda_{K}, \ell) \in \mathcal{X} \,\,\, : \,\,\, \textup{rank } \textup{Hess}(\pi^{*} \ell)|_{\theta_{1} \oplus \lambda_{1}} < p-\max_{\tilde{v}}\dim(\tilde{G}\cdot \tilde{v}) \big{\}}$$

\vspace{0.1em}

\noindent is a Zariski closed subset of $\mathcal{X}$, and a strict subset thanks to the instance above.  So irreducibility of $\mathcal{X}$ implies $\textup{rank} \big{(}\textup{Hess}(\pi^{*} \ell)|_{\theta_{1} \oplus \lambda_{1}} \big{)} \geq p-\max_{\tilde{v}}\dim(\tilde{G}\cdot \tilde{v})$ holds generically, and the claim follows.

We now finish proving correctness.  
Let $\theta_{1} \oplus \lambda_{1}, \ldots, \theta_{K} \oplus \lambda_{K}$ and $\ell$ be generic, in particular satisfying the maximal Hessian rank property.
By generic submersiveness, the column space of $J$ equals $T_{p_{1}}(\mathcal{M}_{d}) + \ldots + T_{p_{k}}(\mathcal{M}_{d})$.  By definition of $\tilde{\mathcal{C}}$, this implies $\tilde{\mathcal{C}} \supset \pi^{-1} \mathcal{C}$, and so

\vspace{-0.6em}

$$T_{\theta_{1} \oplus \lambda_{1}}(\tilde{\mathcal{C}}) \,\,\supset \,\,T_{\theta_{1} \oplus \lambda_{1}}(\pi^{-1}\mathcal{C}).$$  

\vspace{0.3em}

\noindent By upper semi-continuity of fiber dimension \cite[Theorem~14.8]{Eis}, $\dim T_{\theta_{1} \oplus \lambda_{1}}(\pi^{-1}\mathcal{C}) \geq \dim(\tilde{G}) + 1$.  On the other hand, by (\ref{eq:incl}) and the Hessian rank property, $\dim T_{\theta_{1} \oplus \lambda_{1}}(\tilde{\mathcal{C}}) \leq \dim(\tilde{G}) + 1$.  It follows that
$\dim T_{\theta_{1} \oplus \lambda_{1}}(\pi^{-1}\mathcal{C}) = \dim T_{\theta_{1} \oplus \lambda_{1}}(\tilde{\mathcal{C}}) = \dim(\tilde{G}) + 1$.  By $\dim T_{\theta_{1} \oplus \lambda_{1}}(\pi^{-1}\mathcal{C}) = \dim(\tilde{G}) + 1$ and \cite[Theorem~14.6]{Eis} again, we get $\dim T_{p_{1}} \mathcal{C} = 1$.  
Finally, since $\dim \sigma_{K}(\mathcal{M}_{d}) = K \dim(\mathcal{M}_{d})$, we apply Proposition \ref{prop:contact} to conclude $\sigma_{K} (\mathcal{M}_{d})$ has generic unique de-mixing as required.  
The proof of Theorem~\ref{thm:hessian} is complete.
\end{proof}

\subsection{Algorithm to test if $U$ is a separating set}\label{sec:separating-alg}

In this section we prove Proposition~\ref{prop:alg-separate}, which asserts an algorithm to test if a proposed subspace $U$ is a separating set for $\RxG$, given an oracle to find an interior feasible point for an SDP in exact arithmetic. The algorithm is an adaptation of a method of Kemper \cite{kemper2003}, with two changes. First, it uses the real radical membership test of \cite{reid-wang-wolkowicz} in place of an ordinary radical membership test (this is what requires the SDP oracle). Second, it is significantly simpler because our $G$ is linearly reductive (i.e., there exists a Reynolds operator; see Definition~\ref{def:reynolds}). Kemper's method studiously avoids the Reynolds operator, because it was designed as a subroutine in an algorithm to compute generators for the invariant ring of an algebraic group defined over a field of positive characteristic that is reductive, but not linearly reductive. In our context, we have a Reynolds operator, so we can reverse the order: we use Derksen's algorithm \cite{derksen1999} to compute generators as a subroutine for testing a separating set.

For readers without an algebraic background, radicals, real radicals, and the ordinary and real Nullstellensatz are discussed in Appendix~\ref{app:comm-alg}.

We describe the main idea of the algorithm before giving the proof. Recall from Section~\ref{sec:worst-unique} that a separating set for $\RxG$ is a subset $U\subset\RxG$ such that any pair of orbits $\mathfrak{o}_1,\mathfrak{o}_2\subset V(\RR)$ for the action of $G$ which are distinguished by some element $f$ of $\RxG$ (in the sense that $f(\mathfrak{o}_1)\neq f(\mathfrak{o}_2)$) are also distinguished by an element of $U$.\footnote{In fact, because $G$ is compact, any two different orbits in $V(\RR)$ are distinguished by some element of $\RxG$, as discussed in Setup~\ref{rmk:separation-in-terms-of-AG}; thus in our situation we could just as well have required that any two different orbits are distinguished by an element of $U$. However, this is not true if $G$ is not compact, so we have stated the definition of a separating set in the more general way, by only requiring $U$ have the same power to separate orbits that the entire invariant ring $\RxG$ has. In any case this calls attention to the property of separating sets that will actually be used in what follows.} We begin by defining the variety
\[
\mathcal S \defeq \{ (\mathbf u,\mathbf v)\in V\times V : f(\mathbf u)=f(\mathbf v)\text{ for all } f\in \RxG\}.
\]
Here we have introduced new indeterminates $\mathbf u \defeq u_1,\dots,u_p$ and $\mathbf v\defeq v_1,\dots,v_p$ for convenience to work with the product $V\times V$. This variety was introduced in \cite{kemper2003}, and is now called the {\em separating variety} \cite[Section~4.9.1]{dk-book} (and see also \cite{dufresne} where it appears as the {\em separating scheme}). A collection of invariant functions $g_1,\dots,g_n \in \RxG$ is then a separating set if the vanishing set of 
\[
g_1(\mathbf u) - g_1(\mathbf v), \dots, g_n(\mathbf u) - g_n(\mathbf v)
\]
in $V\times V$ is exactly $\mathcal S$. There is a subtlety to note: as $V$, and consequently $V\times V$ and $\mathcal S$, are defined over $\RR$, one can assess this criterion at either of two levels: by considering all $\CC$-points, or only the $\RR$-points. If the set of $\RR$-points at which $g_1(\mathbf u) - g_1(\mathbf v), \dots, g_n(\mathbf u) - g_n(\mathbf v)$ all vanish is exactly $\mathcal S(\RR)$, we have a separating set in the sense needed for this paper. If the set of $\CC$-points at which they all vanish is $\mathcal S(\CC)$, then we have what Dufresne \cite{dufresne} calls a {\em geometric separating set}. A geometric separating set is separating, but not conversely.

By the Nullstellensatz, it follows from the above that $g_1,\dots,g_n$ form a geometric separating set if and only if the ideal $I$ of $\RR[\mathbf u, \mathbf v]$ generated by the polynomials $g_1(\mathbf u) - g_1(\mathbf v), \dots, g_n(\mathbf u) - g_n(\mathbf v)$ has the same radical as the ideal $J$ generated by the polynomials $f(\mathbf u)-f (\mathbf v)$ for all $f\in \RxG$. Likewise, by the Real Nullstellensatz, it follows from the above that the $g_j$'s form a separating set if and only if $I$ has the same {\em real} radical as $J$. We are interested in the latter question. 

We write the real radical of an ideal $I$ as $\sqrt[\RR]{I}$. (See Appendix~\ref{app:comm-alg} for the definition of the real radical.) In view of the real radical membership test recently given by Reid, Wang, and Wolkowicz \cite{reid-wang-wolkowicz}, we can give the algorithm:

\begin{proof}[Proof of Proposition~\ref{prop:alg-separate}]
We are given a basis $g_1,\dots,g_n$ for $U$ and we want to know if it is a separating set for $\RxG$. Using Derksen's algorithm \cite{derksen1999} (see Appendix~\ref{app:generator-alg}), compute an algebra generating set $f_1,\dots,f_m$ for $\RxG$. Then consider the polynomial ring $\RR[\mathbf u,\mathbf v]$ in $2p$ unknowns $u_1,\dots, u_p, v_1,\dots, v_p$, and form the ideal
\[
J\defeq \langle f_1(\mathbf u) - f_1(\mathbf v), \dots, f_m(\mathbf u) - f_m(\mathbf v)\rangle_{\RR[\mathbf u, \mathbf v]}.
\]
We have $\mathcal S (\RR) = V_\RR(J)$ because the $f_1,\dots,f_m$ generate $\RxG$ as an algebra, thus two points $u,v$ of $V(\RR)$ are indistinguishable to all invariants $f$ if and only if they are indistinguishable to $f_1,\dots,f_m$, if and only if the generators of $J$ all vanish when evaluated at $(u,v)$. 

Form the ideal
\[
I\defeq \langle g_1(\mathbf u) - g_1(\mathbf v),\dots, g_n(\mathbf u) - g_n(\mathbf v)\rangle_{\RR[\mathbf u,\mathbf v]}.
\]
Then because $g_1,\dots,g_n \in \RxG = \RR[f_1,\dots,f_m]$, each $g_j(\mathbf u)-g_j(\mathbf v)$ vanishes on the entirety of $\mathcal S(\RR)$. Thus we certainly have $V_\RR(J) = \mathcal S(\RR) \subset V_\RR(I)$. The question is whether we have equality.

By the Real Nullstellensatz (e.g., \cite[Theorem~6.7]{lam1984}), $V_\RR(I)=V_\RR(J)$ if and only if $I$ and $J$ have the same real radical. It is enough to check if $J$ is contained in $I$'s real radical, because this would imply $\sqrt[\RR]{J} \subset \sqrt[\RR]{I}$, whereupon 
\[
V_\RR(J) = V_\RR(\sqrt[\RR]{J})\supset V_\RR(\sqrt[\RR]{I}) = V_\RR(I),
\]
and we already have the reverse inclusion.

Therefore, to test if $U$ is separating, it suffices to use the real radical membership test of Reid, Wang, and Wolkowicz \cite[Algorithm~1 and Remark~III.1]{reid-wang-wolkowicz} to see if $f_j(\mathbf u) - f_j(\mathbf v)$ lies in the real radical of $I$, for each $j=1,\dots,m$.
\end{proof}

\begin{remark}\label{rmk:sdp-numerics}
The real radical membership test of Reid, Wang, and Wolkowicz \cite{reid-wang-wolkowicz} relies on semidefinite programming. Without the hypothesized oracle, this potentially introduces numerical issues. They may be surmountable with current technology. The following approach was suggested to us by Fei Wang: There is recent progress of Henrion, Naldi, and Safey El Din on carrying out semidefinite programming tasks in exact symbolic arithmetic \cite{henrion2021}, but it is not designed to produce interior points. On the other hand, the method of \cite{henrion2021} can carry out the facial reduction steps in \cite{reid-wang-wolkowicz}, and we are left with the problem of finding an interior point in a full-rank feasible set in exact arithmetic. If an interior point method (e.g. \cite{benterki2011finding}) can produce an interior point whose distance from the boundary is controlled, then it can be rounded off to an exact interior point. These considerations go beyond our present scope.
\end{remark}

\begin{remark}
We add a word about the relationship between the algorithms of Propositions~\ref{prop:alg-separate} and \ref{prop:alg-worst-list}. They have the same underlying idea. In each case, there is a geometric object such that the desired condition (respectively ``$U$ is separating" and ``$\RxG$ is module-finite over $\RR[U]$") is detected by whether a certain ideal (respectively the $I\subset \RR[\mathbf u,\mathbf v]$ defined above, and $\langle U\rangle_{\Rx}\subset \Rx$) cuts it out. In the former case, this object is the set of real points of the separating variety $\mathcal S$ in $V\times V$; in the latter case, it is the {\em Hilbert nullcone} $\mathcal N_V\subset V$, defined to be the vanishing set in $V(\CC)$ of the positively graded ideal $\RxG_{>0}$ of $\RxG$. In both cases, a generating set for $\RxG$ allows one to write down polynomials (respectively the polynomials $f_j(\mathbf u) - f_j(\mathbf v)$, or the $f_j$'s themselves) that certainly cut out the desired object. Then, by the (respectively real or ordinary) Nullstellensatz, whether the proposed ideal (respectively $I$ or $\langle U\rangle_{\Rx}$) also cuts it out can be detected by a radical membership test---for real radicals, when we are concerned with real points (as in $\mathcal S(\RR)$), or for ordinary radicals, when we are concerned with all the points (as in $\mathcal N_V$).
\end{remark}

\subsection*{Acknowledgements}
The authors thank the anonymous referees for thoughtful feedback which greatly improved the article.

ASW is grateful for feedback given by several members of the audience when he presented a preliminary version of these results at a workshop dedicated to multi-reference alignment, organized as part of the Simons Collaboration on Algorithms \& Geometry.

BBS is grateful to Alena Pirutka for directing him to~\cite{pop}, and to Jason Bell, Sophie Marques, Daniel Litt, Larry Guth, and especially Fei Wang, Or Hershkovits, and Gregor Kemper for helpful conversations.

\bibliographystyle{alpha}
\bibliography{bib}

\appendix

\section{Algebra and invariant theory primer}\label{app:algprimer}

In an effort to make this paper accessible to as broad an audience as possible, in this appendix we collect together material from algebra, algebraic geometry, and invariant theory that is used in Sections~\ref{sec:algebraic} and \ref{sec:algebraic-proofs} and may be unfamiliar to applied mathematicians working in statistical signal processing.

\subsection{Fields}\label{app:fields}

This section collects together the language and results from field theory that are used in the paper.

\begin{definition}
If $F_1$ is a subfield of $F_2$, we write $F_2/F_1$ and call this a \emph{field extension}. The \emph{degree} of the extension, denoted $[F_2 : F_1]$, is the dimension of $F_2$ as a vector space over $F_1$. If $[F_2:F_1]<\infty$, the extension is \emph{finite} (regardless of whether $F_1,\,F_2$ are finite fields). A field $L$ contained in $F_2$ and containing $F_1$ is called a \emph{subextension}, and we write, ``$L/F_1$ is a subextension of $F_2/F_1$."
\end{definition}

\begin{remark}\label{rmk:finite-is-transitive}
It is straightforward to check that if $L/F_1$ is a subextension of $F_2/F_1$, then $[F_2:L][L:F_1]=[F_2:F_1]$ \cite[Theorem~4.2]{jacobson}. Thus, a finite extension of a finite extension is finite.
\end{remark}

\begin{definition}\label{def:field-and-ring-generation}
If $F$ is a field sitting inside some ambient field $\Omega$, and $S$ is any subset of $\Omega$, then the {\em ring} (or {\em algebra}) {\em generated over $F$ by $S$} is the smallest subring of $\Omega$ containing both $F$ and $S$. It is denoted by $F[S]$. It consists of polynomials in the elements of $S$ with coefficients in $F$. In the same situation, $F(S)$ denotes the {\em field generated over $F$ by $S$}, the smallest subfield of $\Omega$ containing both $F$ and $S$. It consists of rational functions in the elements of $S$ with coefficients in $F$. If $S=\{\alpha\}$ is a singleton, we also write $F[\alpha]$ and $F(\alpha)$ for $F[S]$ and $F(S)$, respectively.
\end{definition}

\begin{remark}
The field $\CC(\mathbf{x})$ of rational functions in the $p$ indeterminates $x_1,\dots,x_p$ with coefficients in $\CC$ serves as an ambient field for every ring and field occurring in the present work.
\end{remark}

\begin{definition}
In the situation of Definition~\ref{def:field-and-ring-generation}, if $\alpha\in \Omega$ is a root of a nonzero univariate polynomial with coefficients in $F$, it is said to be \emph{algebraic over $F$}. There is a unique monic polynomial with coefficients in $F$ of which $\alpha$ is a root and such that $\alpha$ is not the root of a lower-degree polynomial with coefficients in $F$; it is called the \emph{minimal polynomial} of $\alpha$.
\end{definition}

\begin{definition}\label{def:algebraically-independent}
If $F$ is contained in a ring $R$ ($R$ could be an ambient field, for example), and a tuple $(\alpha_1,\dots\alpha_r)$ of elements of $R$ causes a nonzero multivariate polynomial $f$ with coefficients in $F$ to vanish upon substitution (i.e., $f(\alpha_1,\dots,\alpha_r)=0$ in $R$), the elements $\alpha_1,\dots,\alpha_r$ are said to be \emph{algebraically dependent} over $F$, and the polynomial $f$ is said to be an \emph{algebraic relation} between them. If there is no such $f$, then $\alpha_1,\dots,\alpha_r$ are said to be \emph{algebraically independent}. 
\end{definition}

\begin{definition}
A basic fact is that given a field extension $F_2/F_1$, all maximal sets of elements of $F_2$ algebraically independent over $F_1$ have the same cardinality \cite[Theorem~8.35]{jacobsonii}. This common cardinality is the {\em transcendence degree} of $F_2$ over $F_1$, and any such maximal algebraically independent subset is called a {\em transcendence basis} for $F_2$ over $F_1$.
\end{definition}

\begin{remark}\label{rmk:finite-field-extension}
In general, if $\alpha$ is an element of an ambient field containing $F$, then the field $F(\alpha)$ may be bigger than the ring $F[\alpha]$; however, if $\alpha$ is algebraic over $F$, then they coincide. Furthermore, this happens if and only if the extension $F(\alpha)/F$ has finite degree. See \cite[Theorem~4.1]{jacobson} for more details.
\end{remark}

\begin{proposition}\label{prop:field-finite}
Let $\RR \subseteq F_1 \subseteq F_2$ with $F_2$ finitely generated (as a field) over $\RR$. Let $r$ be the transcendence degree of $F_2$ over $\RR$. The field extension $F_2/F_1$ has finite degree if and only if $F_1$ contains $r$ elements algebraically independent over $\RR$.
\end{proposition}
\begin{proof}
If $F_1$ contains $r$ algebraically independent elements then the extension $F_2/F_1$ is both algebraic and finitely generated; induction on the number of generators (in view of Remark~\ref{rmk:finite-field-extension} and Remark~\ref{rmk:finite-is-transitive}) shows that the total extension degree is finite. Otherwise, the extension is transcendental and has infinite degree, again by Remark~\ref{rmk:finite-field-extension}.
\end{proof}

\begin{definition}
The \emph{characteristic} of a field $F$ is the smallest natural number $N$ such that $1+1+\dots+1$ ($N$ times) equals $0$. If there is no such natural number, $F$ is said to be of \emph{characteristic zero}.
\end{definition}

\begin{remark}
The present work only considers fields of characteristic zero, with the passing exception of the discussions in Remark~\ref{rmk:nonmodular-field-gen} and Remark~\ref{rmk:pontryagin-other-fields}.
\end{remark}

\begin{definition}
Suppose $F_2/F_1$ is a field extension and $f$ a univariate polynomial with coefficients in $F_1$ such that $f$ factors into linear factors over $F_2$, and the roots of $f$ in $F_2$ generate the latter over $F_1$. Then $F_2$ is a \emph{splitting field} for $f$ over $F_1$.
\end{definition}

\begin{remark}
Given a univariate polynomial $f$ over a field $F$, a splitting field for $f$ over $F$ always exists \cite[Theorem~4.3]{jacobson}.
\end{remark}

\begin{definition}
An algebraic field extension $F_2/F_1$ is said to be \emph{separable} if each $\alpha\in F_2$ is a root of an irreducible polynomial over $F_1$ all of whose roots (in a splitting field for that polynomial over $F_1$) are distinct.
\end{definition}

\begin{remark}
Separability is only a question in positive characteristic. The reader interested in the main setting of this work can safely assume every field extension is separable. (This is because a multiple root of a polynomial is exactly a root shared with its derivative, but in characteristic zero, an irreducible polynomial is necessarily prime to its derivative because the latter is of lower degree. In positive characteristic, the derivative can be zero.)
\end{remark}

\begin{proposition}[Primitive element theorem]\label{prop:primitive-element}
Any finite separable field extension $F_2/F_1$ is generated by a single element $\alpha\in F_2$.
\end{proposition}

\begin{proof}
See \cite[Corollary to Theorem~4.28]{jacobson}.
\end{proof}

The primitive element theorem is used in the algebraic geometry-free proof of Theorem~\ref{thm:bound-list} in Section~\ref{sec:generic-proof}.

\begin{definition}
A finite field extension $F_2/F_1$ is said to be \emph{Galois} if it is separable and $F_2$ is a splitting field for some polynomial over $F_1$.
\end{definition}

\begin{remark}
Every finite separable field extension $F_2/F_1$ is contained in a Galois extension $L/F_1$. If $\alpha$ is a generator for $F_2$ over $F_1$, and $f$ is its minimal polynomial, then $L$ can be taken to be a splitting field for $f$ over $F_1$.
\end{remark}

\begin{definition}\label{def:galois-group}
Given a finite Galois extension $F_2/F_1$, the group of automorphisms of $F_2$ that fix $F_1$ pointwise is called the {\em Galois group} of $F_2/F_1$, and denoted $\operatorname{Gal}(F_2/F_1)$.
\end{definition}

\begin{proposition}[Fundamental theorem of Galois theory]\label{prop:fund-thm-galois-thy}
In the situation of Definition~\ref{def:galois-group}, the group $\Gamma \defeq\operatorname{Gal}(F_2/F_1)$ is finite and of order equal to the extension degree $[F_2:F_1]$. Furthermore, there is a bijective, inclusion reversing correspondence between the subgroups of $\Gamma$ and the subextensions of $F_2/F_1$. The correspondence sends a subextension $L/F_1$ to the subgroup of $\Gamma$ whose elements fix each element of $L$ pointwise, and sends a subgroup $H$ to the set of elements of $F_2$ that each element of $H$ fixes pointwise.
\end{proposition}

\begin{proof}
See \cite[p.~240]{jacobson}.
\end{proof}

The fundamental theorem of Galois theory is used in the proof of Theorem~\ref{thm:reg-rep-galois}.

\subsection{Some algebraic geometry}\label{app:someAG}

In this section we give an informal overview of the construction of the parameter spaces $X$ and $Y$ and the maps $\pi$ and $\varphi$ used in the proofs of Theorems~\ref{thm:generic-list} and \ref{thm:bound-list}, and Propositions~\ref{prop:worst-case-unique} and \ref{prop:worst-list} given in Section \ref{sec:algebraic}. It is intended for readers without an algebraic geometry background but interested in those proofs. Systematic references for this material include \cite{mumford},   \cite{shaf}, and Chapter 1 of \cite{springer}.

We stress that this discussion is informal. The construction we are about to give of the algebraic variety associated to a finitely generated subalgebra of $\Rx$ elides the distinction between an affine $\RR$-variety and its base change to $\CC$, for example, as well as between the variety and a particular embedding of it in some $\CC^N$. Nonetheless, the picture painted here is adequate to support the argumentation used in the present work. If one visualizes affine $\RR$-varieties as vanishing sets in $\CC^N$ of polynomials whose coefficients lie in $\RR$, and keeps track of (i) which points have only real coordinates (i.e., the $\RR$-points), and (ii) which $\CC$-valued functions on the ambient space $\CC^N$ are given by polynomials with real coefficients, then one is storing all the relevant data.

Let $A$ be a finitely generated subalgebra of $\Rx$, for example $\RxG$ or $\RR[U]$. Choose a finite set of generators $f_1,\dots,f_N$ for $A$ over $\RR$. In general, these will not be algebraically independent: there will be many polynomials $g\in \RR[x_1,\dots,x_N]$ such that
\[
g(f_1,\dots,f_N) = 0
\]
as polynomials. Let $X_A$ be the common zero set in $\CC^N$ of all of these polynomials $g$. This is the {\em affine $\RR$-variety} associated to the algebra $A$. The algebra $A$ is, in turn, the {\em coordinate ring} of the variety $X_A$.

The polynomial map $\CC^p \rightarrow \CC^N$ given by sending a point to the $N$-tuple obtained by evaluating $f_1,\dots,f_N$ at that point,
\[
(x_1,\dots,x_p)\mapsto (f_1(x_1,\dots,x_p),\dots,f_N(x_1,\dots,x_p)),
\]
has image contained in $X_A$. This is the {\em morphism of varieties induced by} the ring inclusion $A\subset \Rx$. 

Note that, because the coefficients of $f_1,\dots,f_N$ are real, the image of $\RR^p$ under this map lands in $X_A(\RR)$, the set of points of $X_A$ with real coordinates. 

Taking $A$ as $\RxG$, this construction yields the variety $X$ and the map $\pi:V\rightarrow X$ used in the proofs in Section~\ref{sec:algebraic}. Taking $A$ as $\RR[U]$, we get the variety $Y$. The map $\varphi:X\rightarrow Y$ will be discussed shortly.

This construction depends on a choice of generators $f_1,\dots,f_N$ for $A$ over $\RR$, which will in general not be unique. Not even the number $N$ of generators is uniquely determined by $A$. However, given two different choices of generators, the two different zero sets yielded by the above construction---call them $X_A\subset \CC^N$ and $X'_A\subset \CC^{N'}$---will be isomorphic in the sense that there will exist polynomial maps $\Phi:\CC^N\rightarrow \CC^{N'}$ and $\Psi:\CC^{N'}\rightarrow \CC^N$ that induce a bijection of $X_A$ with $X'_A$. It is because of this that we can speak of the $X_A$ yielded by the above construction as ``the" variety associated to the algebra $A$. Different choices of generators can be seen as yielding different embeddings of $X_A$ in various $\CC^N$'s.

One obtains the maps $\Phi,\Psi$ of the last paragraph by expressing each of the sets of generators for $A$ (used to construct $X_A$ and $X'_A$) as polynomials in the other set. Since all of these generators lie within the $\RR$-algebra $A$, $\Phi$ and $\Psi$ will be defined by polynomials with coefficients in $\RR$, so they will also biject the $\RR$-points $X_A(\RR)=X_A\cap \RR^N$ of $X_A$ with the $\RR$-points $X'_A(\RR) = X'_A\cap \RR^{N'}$ of $X'_A$. Therefore, the set of $\RR$-points of the variety $X_A$ is determined uniquely up to isomorphism, regardless of the embedding determined by a choice of generators.

The map $\varphi:X\rightarrow Y$ used in Section~\ref{sec:algebraic} is obtained as follows. Suppose we have chosen a set of generators $f_1,\dots,f_m$ for $\RxG$, yielding the construction of $X$ as a subset of $\CC^m$, as above, and also a set of generators $g_1,\dots,g_n$ for $\RR[U]$, yielding $Y$ as a subset of $\CC^n$. Because $\RR[U]$ is contained in $\RxG$, which is generated by $f_1,\dots,f_m$, there exist some polynomials $G_1,\dots,G_n$ in $\RR[x_1,\dots,x_m]$ such that
\[
g_i = G_i(f_1,\dots,f_m)
\]
for each $i=1,\dots,n$. Then the polynomial map $\CC^m\rightarrow\CC^n$ given by
\[
(x_1,\dots,x_m)\mapsto (G_1(x_1,\dots,x_m),\dots,G_n(x_1,\dots,x_m))
\]
restricts to a map $X\rightarrow \CC^n$ whose image lands inside $Y$. This is the map $\varphi$ used in Section~\ref{sec:algebraic}. Furthermore, the explicit description of it given here is also used in Section~\ref{sec:generic-converse} to give a differential-geometric proof of the second part of Theorem~\ref{thm:generic-list}. One says that $\varphi$ is the {\em morphism of varieties induced by} the ring inclusion $\RR[U]\subset \RxG$.

We now describe the topology on varieties used in algebraic geometry.

\begin{definition}\label{def:zariski}
The {\em Zariski topology on $\CC^N$} is the topology with basis the nonvanishing sets of polynomial functions on $\CC^N$. For any variety $X_A\subset \CC^N$ embedded in $\CC^N$ as above, the {\em Zariski topology on $X_A$} is the subspace topology induced by the Zariski topology on $\CC^N$. Similarly, the {\em Zariski topology on $\RR^N$} is the subspace topology induced from the Zarkiski topology on $\CC^N$ via $\RR^N$'s natural inclusion in $\CC^N$.
\end{definition}

For example, because polynomials on $\CC$ vanish in at most finitely many points, the Zariski topology on $\CC^1$ (viewed as an affine algebraic variety) has as open sets exactly the cofinite sets and the empty set.

\begin{remark}
A comment on the Zariski topology: it is extremely coarse. The open sets are huge: since polynomial vanishing sets have Lebesgue measure zero, all nonempty open sets have full measure---a fact we use several times. Thus Zariski-density does not conform to intuition trained on the Euclidean topology about the meaning of the word ``density". For example, a set may be both Zariski-dense and discrete for the Euclidean topology. To illustrate, in $\CC^1$ , the ordinary integers form a Zariski-dense subset. However, a Zariski-dense subset of a positive-dimensional $\RR$-variety is necessarily {\em infinite}, because if it were finite, one could contain it in the vanishing set of a polynomial. This is enough for our purposes.
\end{remark}

With the Zariski topology defined, we can explain what is meant by saying that $X,Y$ are irreducible and $\varphi,\psi$ are dominant.

An $\RR$-variety is {\em geometrically irreducible} if its set of $\CC$-points is not the union of proper Zariski-closed subsets. It is {\em irreducible} if this set is not the union of proper Zariski-closed subsets that are stable under the action of complex conjugation on the coordinates. Clearly geometric irreducibility implies irreducibility.

For $A$ a finitely generated subalgebra of $\Rx$, the associated variety $X_A$ defined above is always geometrically irreducible. This is because geometric irreducibility of $X_A$ is equivalent to the complexification of the coordinate ring $A$ being an integral domain; see, e.g., \cite[pp.~34--35]{shaf}. (An {\em integral domain} is a ring with the property that $fg=0$ implies $f=0$ or $g=0$.) As $\Cx$ is an integral domain, all its subrings are as well. 
It follows in particular that the $X$ and $Y$ of Section~\ref{sec:algebraic} are geometrically irreducible, and therefore irreducible.

A map of algebraic varieties, such as $\varphi$ or $\pi$ above, is said to be {\em dominant} if the image of (the set of $\CC$-points of) the domain is dense in the codomain, with respect to the Zariski topology. This is equivalent to the associated ring map being injective; see, e.g., \cite[Lemma~1.9.1(ii)]{springer}. Because $\RR[U]\subset \RxG$ and $\RxG\subset\Rx$ are ring inclusions, $\varphi$ and $\pi$ are both dominant.

\subsection{Some commutative algebra}\label{app:comm-alg}

This section draws out the principles of commutative algebra used in various proofs in Section~\ref{sec:algebraic}.

\subsubsection{Height, Krull dimension, integrality, Noether normalization}

The following ideas are used in the proof of Proposition~\ref{prop:trdeg-het}.

\begin{definition}\label{def:height}\label{def:krull-dim}
The {\em height} of a prime ideal $\mathfrak{p}$ in a commutative ring $A$ is the length $h$ of the longest strictly ascending chain of prime ideals 
\[
\mathfrak{p}_0\subset\mathfrak{p}_1\subset\dots\subset\mathfrak{p}_h = \mathfrak{p}
\]
contained in $\mathfrak{p}$. The {\em height} of an arbitrary ideal $I$ is the infimum of the heights of all prime ideals containing $I$. The {\em Krull dimension} of a ring $A$ is the supremum of the heights of all prime ideals of $A$. It is written $\dim_{\mathrm{Krull}} A$, or just $\dim A$ when the meaning is clear in context.
\end{definition}

A field has Krull dimension zero, since its only prime ideal is the zero ideal. A polynomial ring such as $\Rx = \RR[x_1,\dots,x_p]$ has Krull dimension $p$, and $\langle 0\rangle\subset \langle x_1\rangle\subset \langle x_1,x_2\rangle\subset\dots\subset \langle x_1,\dots,x_p \rangle$ is a chain of prime ideals realizing the maximum length. In the notation of Appendix~\ref{app:someAG}, the Krull dimension of a finitely generated subalgebra $A$ of $\Rx$ conforms to geometric intuition about the dimension of the associated variety $X_A$; specifically, $X_A$'s nonsingular locus is a complex manifold of complex dimension equal to $\dim_{\mathrm{Krull}} A$. The definition has the same idea behind it as defining the dimension of a vector space as the length of a maximal flag.

(Note that Krull dimension of an $\RR$-algebra does not coincide with dimension as an $\RR$-vector space. Indeed, a commutative $\RR$-algebra that is finite-dimensional as an $\RR$-vector space must have Krull dimension zero.)

\begin{definition}
A ring extension $A\subset B$ is {\em integral} if every element of $B$ satisfies a monic polynomial with coefficients in $A$.
\end{definition}

Integral extensions preserve Krull dimension:

\begin{lemma}
If a ring extension $A\subset B$ is integral, then $\dim_{\mathrm{Krull}}A = \dim_{\mathrm{Krull}}B$. 
\end{lemma}

\noindent This is a consequence of standard facts about integral extensions, see for example \cite[Section~4.4]{Eis}. (The fact that $\dim_{\mathrm{Krull}}A\geq \dim_{\mathrm{Krull}}B$ follows from the incomparability theorem \cite[Corollary~4.18]{Eis}, while the fact that $\dim_{\mathrm{Krull}}A\leq \dim_{\mathrm{Krull}}B$ follows from lying-over and going-up \cite[Proposition~4.15]{Eis}.)

Given a finite group $G$ acting on a ring $A$, $A$ is always an integral extension of the invariant ring $A^G$, because for any $f\in A$, the polynomial $\prod_{g\in G}(X - g(f))$ is monic with coefficients in $A^G$, and $f$ satisfies it. Thus taking invariants by a finite group never affects Krull dimension.

Many questions about the Krull dimension of a finitely generated algebra over a field are illuminated by the {\em Noether normalization lemma}:

\begin{proposition}[Noether normalization lemma]
If $A$ is a finitely generated algebra over a field $k$, of Krull dimension $d$, then there exist $d$ algebraically independent elements $f_1,\dots,f_d\in A$, such that $A$ is finitely generated as a module over the subring $k[f_1,\dots,f_d]$.
\end{proposition}

See \cite[Chapter I, \S~5, Theorem 10]{shaf} for a geometric argument or \cite[Theorem~13.3]{Eis} for a sharper version. Note that the subring $k[f_1,\dots,f_d]$ is (isomorphically) a polynomial ring, due to the algebraic independence of $f_1,\dots,f_d$.

Module-finite extensions are automatically integral, thus the Noether normalization lemma guarantees that any finitely generated algebra over a field is integral over a polynomial ring. It follows from this that Krull dimension is equal to transcendence degree for such a ring.

If $A$ and $B$ are two finitely generated $k$-algebras, say of Krull dimensions $d$ and $e$, then by Noether normalization, they are module-finite over polynomial rings in $d$ and $e$ indeterminates respectively. Then the tensor product $A\otimes_k B$ is module-finite (and therefore integral) over the tensor product of the two polynomial rings. But the tensor product of polynomial rings $k[x_1,\dots,x_d]$ and $k[y_1,\dots,y_e]$ over $k$ is exactly the polynomial ring $k[x_1,\dots,x_d,y_1,\dots,y_e]$, of Krull dimension $d+e$. As $A\otimes_kB$ is integral over this, it has the same dimension. In this way, we see that taking the tensor product of finitely generated $k$-algebras has the effect of adding the dimensions. (In the ring-variety correspondence described in Appendix~\ref{app:someAG}, tensor products of algebras correspond to products of the associated varieties, see for example \cite[pp.~24--25]{shaf}. Then the fact that taking the tensor product adds Krull dimensions corresponds geometrically to the fact that the product of varieties of dimension $d$ and $e$ has dimension $d+e$.)

\subsubsection{Invariant rings of normal subgroups}

We mention one more principle used in the proof of Proposition~\ref{prop:trdeg-het}. If $A$ is a ring with an action of a group $G$, and $N$ is a normal subgroup of $G$, then the invariant ring $A^N$ of $N$ is stable as a set under the action of $G$, because for $f\in A^N$, $g\in G$, and $n\in N$, we have $ng(f) = gn'(f)$ for some $n'\in N$ (as $N$ is normal), and $gn'(f) = g(f)$ because $f\in A^N$; therefore $g(f)$ is invariant under the action of $N$, i.e., $g(f)\in A^N$. Furthermore, because $N$ acts trivially on $A^N$, the action of $G$ on $A^N$ induces an action of the quotient group $G/N$, and invariance of an element under $G$ is equivalent to invariance under $N$ (i.e., membership in $A^N$) plus invariance under the action of $G/N$ on $A^N$. In other words, 
\[
A^G = \left(A^N\right)^{G/N}.
\]

\subsubsection{Radical, real radical, Nullstellensatz}

The following ideas are used in the proofs in Sections~\ref{sec:worst-unique} and \ref{sec:separating-alg}.

\begin{definition}
Given an ideal $I$ in a ring $A$, the {\em radical} of $I$, written $\sqrt{I}$, is the ideal
\[
\{f\in A : f^k\in I\text{ for some }k\in\mathbb{N}\}.
\]
An ideal is {\em radical} if it is equal to its own radical.
\end{definition}

\noindent The radical is a {\em closure operation}, i.e., it contains $I$, it respects containments, and $\sqrt{\sqrt{I}} = \sqrt{I}$. If the ring in question is $\Rx$ or $\Cx$, the radical plays a decisive role in understanding the {\em vanishing set}
\[
V(I) \defeq \{ x\in \CC^p: f(x)=0\text{ for all }f\in I\}.
\]
One begins to intuit this because, while in general, enlarging an ideal $I$ will cause the vanishing set $V(I)$ to shrink, passage from $I$ to $\sqrt{I}$ does not have this effect: one has $V(I) = V(\sqrt{I})$ because if $f^k$ vanishes at $x$ for some $k$, then $f$ itself vanishes at $x$. In fact, the relationship is very tight:

\begin{theorem}[Hilbert's Nullstellensatz]
The ideal $I(V(I))$ of $\Rx$ consisting of polynomials that vanish on the set $V(I)\subset\CC^p$ is exactly $\sqrt{I}$.
\end{theorem}

It follows from this that two ideals have the same vanishing set if and only if they have the same radical. This theorem is foundational in algebraic geometry: it is the starting point of the algebra-geometry duality described in Appendix~\ref{app:someAG}. For a proof, see \cite[pp.~281--282]{shaf}. 

There is an analogous concept, and an analogous theorem, when one wants to focus on the set of {\em real} points where some polynomials vanish:

\begin{definition}
Given an ideal $I$ in an $\RR$-algebra $A$, the {\em real radical} of $I$ is the ideal
\[
\sqrt[\RR]{I} \defeq \{f\in A: f^{2k}+\sigma\in I\text{ for some }k\in\mathbb{N}\text{ and }\sigma\in\Sigma A^2\},
\]
where $\Sigma A^2$ denotes the set of elements of $A$ that are finite sums of squares.
\end{definition}

\noindent Given an ideal $I\subset \Rx$, we define the {\em real vanishing set}
\[
V_\RR(I) \defeq\{ x\in \RR^p : f(x)=0\text{ for all }f\in I\}
\]
in analogy with the (complex) vanishing set $V(I)$.

\begin{theorem}[Krivine-Dubois-Risler Real Nullstellensatz]
If $I$ is an ideal in $\Rx$, then the ideal $I(V_\RR(I))$ consisting of polynomials that vanish on the set $V_\RR(I)\subset \RR^p$ is exactly $\sqrt[\RR]{I}$.
\end{theorem}

See \cite[Section~6]{lam1984} for a proof. As with the ordinary Nullstellensatz, it follows that two ideals define the same real vanishing set if and only if they have the same real radical.

\subsection{Some invariant theory}\label{app:invthybasics}

The following basic facts from invariant theory were exploited in Section~\ref{sec:algebraic} (see Setup~\ref{rmk:AGsetup} and \ref{rmk:separation-in-terms-of-AG}).

\begin{theorem}\label{thm:finite-gen}
The invariant ring $\RxG$ is finitely generated as an $\RR$-algebra. In other words, there exist generators $f_1, \ldots, f_m \in \RxG$ such that $\RR[f_1,\ldots,f_m] = \RxG$.
\end{theorem}
The following is a sketch of the standard argument (see for example \cite[Theorem~2.2.10]{dk-book}, \cite[Theorem~4-1.3]{kac-notes}).
\begin{proof}[Proof sketch]
Because $\Rx$ is a noetherian ring by the Hilbert basis theorem, the ideal $I$ of $\Rx$ generated by the homogeneous $G$-invariants of positive degree is in fact generated by a finite subset $f_1,\dots,f_m$ of these. This finite subset is actually also an algebra generating set for $\RxG$. This is seen by an induction argument. It is evidently sufficient to show that an arbitrary homogeneous element $f\in \RxG$ of positive degree lies in $\RR[f_1,\dots,f_m]$. One first represents $f$ as a linear combination of $f_1,\dots,f_m$ with coefficients in $\Rx$, possible because $f\in I$. Without loss of generality, the coefficients can be taken to be homogeneous because $f$ and the $f_i$ are. Then one applies the Reynolds operator to this expression, obtaining an expression for $f$ as a linear combination of the $f_i$ with {\em invariant} homogeneous coefficients. Because the $f_i$ have positive degrees, the degrees of the coefficients are lower than that of $f$. Then one applies induction on the degree, so the coefficients lie in $\RR[f_1,\dots,f_m]$. This completes the proof.
\end{proof}

\noindent Furthermore, there is an algorithm to find a generating set; see Appendix~\ref{app:generator-alg}.

\begin{theorem}\label{thm:RxG-separates}
The full invariant ring $\RxG$ resolves every $\theta \in V$.
\end{theorem}
This was mentioned in Setup~\ref{rmk:separation-in-terms-of-AG}, were we cited \cite[Chapter~3, \S~4, Theorem~3]{on-vin}. The following proof is taken from \cite[Theorem~6-2.2]{kac-notes}.
\begin{proof}
Let $\mathfrak{o}_1, \mathfrak{o}_2 \in V/G$ be distinct (and therefore disjoint) orbits. Since $G$ is compact and acts continuously, $\mathfrak{o}_1$ and $\mathfrak{o}_2$ are compact subsets of $V$. Thus by Urysohn's lemma there exists a continuous function $\tilde f: V \to \RR$ such that $\tilde f(\tau) = 0 \;\; \forall \tau \in \mathfrak{o}_1$ and $\tilde f(\tau) = 1 \;\; \forall \tau \in \mathfrak{o}_2$. The Stone--Weierstrass theorem states that a continuous function on a compact domain can be uniformly approximated to arbitrary accuracy by a polynomial. This means there is a polynomial $f \in \Rx$ with $f(\tau) \le 1/3 \;\; \forall \tau \in \mathfrak{o}_1$ and $f(\tau) \ge 2/3 \;\; \forall \tau \in \mathfrak{o}_2$. Applying the Reynolds operator to $f$, one obtains an invariant polynomial $h = \mathcal{R}(f)$ that separates the two orbits: $h(\mathfrak{o}_1) \le 1/3$ and $h(\mathfrak{o}_2) \ge 2/3$.
\end{proof}

This does not hold when the group $G$ is not compact---the standard example is given in footnote~\ref{note:noncompact-do-not-separate}.

\subsection{Gr\"{o}bner bases}\label{app:grobner}

In this section we show how to use Gr\"obner bases to test various algebraic conditions. The ideas from this section are mostly standard in the theory of Gr\"obner bases; see, e.g.,\ \cite{iva} for a reference. We also (in \ref{app:generator-alg}) describe the algorithm of Derksen for computing generators for the invariant ring $\RxG$, which relies on Gr\"obner bases. The ideas in this section are used in Section~\ref{sec:algebraic} to prove many of the algorithmic results.

\begin{definition}\label{def:monomial-order}
A \emph{monomial order} on $\RR[{\bf x}]$ is a well-ordering on the set $\mathcal{M}$ of all (monic) monomials, satisfying $M \le N \Leftrightarrow MP \le NP$ for all $M,N,P \in \mathcal{M}$. We say that a monomial order \emph{favors} an indeterminate $x_i$ \emph{over} a subset $S$ of other indeterminates if the monomial $x_i$ is larger (with respect to the monomial order) than any monomial supported on $S$; we say it \emph{favors} $x_i$ (no qualification) if it favors $x_i$ over all other indeterminates. We write $\mathrm{LM}(f)$ to denote the leading monomial of a polynomial $f$, i.e.,\ the monomial occurring in $f$ that is largest with respect to the monomial order; $\mathrm{LM}(f)$ does not include the coefficient.
\end{definition}

\begin{definition}
A \emph{Gr\"obner basis} of an ideal $I \subseteq \RR[{\bf x}]$ is a finite subset $B \subseteq I$ such that for every $f \in I$ there exists $b \in B$ such that $\mathrm{LM}(f)$ is a multiple of $\mathrm{LM}(b)$. A Gr\"{o}bner basis $B$ is \emph{reduced} if all its elements are monic and it has the additional property that for every pair of distinct $b,b' \in B$, no monomial occurring in $b$ is a multiple of $\mathrm{LM}(b')$.
\end{definition}

\noindent The following basic facts about Gr\"obner bases are proved in \cite{iva}. A Gr\"obner basis is indeed a basis, in that it generates the ideal. Every ideal $I \subseteq \Rx$ has a Gr\"obner basis, and has a unique reduced Gr\"obner basis. \emph{Buchberger's algorithm} computes the reduced Gr\"obner basis of an ideal $I = \langle f_1, \ldots, f_m \rangle$, given a list of generators $f_i$. (It is not a polynomial-time algorithm, however.)

\begin{paragraph}{Membership in an ideal.}
Suppose we want to know whether a given polynomial $f$ is contained in an ideal $I$ of $\Rx$ generated by some other given polynomials. This situation generalizes the question of whether a given univariate polynomial is a multiple of another given univariate polynomial, which can be solved with simple polynomial division. The following algorithm generalizes the polynomial division algorithm to the multivariate setting. 

Compute a Gr\"{o}bner basis $B=b_1,\dots,b_m$ of $I$. Note that by the definition of a Gr\"{o}bner basis, it is impossible for $f$ to be in $I$ unless $LM(f)$ is a multiple of some $LM(b_j)$. Proceed as follows: if $LM(f)$ is a multiple of any $LM(b_j)$, then compute $f-cb_j$, where $c$ is an appropriate scalar so that the leading term of $f$ is canceled by $cb_j$. Update $f\defeq f-cb_j$. Repeat until either $f=0$ (in which case the original $f$ was in $I$) or $LM(f)$ is not a multiple of any $LM(b_j)$ (in which case the original $f$ was not in $I$). The algorithm terminates because at each stage, $f-cb_j$ has a lower leading monomial than $f$ with respect to the chosen monomial order, and the latter is a well-order.
\end{paragraph}

Suppose we are interested in the relations between polynomials $f_1,\ldots,f_m \in \RR[{\bf x}]$. Introduce additional variables ${\bf t} = (t_1,\ldots,t_m)$ and consider the ideal $I \defeq \langle f_1({\bf x}) - t_1, \ldots, f_m({\bf x}) - t_m \rangle \subseteq \RR[{\bf x}, {\bf t}]$. Given $f_1,\ldots,f_m$ there is an algorithm to compute a Gr\"obner basis for the \emph{elimination ideal}
$$J \defeq \langle f_1({\bf x}) - t_1, \ldots, f_m({\bf x}) - t_m \rangle \cap \RR[{\bf t}].$$
In fact, the algorithm is simply to compute a Gr\"obner basis for $I$ using a particular monomial order and then keep only the elements that depend only on ${\bf t}$ (see Chapter~3 of \cite{iva}). The elimination ideal consists precisely of the polynomial relations among $f_1,\ldots,f_m$:

\begin{lemma}\label{lem:rel}
For any polynomial $P \in \RR[{\bf t}]$ we have: $P \in J$ if and only if $P(f_1({\bf x}), \ldots, f_m({\bf x})) \equiv 0$.
\end{lemma}
\begin{proof}
The direction `$\Rightarrow$' is clear because if we let $t_i = f_i({\bf x})$ for all $i$ then the generators of $I$ vanish and so every element of $I$ vanishes. To show the converse, it suffices to show that for any polynomial $P \in \RR[{\bf t}]$, $P(f_1({\bf x}), \ldots, f_m({\bf x})) - P(t_1,\ldots,t_m) \in I$. This can be shown inductively using the following key idea:
$$X_1 X_2 - t_1 t_2 = (X_1 - t_1)X_2 + (X_2 - t_2)t_1,$$
and so $X_1 X_2 - t_1 t_2 \in \langle X_1 - t_1, X_2 - t_2 \rangle$. By induction, a product of $X_i$'s minus a corresponding product of $t_i$'s lies in $\langle X_i-t_i\rangle_{i=1,\dots,m}$. Since ideals are vector spaces over the ground field, it follows that $P(X_1,\ldots, X_m) - P(t_1,\ldots,t_m)\in \langle X_i-t_i\rangle_{i=1,\dots,m}$ for any polynomial $P$, and substituting $f_i(\bf x)$ for $X_i$ (for each $i$) then leads to the desired conclusion.
\end{proof}

\begin{paragraph}{Membership in an $\RR$-algebra.}
Suppose we want to know whether $f_m \in \RR[f_1,\ldots,f_{m-1}]$. This is equivalent to asking whether there exists $P \in J$ of the form
\begin{equation}\label{eq:grob-alg}
P({\bf t}) = t_m - Q(t_1,\ldots,t_{m-1})
\end{equation}
for some $Q \in \RR[t_1,\ldots,t_{m-1}]$. Suppose that $J$ contains an element $P$ of the form (\ref{eq:grob-alg}). Compute a Gr\"obner basis $B$ for $J$ with respect to a monomial order that favors $t_m$. The leading monomial of $P$ is $t_m$ so by the definition of a Gr\"obner basis there must be an element $b \in B$ whose leading monomial divides $t_m$. Since $1 \notin J$ (by Lemma~\ref{lem:rel}), the leading monomial of $b$ is exactly $t_m$ and so $b$ takes the form (\ref{eq:grob-alg}). Therefore, $f_m \in \RR[f_1,\ldots,f_{m-1}]$ if and only if $B$ contains an element of the form (\ref{eq:grob-alg}).
\end{paragraph}

\begin{paragraph}{Membership in a field.}
Suppose we want to know whether $f_m \in \RR(f_1,\ldots,f_{m-1})$. This is equivalent to asking whether $f_m$ can be expressed as a rational function of $f_1,\ldots,f_{m-1}$ (with coefficients in $\RR$), which is equivalent (by multiplying through by the denominator) to asking whether there exists $P \in J$ of the form
\begin{equation}\label{eq:grob-field}
P({\bf t}) = t_m Q_1(t_1,\ldots,t_{m-1}) - Q_2(t_1,\ldots,t_{m-1}) \quad\text{with } Q_1 \notin J.
\end{equation}
Suppose that $J$ contains some element $P$ of the form (\ref{eq:grob-field}). Compute a \emph{reduced} Gr\"obner basis $B$ for $J$ with respect to a monomial order that favors $t_m$; we claim that an element of the form \eqref{eq:grob-field} must occur in $B$, so that $f_m$'s membership in the field generated by $f_1,\dots,f_{m-1}$ is detected by computing a reduced Gr\"obner basis for $J$ and looking for elements of the form \eqref{eq:grob-field}. We see this as follows:

It is a basic property of Gr\"obner bases that $P$ can be written as
$$P({\bf t}) = \sum_i p_i({\bf t}) b_i({\bf t})$$
where $p_i \in \RR[{\bf t}]$ and $b_i \in B$ with $\mathrm{LM}(p_i) \le \mathrm{LM}(P)$ and $\mathrm{LM}(b_i) \le \mathrm{LM}(P)$. Because $P$'s leading term has degree $1$ in $x_m$, any $b_i$'s that occur in this expression cannot be of degree greater than $1$ in $x_m$. We claim at least one of them has degree exactly $1$. For if none of them involve $x_m$ at all, then they are all homogeneous with respect to the grading that assigns degree $1$ to $x_m$ and degree zero to all other indeterminates; thus the ideal they generate is graded with respect to this grading. It follows that $t_mQ_1$ and $Q_2$ are separately contained in this ideal (and therefore in $J$) as they are homogeneous with respect to this grading. But then in an expression $t_mQ_1 = \sum q_ib_i$ for $t_mQ_1$, each term on the right side can be taken homogeneous with respect to this same grading, and then each $q_i$ must be divisible by $t_m$. Dividing it out, we obtain $Q_1=\sum q_i'b_i$, thus $Q_1\in J$. This is a contradiction of the assumption \eqref{eq:grob-field} about $P$. It follows that some $b_i$ is degree exactly $1$ in $t_m$, so can be written $b_i = t_mQ_1'-Q_2'$ with $Q_1',Q_2'$ polynomials in $t_1,\dots,t_{m-1}$. If $Q_1'\in J$, then some $b_j$ must have leading monomial dividing the leading monomial of $Q_1'$ (and $j\neq i$ because $b_i$ is the wrong degree in $t_m$ for this); but then $b_j$ also divides the leading term of $b_i$, contradicting that $B$ is reduced. Therefore $b_i$ has the claimed form \eqref{eq:grob-field}.
\end{paragraph}

\begin{remark}
An alternative to using Gr\"obner bases for algebra and field membership is to solve a (very large) linear system in order to find the minimal relation among a set of polynomials. There are bounds on the maximum possible degree of such a relation (if one exists) \cite{complexity-ann}.
\end{remark}

\subsubsection{Algorithm for generators of $\RxG$.}
\label{app:generator-alg}

We briefly describe the beautiful algorithm due to Derksen \cite{derksen1999}, as presented by Derksen and Kemper \cite{dk-book}, to compute a generating set for $\RxG$. As discussed in Setup~\ref{rmk:X(C)etc}, $G$ is necessarily a linear algebraic group defined over $\RR$. Furthermore, the complexification $G_\CC$ of $G$ is reductive, and it has the same invariant ring as $G$ up to base change, i.e., 
\[
\CC[\mathbf{x}]^{G_\CC} = \CxG \cong \CC\otimes_\RR\RxG.
\]
In characteristic zero, reductive implies linearly reductive. Thus \cite[Algorithm~4.1.9]{dk-book} can be applied in the present situation to obtain a generating set for $\CxG$ over $\CC$ (and, since all the computations happen over $\RR$, it will in fact provide a generating set for $\RxG$ over $\RR$). The algorithm takes as input (i) an explicit description of $G$ as a linear algebraic group, given by generators $h_1,\dots,h_t$ for an ideal in a polynomial ring $\RR[z_1,\dots,z_\ell]$ that defines $G$ as an affine $\RR$-variety, (ii) a description of the action of $G$ on $V$ as a $d\times d$ matrix whose entries $a_{ij}$ are polynomials in the coordinate functions $z_1,\dots,z_\ell$ on $G$, and (iii) a means to compute the Reynolds operator (Definition~\ref{def:reynolds}). 

\begin{example}
To illustrate the data (i) and (ii), consider the group $G=O(2)$. It consists of matrices 
\[
P = \begin{pmatrix} x&z \\ y&w\end{pmatrix}
\]
satisfying the equations $x^2+y^2=1$, $z^2+w^2=1$, and $xz+yw=0$; thus we can take $G$ to be defined by the ideal 
\[
\langle x^2+y^2-1, z^2+w^2-1, xz+yw\rangle
\]
in the polynomial ring $\RR[x,y,z,w]$. Let $V$ be the space of quadratic forms on two variables $X,Y$. A quadratic form $aX^2+bXY+cY^2$ can be specified by the matrix
\[
Q = \begin{pmatrix}
a & b/2 \\
b/2 & c
\end{pmatrix},
\]
whereupon $G$ acts on quadratic forms via $Q\mapsto PQP^T$. With respect to the basis $X^2, XY, Y^2$, the matrix describing the action of $U$ on $V$ is then
\[
\begin{pmatrix}
x^2&xz&z^2\\
2xy&yz+xw&2zw\\
y^2&yw&w^2
\end{pmatrix}.
\]
The entries of this matrix are the polynomials $a_{ij}$ mentioned above.
\end{example}

Derksen's algorithm for computing a generating set for $\RxG$ proceeds as follows:
\begin{enumerate}
    \item In the polynomial ring $\RR[x_1,\dots,x_n,y_1,\dots,y_n,z_1,\dots,z_\ell]$, choose a monomial ordering that favors the $z$'s over the $x$'s and $y$'s, and with respect to this order, compute a Gr\"obner basis $B$ for the ideal
    \[
    \left\langle h_1,\dots,h_t,\left\{y_i - \sum_{j=1}^n a_{ji}x_j: i=1,\dots, d\right\}\right\rangle.
    \]
    \item Discard any elements of $B$ containing any of the $z$'s.
    \item In the remaining elements, set all $y$'s to zero.
    \item Apply the Reynolds operator to each resulting polynomial. This is a generating set for $\RxG$.
\end{enumerate}
See \cite[Section~4.1]{dk-book} for an illuminating discussion including proof of correctness and several examples.

\subsection{Pontryagin duality}\label{app:pontryagin}

The Pontryagin duality theorem is used in the proof of Theorem~\ref{thm:reg-rep-galois}. It is an abstract generalization of the relationship between the time domain and the frequency domain in Fourier analysis.

\begin{definition}
Suppose $G$ is a locally compact, abelian Hausdorff topological group. A continuous homomorphism from $G$ to the circle group $S^1\defeq \{z:|z|=1\}\subset \CC$ (with its Euclidean topology) is called a {\em character} of $G$. The set of characters forms a group under pointwise multiplication, called the \emph{character group} or the \emph{Pontryagin dual} of $G$, and denoted by $\hat G$.
\end{definition}

For example, a character of $S^1$ itself has the form $z\mapsto z^n$ for $n\in\ZZ$, and the pointwise product of two of these characters adds the exponents; thus $\widehat{S^1}$ is isomorphic to the additive group of integers. In classical Fourier analysis, a periodic function, i.e., a function on $S^1$, is analyzed in terms of its Fourier coefficients, which are indexed by the integers. Thus the Fourier transform of a function on $S^1$ is a function on the character group of $S^1$. This idea generalizes to any locally compact and Hausdorff abelian group.

\begin{remark}
The character group $\hat G$ is itself a locally compact, abelian Hausdorff topological group under the compact-open topology \cite[Theorem~1.2.6]{rudin}.
\end{remark}

Fixing an element $g\in G$, we get a corresponding homomorphism
\begin{align*}
\hat G &\rightarrow S^1\\
\chi&\mapsto \chi(g),
\end{align*}
thus there is a natural map from $G$ to the character group of $\hat G$ (i.e., the double [Pontryagin] dual of $G$) obtained by mapping $g$ to the corresponding homomorphism.

\begin{proposition}[Pontryagin duality]
The map just described is an isomorphism of topological groups.
\end{proposition}

\begin{proof}
See \cite[Theorem~1.7.2]{rudin}.
\end{proof}

\begin{remark}\label{rmk:pontryagin-other-fields}
If $G$ is finite (say of order $N$), the theory empties itself of topological content, so it generalizes to other ground fields than $\CC$. (We mention this in view of Remark~\ref{rmk:nonmodular-field-gen}.) In this case, a continuous homomorphism $G\rightarrow S^1$ is the same as an abstract homomorphism $G\rightarrow\mu_N^\times$, where $\mu_N^\times$ is the group of $N$th roots of unity. Given any field $F$ of characteristic not dividing $N$, the group $\mu_N^\times$ can also be found inside some finite extension $K$ of $F$, and any homomorphism $G\rightarrow K^\times$ must land in $\mu_N^\times$ because every element of $G$ has order dividing $N$. Thus the character group can also be seen as the group of homomorphisms $G\rightarrow K^\times$, and Pontryagin duality holds in this context. 
\end{remark}

\section{Spherical harmonics and $\mathrm{SO}(3)$ invariants}\label{app:so3}

\subsection{Spherical harmonics}

We follow the conventions of \cite{real-spherical}. Parametrize the unit sphere by angular spherical coordinates $(\theta,\phi)$ with $\theta \in [0,\pi]$ and $\phi \in [0,2\pi)$. (Here $\theta = 0$ is the north pole and $\theta = \pi$ is the south pole.) For integers $\ell \ge 0$ and $-\ell \le m \le \ell$, define the complex spherical harmonic
$$Y_{\ell m}(\theta,\phi) = (-1)^m N_{\ell m} P^m_\ell(\cos \theta) \mathrm{e}^{i m \phi}$$
with normalization factor
$$N_{\ell m} = \sqrt{\frac{(2\ell+1)(\ell-m)!}{4\pi (\ell+m)!}}$$
where $P^m_\ell(x)$ are the associated Legendre polynomials
$$P_\ell^m(x) = \frac{1}{2^\ell \ell!} (1-x^2)^{m/2} \frac{\dee^{\ell+m}}{\dee x^{\ell+m}}(x^2-1)^\ell.$$

In the $S^2$ registration problem we are interested in representing a real-valued function on the sphere, in which case we use an expansion (with real coefficients) in terms of the real spherical harmonics:
$$S_{\ell m}(\theta,\phi) = \left\{\begin{array}{ll}
\frac{(-1)^m}{\sqrt 2} (Y_{\ell m}(\theta,\phi) + \overline{Y_{\ell m}}(\theta,\phi)) = \sqrt{2} N_{\ell m} P^m_\ell(\cos \theta) \cos(m \phi) & m > 0, \\
Y_{\ell 0}(\theta,\phi) = N_{\ell 0}P^0_\ell(\cos \theta) & m = 0, \\
\frac{(-1)^m}{i\sqrt{2}} (Y_{\ell |m|}(\theta,\phi) - \overline{Y_{\ell |m|}}(\theta,\phi)) = \sqrt{2}N_{\ell |m|} P_\ell^{|m|}(\cos\theta)\sin(|m|\phi) & m < 0.
\end{array}\right.$$

\noindent Here $\overline{Y_{\ell m}}$ is the complex conjugate of $Y_{\ell m}$, which satisfies the identity 
\begin{equation}\label{eq:conj-identity}
\overline{Y_{\ell m}}(\theta,\phi) = (-1)^m Y_{\ell (-m)}(\theta,\phi).
\end{equation}

\noindent Above we have also used the identity $P^{-m}_\ell = (-1)^m \frac{(\ell-m)!}{(\ell+m)!} P^m_\ell$, which implies $N_{\ell (-m)} P^{-m}_\ell = (-1)^m N_{\ell m} P^m_\ell$.

In the cryo-EM problem we are instead interested in representing the Fourier transform of a real-valued function. Such a function $f$ has the property that if $\vec r$ and $-\vec r$ are antipodal points on the sphere, $f(-\vec r) = \overline{f(\vec r)}$. For this type of function we use an expansion (with real coefficients) in terms of a new basis of spherical harmonics:
$$H_{\ell m}(\theta,\phi) = \left\{\begin{array}{ll}
\frac{1}{\sqrt 2} (Y_{\ell m}(\theta,\phi) + (-1)^{\ell+m}Y_{\ell (-m)}(\theta,\phi)) & m > 0, \\
i^{\ell} Y_{\ell 0}(\theta,\phi) & m = 0, \\
\frac{i}{\sqrt 2} (Y_{\ell |m|}(\theta,\phi) - (-1)^{\ell+m}Y_{\ell (-|m|)}(\theta,\phi)) & m < 0.
\end{array}\right.$$

\noindent One can check that $H_{\ell m}(-\vec r) = \overline{H_{\ell m}(\vec r)}$ using (\ref{eq:conj-identity}) along with the fact $Y_{\ell m}(-\vec r) = (-1)^\ell Y_{\ell m}(\vec r)$ which comes from $P^m_\ell(-x) = (-1)^{\ell+m}P^m_\ell(x)$.

\subsection{Wigner D-matrices}

We will mostly work in the basis of complex spherical harmonics $Y_{\ell m}$ since the formulas are simpler. The analogous results for the other bases can be worked out by applying the appropriate change of basis.

Let $V_\ell \simeq \CC^{2\ell+1}$ be the vector space consisting of degree-$\ell$ complex spherical harmonics represented in the basis $\{Y_{\ell m}\}_{-\ell \le m \le \ell}$, i.e.\ $v \in \CC^{2\ell+1}$ encodes the spherical harmonic $\sum_{m=-\ell}^\ell v_m Y_{\ell m}$. These $V_\ell$ (for $\ell = 0,1,2,\ldots$) are the irreducible representations of $\mathrm{SO}(3)$. Each can also be defined over the real numbers by changing basis to the real spherical harmonics $S_{\ell m}$.

A group element $g \in \mathrm{SO}(3)$ acts on a (spherical harmonic) function $f: S^2 \to \mathbb{R}$ via $(g \cdot f)(x) = f(g^{-1} x)$. The action of $g$ on $V_\ell$ is given by the \emph{Wigner D-matrix} $D^\ell(g) \in \CC^{(2\ell+1) \times (2\ell+1)}$ defined using the conventions of \cite{real-spherical}.

We will need the following orthogonality properties of the Wigner D-matrices. First, the standard Schur orthogonality relations from representation theory yield
$$\Ex_{g \sim \Haar(\mathrm{SO}(3))} \overline{D^\ell_{mk}(g)} D^{\ell'}_{m'k'}(g) = \frac{1}{2\ell+1} \one_{\ell=\ell'} \one_{m = m'} \one_{k = k'}.$$
We also have \cite{rose}
$$D^\ell_{mk}(g) D^{\ell'}_{m'k'}(g) = \sum_{L = |\ell-\ell'|}^{\ell+\ell'} \langle \ell \, m \, \ell' \, m' | L \, (m + m') \rangle \langle \ell \, k \, \ell' \, k' | L \, (k + k') \rangle D^L_{(m+m')(k+k')}(g)$$
where $\langle \ell_1 \, m_1 \, \ell_2 \, m_2 | \ell \, m \rangle$ is a \emph{Clebsch-Gordan coefficient}. There is a closed-form expression for these coefficients \cite{qm-book}:
\begin{align*}
\langle \ell_1 \, m_1 \, \ell_2 \, m_2 | \ell \, m \rangle = &\one_{m = m_1 + m_2} \sqrt{\frac{(2\ell+1)(\ell+\ell_1-\ell_2)!(\ell-\ell_1+\ell_2)!(\ell_1+\ell_2-\ell)!}{(\ell_1+\ell_2+\ell+1)!}} \;\times \\
&\sqrt{(\ell+m)!(\ell-m)!(\ell_1-m_1)!(\ell_1+m_1)!(\ell_2-m_2)!(\ell_2+m_2)!} \;\times \\
&\sum_k \frac{(-1)^k}{k!(\ell_1+\ell_2-\ell-k)!(\ell_1-m_1-k)!(\ell_2+m_2-k)!(\ell-\ell_2+m_1+k)!(\ell-\ell_1-m_2+k)!}
\end{align*}
where the sum is over all $k$ for which the argument of every factorial is nonnegative.

\subsection{Moment tensor}\label{app:so3-moment}

Let $\mathcal{F}$ be a multi-set of frequencies from $\{1,2,\ldots\}$ and consider the action of $G = \mathrm{SO}(3)$ on $V = \oplus_{\ell \in \mathcal{F}} V_\ell$. Recall that we want an explicit formula for $T_d({\bf x}) = \EE_g[(\Pi(g \cdot {\bf x}))^{\otimes d}]$ with $g \sim \Haar(G)$ (where $\Pi$ can be the identity in the case of no projection). We have
$$\EE_g[(\Pi(g \cdot {\bf x}))^{\otimes d}] = \Pi^{\otimes d} \EE_g[\rho(g)^{\otimes d}] {\bf x}^{\otimes d}$$
(where ${\bf x}^{\otimes d}$ is a column vector of length $\dim(V)^d$) and so we need an explicit formula for the matrix $\EE_g[\rho(g)^{\otimes d}]$. Here $\rho(g)$ is block diagonal with blocks $D^\ell(g)$ for $\ell \in \mathcal{F}$. There are no degree-1 invariants since we have excluded the trivial representation $\ell = 0$. For the degree-2 invariants $\EE_g[\rho(g)^{\otimes 2}]$, consider a particular block $\EE_g[D^{\ell_1}(g) \otimes D^{\ell_2}(g)]$ for some pair $(\ell_1,\ell_2)$. The entries in this block can be computed using the above orthogonality relations (and using $D^0_{00}(g) = 1$):
\begin{align*}
\EE_g[(D^{\ell_1}(g))_{m_1 k_1} (D^{\ell_2}(g))_{m_2 k_2}] &= \one_{\ell_1 = \ell_2} \one_{m_1 = -m_2} \one_{k_1 = -k_2} \langle \ell_1 \, m_1 \, \ell_2 \, m_2 | 0 \, 0 \rangle \langle \ell_1 \, k_1 \, \ell_2 \, k_2 | 0 \, 0 \rangle\\
&= \one_{\ell_1 = \ell_2} \one_{m_1 = -m_2} \one_{k_1 = -k_2} \frac{(-1)^{m_1+k_1}}{2 \ell_1+1}
\end{align*}
using the special case $\langle \ell_1 \, m_1 \, \ell_2 \, m_2 | 0 \, 0 \rangle = \one_{\ell_1 = \ell_2} \one_{m_1 = -m_2} \frac{(-1)^{\ell_1+m_1}}{\sqrt{2\ell_1+1}}$.

Similarly, for degree-3 we have
\begin{align*}
\EE_g&[(D^{\ell_1}(g))_{m_1 k_1} (D^{\ell_2}(g))_{m_2 k_2} (D^{\ell_3}(g))_{m_3 k_3}] = \\ &\one_{|\ell_2 - \ell_3| \le \ell_1 \le \ell_2 + \ell_3} \one_{m_1 + m_2 + m_3 = 0} \one_{k_1 + k_2 + k_3 = 0} \frac{(-1)^{m_1 + k_1}}{2\ell_1 + 1} \langle \ell_2 \, m_2 \, \ell_3 \, m_3 | \ell_1 \, (m_2 + m_3) \rangle \langle \ell_2 \, k_2 \, \ell_3 \, k_3 | \ell_1 \, (k_2 + k_3) \rangle.
\end{align*}

\subsection{Projection}\label{app:so3-proj}

Let $V = \oplus_{\ell \in \mathcal{F}} V_\ell$ with $\mathcal{F}$ a subset of $\{1,2,\ldots\}$. Let $\Pi: V \to W$ be the projection that takes a complex spherical harmonic function and reveals only its values on the equator $\theta = \pi/2$. In cryo-EM this projection is applied separately to each shell (see Section~\ref{sec:ex-cryo}). Letting $L = \max_{\ell \in \mathcal{F}} \ell$, the functions $b_{-L}, b_{-L+1}, \ldots, b_L$ (from the circle $S^1$ to $\RR$) form a basis for $W$, where $b_m(\phi) = e^{im\phi}$. The projection $\Pi$ takes the form
$$\Pi(Y_{\ell m}) = (-1)^m N_{\ell m} P^m_\ell(0) b_m$$ extended by linearity. By taking a binomial expansion of $(x^2-1)^\ell$ it can be shown that
\begin{equation}\label{eq:P0}
P_\ell^m(0) = \left\{\begin{array}{ll} 0 & (\ell+m) \text{ odd}, \\ \frac{(-1)^{(\ell-m)/2}}{2^\ell \ell!} \binom{\ell}{(\ell+m)/2} (\ell + m)! & (\ell+m) \text{ even}. \end{array}\right.
\end{equation}

For cryo-EM, if we use the basis $H_{\ell m}$ so that the expansion coefficients are real, the output of the projection can be expressed (with real coefficients) in the basis
$$h_m(\phi) = \left\{\begin{array}{ll} \frac{1}{\sqrt 2}(e^{im\phi} + (-1)^m e^{-im\phi}) & m > 0, \\ 1 & m = 0, \\ \frac{i}{\sqrt 2}(e^{i|m|\phi} - (-1)^m e^{-i|m|\phi}) & m < 0, \end{array}\right.$$
where the projection $\Pi$ takes the form
$$\Pi(H_{\ell m}) = (-1)^m N_{\ell |m|} P^{|m|}_\ell(0) h_m$$
extended by linearity.

\subsection{Explicit construction of invariants}\label{app:cryo-inv}

Consider the cryo-EM setup (see Section~\ref{sec:ex-cryo}) with $S$ shells and $F$ frequencies. We will cover $S^2$ registration as the special case $S = 1$ (without projection). Use the basis of complex spherical harmonics, with corresponding variables $x_{s \ell m}$ with $1 \le s \le S$, $1 \le \ell \le F$, and $-\ell \le m \le \ell$. One can change variables to $S_{\ell m}$ or $H_{\ell m}$ but for our purposes of testing the rank of the Jacobian it is sufficient to just work with $Y_{\ell m}$ (since the change of variables has no effect on the rank of the Jacobian).

Recall that in Appendix~\ref{app:so3-moment} we computed expressions for the matrices $\EE_g[D^{\ell_1}(g) \otimes D^{\ell_2}(g)]$ and $\EE_g[D^{\ell_1}(g) \otimes D^{\ell_2}(g) \otimes D^{\ell_3}(g)]$, and in particular they are rank-1. Using this we can explicitly compute the entries of $T_d({\bf x})$ and thus extract a basis for $U^T_2$ and $U^T_3$. We present the results below.

\subsubsection{Degree-2 invariants} \label{app:deg2invs}

Without projection, the degree-2 invariants are
$$\mathcal{I}_2(s_1,s_2,\ell) = \frac{1}{2\ell+1}\sum_{|k| \le \ell} (-1)^k x_{s_1 \ell k} x_{s_2 \ell (-k)}$$
for $s_1,s_2 \in \{1,\ldots,S\}$ and $\ell \in \{1,\ldots,F\}$. Swapping $s_1$ with $s_2$ results in the same invariant, so take $s_1 \le s_2$ to remove redundancies.

With projection, the degree-2 invariants are
\begin{equation}\label{eq:P2}
\mathcal{P}_2(s_1,s_2,m) = (-1)^m \sum_{\ell \ge |m|} N_{\ell m} N_{\ell (-m)} P^m_\ell(0) P^{-m}_\ell(0) \mathcal{I}_2(s_1,s_2,\ell)
\end{equation}
with $s_1,s_2 \in \{1,\ldots,S\}$ and $m \in \{-F,\ldots,F\}$. Negating $m$ or swapping $s_1$ with $s_2$ results in the same invariant (up to sign) so take $s_1 \le s_2$ and $m \ge 0$ to remove redundancies. Recall the expression (\ref{eq:P0}) for $P^m_\ell(0)$.

It turns out that no degree-2 information is lost in the projection, i.e.\ due to the triangular structure of (\ref{eq:P2}), the $\mathcal{I}_2$ can be recovered from the $\mathcal{P}_2$. To see this, for $m = F,F-1,\ldots,0$ (in that order), use $\mathcal{P}_2(s_1,s_2,m)$ to solve for $\mathcal{I}_2(s_1,s_2,m)$. Note that the coefficient $N_{\ell m} N_{\ell (-m)} P^m_\ell(0) P^{-m}_\ell(0)$ is nonzero iff $m$ and $\ell$ have the same parity.

\subsubsection{Degree-3 invariants} \label{app:deg3invs}

Let $\Delta(\ell_1,\ell_2,\ell_3)$ denote the predicate $|\ell_2 - \ell_3| \le \ell_1 \le \ell_2 + \ell_3$ (which captures whether $\ell_1,\ell_2,\ell_3$ can be the side-lengths of a triangle). Without projection, the degree-3 invariants are
$$\mathcal{I}_3(s_1,\ell_1,s_2,\ell_2,s_3,\ell_3) = \frac{1}{2\ell_1+1} \sum_{\substack{k_1+k_2+k_3 = 0 \\ |k_i| \le \ell_i}} (-1)^{k_1} \langle \ell_2 \, k_2 \, \ell_3 \, k_3 | \ell_1 (-k_1) \rangle x_{s_1 \ell_1 k_1} x_{s_2 \ell_2 k_2} x_{s_3 \ell_3 k_3}$$
for $s_i \in \{1,\ldots,S\}$ and $\ell_i \in \{1,\ldots,F\}$ satisfying $\Delta(\ell_1,\ell_2,\ell_3)$. There are some redundancies here. First, permuting the three $(s_i,\ell_i)$ pairs (while keeping each pair in tact) results in the same invariant (up to scalar multiple). Also, some of the above invariants are actually zero; specifically, this occurs when $(s_1,\ell_1) = (s_2,\ell_2) = (s_3,\ell_3)$ with $\ell_1$ odd, or when $(s_1,\ell_1) = (s_2,\ell_2) \ne (s_3,\ell_3)$ with $\ell_3$ odd (or some permutation of this case).

With projection, the degree-3 invariants are
\begin{align*}
\mathcal{P}_3&(s_1,m_1,s_2,m_2,s_3,m_3) = \\
&(-1)^{m_1} \hspace{-10pt}\sum_{\ell_1,\ell_2,\ell_3\,:\,\Delta(\ell_1,\ell_2,\ell_3)} \hspace{-10pt} N_{\ell_1 m_1} N_{\ell_2 m_2} N_{\ell_3 m_3} P^{m_1}_{\ell_1}(0) P^{m_2}_{\ell_2}(0) P^{m_3}_{\ell_3}(0) \langle \ell_2 \, m_2 \, \ell_3 \, m_3 | \ell_1 (-m_1) \rangle \mathcal{I}_3(s_1,\ell_1,s_2,\ell_2,s_3,\ell_3)
\end{align*}
for $s_i \in \{1,\ldots,S\}$ and $m_i \in \{-F,\ldots,F\}$ such that $m_1 + m_2 + m_3 = 0$. There are again redundancies under permutation: permuting the three $(s_i,m_i)$ pairs results in the same invariant. Negating all three $m$'s also results in the same invariant. There are additional non-trivial linear relations (see Appendix~\ref{app:so3-count} below).

\subsection{Counting the number of invariants}\label{app:so3-count}

\subsubsection{$S^2$ registration}

For the case of $S^2$ registration ($S=1$) the above degree-2 and degree-3 invariants without projection (with redundancies removed as discussed above) form a basis for $\RxG_2 \oplus \RxG_3$ (although we have not made this rigorous). Thus, counting these invariants gives a combinatorial analogue of Proposition~\ref{prop:S2-dim}.

\subsubsection{Cryo-EM}

In this section we give a formula for $\trdeg(U^T_{\le 3})$ for (heterogeneous) cryo-EM (with projection), valid for all $K \ge 1$, $S \ge 1$ and $F \ge 2$. The formula is conjectural but has been tested (via the Jacobian criterion) for various small values of $K, S, F$.

The number of algebraically independent degree-2 invariants turns out to be the number of distinct $\mathcal{I}_2$ invariants (i.e.\ without projection), since the projected invariants $\mathcal{P}_2$ are linear combinations of these. The number of such invariants is $\frac{1}{2} S(S+1)F$.

For degree-3, things are more complicated because the projected invariants $\mathcal{P}_3$ have smaller dimension than the $\mathcal{I}_3$. We start by counting the number of distinct (up to scalar multiple) $\mathcal{P}_3$ invariants. To this end, let $\mathcal{X}(S,F)$ be the set of equivalence classes of tuples $(s_1,m_1,s_2,m_2,s_3,m_3)$ with $s_i \in \{1,\ldots,S\}$ and $m_i \in \{-F,\ldots,F\}$, modulo the relations
$$(s_1,m_1,s_2,m_2,s_3,m_3) \sim (s_2,m_2,s_1,m_1,s_3,m_3) \sim (s_1,m_1,s_3,m_3,s_2,m_2) \qquad \text{(permutation)}$$
$$(s_1,m_1,s_2,m_2,s_3,m_3) \sim (s_1,-m_1,s_2,-m_2,s_3,-m_3) \qquad \text{(negation).}$$
There are some non-trivial linear relations among the distinct $\mathcal{P}_3$ invariants, which we must now account for. The number of such relations is
$$E(S) \defeq 2S + 4S(S-1) + S(S-1)(S-2).$$
This can be broken down as follows. For each $k \in \{1,2,3\}$ there are $2k$ relations for each size-3 multi-subset $\{s_1,s_2,s_3\}$ of $\{1,\ldots,S\}$ with exactly $k$ distinct elements. (We do not currently have a thorough understanding of what exactly the linear relations \emph{are}, but we have observed that the above pattern holds.)

We can now put it all together and state our conjecture. We will also use the formula (\ref{eq:cryo-trdeg}) for $\trdeg(\RxG)$, extended to the heterogeneous case via Proposition~\ref{prop:trdeg-het}.

\begin{conjecture} \label{conj:het-cryo}
Consider heterogeneous cryo-EM with $F \ge 2$ frequencies.
\begin{itemize}
    \item $\trdeg(\RxG) = K[S(F^2 + 2F) - 3] + K-1$,
    \item $\dim(U^T_2) = \frac{1}{2}S(S+1)F$,
    \item $\dim(U^T_3) = |\mathcal{X}(S,F)| - E(S)$,
    \item generic list recovery is possible at degree 3 if and only if $\dim(U^T_2) + \dim(U^T_3) \ge \trdeg(\RxG)$,
    \item generic unique de-mixing holds at degree 3 if and only if $\dim(U^T_2) + \dim(U^T_3) > \trdeg(\RxG)$.
\end{itemize}
\end{conjecture}
\noindent When $S$ and $F$ are large, the dominant term in $\dim(U^T_2) + \dim(U^T_3)$ is $|\mathcal{X}(S,F)| \approx S^3 F^2/4$ and so generic list recovery is possible when (roughly) $K \le S^2/4$.

\begin{remark}
When $S$ is large compared to $F$ we have $\dim(U_2^T) > \trdeg(\RxG)$ and so we might expect generic list recovery to be possible at degree 2. However, this appears to not be the case because unexpected algebraic dependencies are encountered in this regime, i.e.\ $\trdeg(U_2^T) < \trdeg(\RxG) < \dim(U_2^T)$. We have not observed instances where such unexpected algebraic dependencies affect the feasibility of generic list recovery at degree 3.
\end{remark}

\end{document}